\documentclass[12pt,reqno]{amsart}

\usepackage{amsmath}
\usepackage{amsthm}
\usepackage{enumerate}
\usepackage{hyperref}
\usepackage[margin=1.2in]{geometry}
\usepackage[all]{xy}
\usepackage[dvipsnames]{xcolor} 
\usepackage{amssymb} 
\usepackage{graphicx}
\usepackage{subfigure}
\usepackage{etoolbox}

\usepackage[percent]{overpic}

\makeatletter
\def\thm@space@setup{%
  \thm@preskip=12pt plus 4pt minus 6pt
  \thm@postskip=15pt
}
\makeatother

\makeatletter

\def\chaptermark#1{}

\def\chapter{%
  \thispagestyle{headings}
  \@afterindenttrue \secdef\@chapter\@schapter
}

\def\@chapter[#1]#2{\refstepcounter{chapter}%

  \@makechapterhead{#2}\@afterheading}
\def\@schapter#1{\typeout{#1}%
  \let\@secnumber\@empty

  \@makeschapterhead{#1}\@afterheading}
\newcommand\chaptername{\bf Chapter}

\def\@makechapterhead#1{
  \begingroup
  \fontsize{\@xivpt}{16}\bfseries\centering
    \ifnum\c@secnumdepth>\m@ne
      \leavevmode 
            \vskip40pt\centering
      \rlap{\vbox to\z@{\vss\centering
          \centerline{\normalsize\mdseries\centering
              \uppercase\@xp{\chaptername}
              \thechapter}
          \vskip10pt}}
     #1\par \endgroup\vskip 20pt
  }
\def\@makeschapterhead#1{
  \fontsize{\@xivpt}{16}\bfseries\centering
  #1\par \endgroup
 }

\newcounter{chapter}

\newif\if@openright

\makeatother

\numberwithin{section}{chapter}

\theoremstyle{plain}
\newtheorem{theorem}{Theorem}[section]
\newtheorem*{theorem*}{Theorem}
\newtheorem{proposition}[theorem]{Proposition}
\newtheorem{lemma}[theorem]{Lemma}
\newtheorem{corollary}[theorem]{Corollary}
\newtheorem{claim}[theorem]{Claim}

\theoremstyle{definition}
\newtheorem{definition}[theorem]{Definition}

\newtheorem{remark}[theorem]{Remark}

\newtheorem*{fact*}{Fact}
\newtheorem*{claim*}{Claim}
\newtheorem*{remark*}{Remark}

\newtheorem{assumption}[theorem]{Assumption}
\newtheorem*{acknowledgements}{Acknowledgements}

\newtheorem{thm}{Theorem}
\newtheorem{conj}[thm]{Conjecture}
\newtheorem{ques}[thm]{Question}

\newcommand{\NN}{\mathcal{N}}

\newcommand{\ssm}{\smallsetminus}
\newcommand{\bea}{\begin{eqnarray*}}
\newcommand{\eea}{\end{eqnarray*}}

\newcommand{\Int}{\mathrm{int}}

\def\less{<}

\def\int{\text{interior}}

\def\reals {\hbox {\rm {R \kern -2.8ex I}\kern 1.15ex}}
\def\integers {\hbox {\rm { Z \kern -2.8ex Z}\kern 1.15ex}}
\def\naturals {\hbox {\rm {N \kern -2.8ex I}\kern 1.20ex}}
\def\rationals {\hbox {\rm { Q \kern -2.2ex l}\kern 1.15ex}}
\def\hyp {\hbox {\rm {H \kern -2.7ex I}\kern 1.25ex}}

\makeatletter
\patchcmd{\chapter}{\if@openright\cleardoublepage\else\clearpage\fi}{}{}{}
\makeatother

\title {Tunnel number one knots satisfy the Berge Conjecture}

\author{Tao Li}

\author{Yoav Moriah}

\author{Tali Pinsky}

\thanks{The first author is partially supported by an NSF grant. The second author 
wishes to thank the University of British Columbia UBC for its generous hospitality}

\date{\today}

\subjclass[2010]{Primary 57M99}

\keywords{Berge Conjecture, doubly primitive, Heegaard splittings, lens spaces,  
genus reducing surgery}

\address{Department of Mathematics \\
 Boston College \\
 Chestnut Hill, MA 02467}
\email{taoli@bc.edu}

\address{Department of Mathematics \\
Technion \\
Haifa, 32000 Israel}
\email{ymoriah@tx.technion.ac.il} 

\address{Department of Mathematics \\
Technion \\
Haifa, 32000 Israel}
\email{talipi@technion.ac.il}

\begin{document}

\begin{abstract} Let $K$ be a tunnel number one knot in $M$ with irreducible knot exterior, where $M$ is either $S^3$, 
 or a connected sum of $S^2\times S^1$ with any lens space. (In particular, this includes 
 $M = S^2\times S^1$.)  We prove that if a non-trivial Dehn surgery on $K$ yields a lens space, then $K$ is a doubly primitive knot in $M$.  For  $M  = S^3$ this resolves the 
 tunnel number one Berge Conjecture.  For $M  = S^2\times S^1$  this resolves a conjecture of 
 Greene and Baker-Buck-Lecuona for tunnel number one knots.  
 \end{abstract}

\maketitle
\section* {\bf Introduction}

One of the main goals of low dimensional topology is to determine which manifolds are obtained 
by  Dehn fillings on knots in the $3$-sphere.  D. Gabai proved, in \cite {Ga},  the Property $R$ Conjecture,
 i.e. if a Dehn surgery on a knot in $S^3$ yields  $S^2 \times S^1$, then the knot must be a trivial knot. 
Another major achievement in the 80's was the proof by C. Gordon  and J. Luecke, in \cite{GL},  that if 
non-trivial  surgery on a knot yields $S^3$ then the knot  must be the trivial knot.  Both $S^2 \times S^1$ 
and $S^3$ are  lens spaces.  So a natural question is:

\begin{ques} Which knots in $S^3$ other than the unknot have non-trivial surgery 
resulting in a  lens space which is not $S^3$ or $S^2 \times S^1$?
\end{ques}

This question was first raised by L. Moser in 1971 in \cite{Mos}.  She determined the 
surgeries on torus knots yielding lens spaces. Subsequently, Bleiler  and Litherland \cite{BL}, 
Wang \cite{Wa} and Wu \cite{Wu},  independently, characterised  the surgeries on satellite knots 
in $S^3$  which  result in lens spaces. The question is still  unresolved for hyperbolic  knots.

In his celebrated notes, see \cite{Th},  W. Thurston proved that each hyperbolic knot has  only finitely 
many surgery slopes so that the manifolds obtained by Dehn fillings along these slopes fail to be hyperbolic 
(this number is uniformly bounded over all hyperbolic knots, see \cite{A,BlHo,L,LM}).   There has been a  huge effort by many mathematicians to precisely determine these fillings.
There was a huge effort by many mathematicians to precisely determine these fillings. As lens spaces are 
non-hyperbolic, 
the above question can be viewed  as a special case of this huge task.

John Berge observed, see \cite{Be1}, that if a knot $K$ is doubly primitive \footnote{Berge calls  
these knots  ``double primitive''.} (see definitions below), then a Dehn surgery on $K$ yields a 
lens space.  He compiled a list of twelve families of doubly primitive knots $K \subset S^3$ 
including known cases of torus and satellite knots and  asked whether a knot $K\subset S^3$ 
is on this list if and only if $K$ is doubly primitive?  (This question was later formally conjectured 
by C. Gordon  (see \cite {Ki} Problem  1.78)).
\footnote{Ultimately the proof that Berge's list for knots in $S^3$ is complete is a consequence of J. Greene's
 Lens space realization paper (see \cite{Gr}).}  The question became to be known as  
{\it The Berge Conjecture}:

\begin{conj}[The Berge Conjecture] \label{con:BergeConjecture}
Let $K \subset S^3$ be a non-trivial knot which has Dehn surgery resulting in a lens space. 
Then $K$ is doubly primitive.
\end{conj}

A simple closed curve on the boundary of a genus two handlebody is {\it primitive} if the handlebody has a 
compressing disk that transversely intersects this curve in a single point.  A knot $K \subset M$ is {\it doubly primitive} if 
$M$ has a genus two Heegaard splitting and $K$ can be isotoped onto the Heegaard surface so that $K$ 
is primitive in both handlebodies of the Heegaard splitting.  Note that if one adds a 2-handle to a genus two 
handlebody along a primitive curve, then the resulting manifold is a solid torus. Thus a doubly primitive knot
 always has a Dehn surgery that results in a lens space.  
 
 There is a stronger version of the Berge Conjecture,  which says that if a Dehn surgery on $K$ produces a lens space, then $K$ has a doubly primitive presentation such that the framing given by the Heegaard surface is the surgery slope. Note that a knot may have more than one doubly primitive presentation.  By \cite[Theorem 10]{BDH}, non-torus doubly primitive knots in $S^3$ or $S^2\times S^1$ have no unexpected lens space surgery.

Given a knot $K$ in a closed orientable 3--manifold $M$, a {\it tunnel system} for $K$ is a collection of disjoint 
properly embedded arcs $\{t_1,\dots, t_n\}$ in the knot exterior $E(K)=M\ssm \NN(K)$, where $\NN(\cdot)$ 
means a regular neighborhood of $\{\cdot\}$, such that $E(K)\ssm \NN(\cup_{i=1}^nt_i)$ is a handlebody.  The 
{\it tunnel number} of the knot $K$ is the minimal number of arcs in a tunnel system for $K$.  Note that a knot 
$K$ has tunnel number $n$ if and only if $M\ssm \NN(K)$ has  Heegaard genus $n+1$.

Let $K \subset M$ be a doubly primitive knot. If one pushes $K$ into the interior of either genus two handlebody 
of the corresponding Heegaard splitting, then after removing $\NN(K)$, this handlebody becomes a compression 
body and the genus two Heegaard surface for $M$ becomes a Heegaard surface of $M\ssm \NN(K)$. Moreover, 
the 1-handle in the compression body determines an unknotting tunnel for $K$.  This means that all doubly 
primitive knots have tunnel number one.

Therefore, the Berge Conjecture can be divided into two parts:

\begin{conj}[The tunnel number one Berge Conjecture]\label{con:T=1BG}
If $K\subset S^3$ is a tunnel number one knot which admits a Dehn surgery resulting 
in a lens space, then $K$ is doubly primitive.
\end{conj}

\begin{conj}[The lens space Dehn surgery]\label{con:LSDF}
If $K\subset S^3$ is a nontrivial knot which admits a Dehn surgery resulting in a lens space, then $K$ 
is a tunnel number one knot.
\end{conj}

As the notion of doubly primitive knots is not limited to $S^3$, one can ask whether a conjecture similar to the  
Berge Conjecture  holds for other manifolds with Heegaard genus at most $2$:

\begin{ques}\label{ques:general}
Let $K$ be a knot in a 3-manifold $M$ such that a Dehn surgery on $K$ yields a lens space.  Is $K$ necessarily 
doubly primitive with respect to a genus two Heegaard splitting of $M$? 
\end{ques}

The answer to Question~\ref{ques:general} is known to be negative when $M$ is the Poincar\'{e} homology
 sphere~\cite{BH}  or a lens spaces in general \cite{BBL}.   Namely, there are knots in 
 the Poincar\'{e} homology sphere and in non trivial lens spaces so that surgery on them yield lens spaces
 but that the knots are not doubly primitive with respect to a genus two Heegaard splitting of these spaces,

 The next conjecture, which was made by J. Greene  (see  ~\cite[Conjecture 1.8]{Gr}) and by Baker, Buck 
 and Lecuona (see \cite[Conjecture 1.1]{BBL}), says that the answer to Question~\ref{ques:general} is 
 expected to be true if $M$ is $S^2\times S^1$.

\begin{conj}[Berge Conjecture for $S^2 \times S^1$]
\label{con:GreeneConjecture}
If $K$  is a  knot in $S^2 \times S^1$ which admits a lens space Dehn surgery, then $K$ is doubly primitive.
\end{conj}

In Theorem \ref{thm:MainTheorem} we prove Conjecture~\ref{con:T=1BG}. Theorem
~\ref{thm:MainTheorem} also implies that Conjecture~\ref{con:GreeneConjecture} also holds for knots with 
tunnel number one.

\begin{thm}\label{thm:MainTheorem} Let $K\subset M$ be a tunnel number one knot with irreducible knot exterior,
where $M$ is either $S^3$ or $ (S^2 \times S^1) \# L(r,s)$, (where $L(r,s)$ is 
any lens space).  If a non-trivial Dehn surgery on $K$ yields a lens space, then  $K$  
is doubly primitive.
\end{thm}

Note that if $M=(S^2 \times S^1) \# L(r,s)$ and $K$ is a core curve of a Heegaard solid torus in the lens space summand $L(r,s)$, then a Dehn surgery on $K$ may change $L(r,s)$ to $S^3$, thereby transforming $M$ to the lens space $S^2\times S^1$. However, $K$ may not be doubly primitive on the genus-2 Heegaard surface of $M$.  
This special case is the only reason we require the knot exterior $M\setminus K$ to be irreducible in Theorem \ref{thm:MainTheorem}. 
In particular, this hypothesis has no impact when $M$ is $S^3$ or $S^2\times S^1$, since $S^3\setminus K$ is irreducible for all knots, and if $(S^2\times S^1)\setminus K$ is reducible, then $K$ lies in a 3-ball and no nontrivial Dehn surgery produces a lens space. 
The only point in our proof where this hypothesis is required is in Lemma~\ref{lem:integer slope}.

Theorem \ref{thm:MainTheorem}  suggests that the answer to Question~\ref{ques:general} might be 
positive in the case where $M$ is a connected sum of $S^2\times S^1$ and a lens space.

Let $K$ be a knot in $S^3$. We say that $K$ admits {\it genus reducing surgery} if  
$S^3 \ssm \NN(K)$  has a non-trivial surgery resulting in a manifold $M$ so that 
 $g(M) \leq g(S^3 \ssm \NN(K)) - 1$, where $g(X)$ denotes the Heegaard genus 
of the manifold $X$.  It will be called {\it strongly genus reducing surgery} if 
$g(M) \less g(S^3 \ssm \NN(K)) - 1$. In this context, the second author would like to 
make the following conjecture:

\begin{conj}[Strong genus reducing surgery]\label{con:StrongGenusReducingConjecture}
Knots in $S^3$ do not have strongly genus reducing surgery slopes other than the slope of the meridian.
\end{conj}

The following is a weaker conjecture: 

\begin{conj}\label{con:WeakConjecture}   For any knot $K \subset S^3$, no integer slope  is a strongly 
genus reducing  surgery slope.
\end{conj}

It is a consequence of the Cyclic Surgery Theorem \cite{CGLS} that if $K \subset S^3$ is a hyperbolic 
knot only integer slope surgery can yield a lens space.  Thus Conjecture~\ref{con:WeakConjecture}, 
if true,  together with Theorem \ref{thm:MainTheorem} implies  the Berge Conjecture. 

The Berge Conjecture has been studied extensively, e.g., by Berge himself \cite{Be1,Be2, Be3} also by 
D. Gabai \cite{Gabai}, K. Baker,  E. Grigsby, and M. Hedden \cite{BGH}, P. Ozsvath and Z. Szabo \cite{OZ}, 
J. Rasmussen  \cite{RJ},  T. Saito \cite{TS}, M. Tange \cite{Ta1,Ta2}, Y. Ni \cite{Ni}, J. Greene \cite{Gr} and
others. There is some strong evidence for the Berge Conjecture. For example, Yi Ni \cite{Ni} proved that 
if a knot in $S^3$ admits lens space surgery, then the knot must be fibered (all Berge knots are fibered).  
Greene proved that if a lens space can be  obtained by surgery on a knot in $S^3$ then such a lens space  
can be  obtained by a knot in the Berge list (see \cite[Theorem 1.3] {Gr}).

However, very little is currently known on either Conjecture \ref{con:StrongGenusReducingConjecture} or 
Conjecture \ref{con:WeakConjecture}.  Some results in this direction were obtained by K. Baker, 
C. Gordon and J. Luecke, see  \cite[Corollary 1.1] {BGL}.

\vskip5pt

\noindent{\bf The main idea of the proof and outline of the paper:}

\vskip3pt

Assume that  $M$ is a $3$-manifold and $K \subset M$ is a knot. Let $\{U, W\}$ be a Heegaard splitting for 
$M \ssm \NN(K)$ along a genus-$g$ Heegaard surface $\Sigma$, where $U$ is a compression body and 
$W$ is a handlebody, see Definitions \ref{def:CompressionBody}. Then $\Sigma$ is also a Heegaard surface 
of the manifold $M(r)$ obtained by $r$-surgery on $M \ssm \NN(K)$. If the Heegaard splitting of $M(r)$ along 
$\Sigma$ is stabilized, then the handlebody $W$ must contain an essential disk $D$ and the compression body 
$U$ must contain a planar surface $P$ so that $P \cap D$ is a single point in $\Sigma$ and, after the Dehn 
surgery on $K$, $P$ becomes an essential disk $\widehat{P}$. So $(\widehat{P}, D)$ is a genus-reducing 
destabilizing disk pair, see Definitions \ref{def:CompressionBody}  and  \ref{def:PD} below. The pair $(P,D)$ is 
called a $(\mathcal{P},\mathcal{D})$-pair.

As described at the beginning of Section \ref{sec:complexity}, we are dealing with a genus two Heegaard 
splitting of $M \ssm \NN(K)$, where $M$ is either $S^3$ or $(S^2 \times S^1) \# L(r, s)$, which becomes,  
after integer surgery, a reducible genus-two Heegaard splitting of a lens space. Thus the genus-two Heegaard 
splitting must contain a $(\mathcal{P},\mathcal{D})$-pair $(P,D)$. 

Consider a Heegaard diagram of the genus-two Heegaard splitting of $M$, namely 
$\widehat{V} = \{\alpha, \gamma\}$ and $\widehat{W} = \{\delta, \varepsilon\}$, where  $\delta=\partial D$, 
$\varepsilon$ is a meridian of the handlebody $W$ disjoint from $\delta$ and $\partial P$, $\alpha$ is a 
boundary curve of a vertical annulus in the compression body $U$, and $\gamma$ is the boundary of a 
disk in $U$. See the beginning of Section \ref{sec:MeridionalDisksandWhiteheadGraphs}   
Definition \ref{def:WhiteheadGraphs} and  Definition \ref{def:standard} for more details. 

The intuition behind the proof is the following:  We do not know what the curves  $\alpha, \gamma, \delta, \varepsilon$ 
and $\partial P$ look like. However, since the Heegaard surface $\Sigma$ is of genus two and since 
$ \gamma \cap \alpha = \delta \cap \varepsilon = \varepsilon\cap \partial P = \emptyset $ and 
$\partial P  \cap \partial D$  is a single point, the possible configurations of these curves should be 
sufficiently  limited to give the required result. 

\vskip5pt

\noindent {\bf Main idea:} The name of the game is to analyse all possible configurations for 
these curves. Note that if we cut the Heegaard surface $\Sigma$ along $\partial_+ P$, the curves  
$\widehat{V} = \{ \gamma, \alpha\}$ and $\widehat{W}= \{\delta, \varepsilon \}$  
are cut into a collection of arcs on $\Sigma \ssm \NN(\partial_+ P)$. Using this information  
define a notion  of \emph {complexity}  $c(P,D,\alpha,\gamma)$ of a $(\mathcal{P},\mathcal{D})$-pair 
(see Definition \ref{def:complexity}).  In order  to reduce this complexity, we perform  \emph{wave moves}
 (see Definition \ref {def:Shortcuts}) on the above curves which yield new curves and meridional disks
so that either: 
\begin{enumerate}
\item We find a new $(\mathcal{P},\mathcal{D})$-pair or a new $\alpha$ with smaller complexity 
$c(P,D,\alpha,\gamma)$.  If this is the case then, after iterating the argument if need be, 
we use Proposition~\ref{pro:Tao} to conclude that $K$ is  a Berge knot, or

\item The knot $K$ is a Berge-Gabai knot, see \cite{Be2, Gabai}, and hence doubly primitive i.e., 
$K$ is a Berge knot, or

\item The new Heegaard diagram, which is not standard, has no wave. In this case 
Theorems~\ref{thm:HOT} and \ref{thm:NeOk} are contradicted when $M$ is one 
of $S^3$ or $ (S^2 \times S^1) \# L(r,s)$, where $L(r, s)$ is a lens space. 
\end{enumerate}

\vskip5pt

\noindent {\bf Outline of the paper:} In Chapter \ref{cpt:preliminaries}, Sections 1 - 3 we define the necessary 
notions and prove some of the lemmas that are needed throughout the paper. In particular, we define a notion 
of \emph{complexity} and show that if there is a $(\mathcal{P},\mathcal{D})$-pair and an annulus with 
very small complexity then the knot $K$ must be doubly primitive, see Proposition~\ref{pro:Tao}. We also 
study the Whitehead graphs of genus-two Heegaard splittings of $M$. The edges of a Whitehead graph 
of $M$ are restricted by Theorems~\ref{thm:HOT}  and \ref{thm:NeOk}. For example if the Heegaard 
diagram is not ``standard'' (as in Definition \ref{def:standard}) the edges  cannot be \emph{blocking-edges} 
as defined  in  Definition \ref{def: blockingedges}.

The next crucial step, done in  Section \ref{sec:NoCircularArcs} of Chapter \ref{cpt:PlanarSurfacesGenusTwo}, 
is to first cut the Heegaard surface $\Sigma$ open along the curve $\partial_+P$, denoting the resulting
compact surface by $\widehat{\Sigma}$, and then fix a ``standard'' decomposition for $\widehat{\Sigma}$.  
The surface $\widehat{\Sigma}$ is composed of two annuli, denoted by $\mathcal{A}_r$ and $\mathcal{A}_l$, 
and two rectangles denoted by $\mathcal{R}^u$ and $\mathcal{R}^d$. The arcs $ \gamma \ssm \partial_+P $, 
$\alpha \ssm \partial_+P$, $ \delta \ssm \partial_+P ,$ and the curve  $\varepsilon$ in 
$\widehat{\Sigma} = \mathcal{A}_r \cup \mathcal{R}^u \cup \mathcal{R}^d \cup \mathcal{A}_l$  
have certain properties and restrictions, which are discussed in lemmas,  remarks  and propositions labeled 2.1.1 - 2.1.9.

In Section \ref{sec:Delta} we concentrate on the $\{\varepsilon, \delta\}$  and $\partial_+P$ curves.  
In order to describe the curves we define the notion of \emph{paths}, more specifically  \emph{short} and  
\emph{long} paths, which the curves $\{\varepsilon, \delta\}$ must follow  (take) (see  Definition \ref{def:Path}). 
Each of the curves can take a path many times. Since the curves $\{\varepsilon, \delta\}$  after being cut 
along $\partial_+P$ can have many parallel arcs which take the same path, it is convenient to think  of them 
as  \emph{train tracks}. This is done  in the discussion  below Definition \ref{def:Path}. These  train tracks will 
play a crucial role in Chapter \ref{cpt:ObtainingTheContradiction}. 

Chapter \ref{cpt:ObtainingTheContradiction}, which is the essence of the proof, is a careful and 
detailed analysis of the curve $\delta $  in $\widehat{\Sigma}$ depending on  the paths that 
$\delta$ takes. It turns out that there are four 
possibilities, i.e.,

\begin{enumerate}
\item $\delta$ takes one short path in $\mathcal{A}_l$ and one short path in $\mathcal{A}_r$.
\vskip3pt
\item $\delta$ takes one short path in $\mathcal{A}_l$ or $\mathcal{A}_r$ and two short paths 
in the other annulus.
\vskip3pt
\item $\delta$ takes two short paths in $\mathcal{A}_l$ and two short paths in $\mathcal{A}_r$.
\vskip3pt
\item $\delta$ takes a long path in one (or both) of the annuli.
\end{enumerate}
This analysis gives the three possibilities as in the main idea above, which proves the theorem.

\begin{acknowledgements} 
The authors wish to thank John Berge for important remarks regarding $(1,1)$-knots and for pointing out the 
importance of the papers by T. Homa, M. Ochiai and M. Takahashi \cite{HOT} and M. Ochiai \cite{Och} to the 
Berge Conjecture. The first author would like to thank Ken Baker and Alex Zupan for helpful communications. 
The second and third authors wish to thank Dale Rolfsen for his support. Also we wish to thank the referee 
for numerous corrections and suggestions.
\end{acknowledgements}

 
 
\chapter{Preliminaries}\label{cpt:preliminaries}

Throughout the paper, we assume that $M$ is either $S^3$ or $ (S^2 \times S^1) \# L(r, s)$, where $L(r, s)$ 
is a  lens space. Since $S^3$ is a lens space the case $M =  S^2 \times S^1$ is included. For any space $Y$, 
we use $|Y|$ to denote the number of components of $Y$,  use $\Int(Y)$ to denote the interior of $Y$ and use 
$\NN(Y)$ to denote a small regular neighborhood of $Y$. Finally, throughout the paper the hyperelliptic involution 
will be denoted by $\pi$.

\section{Primitive and doubly primitive}\label{sec:Primitive}

We begin with some definitions:

\begin{definition}\label{def:CompressionBody}
 A {\it compression body} is a connected $3$-manifold $U$ obtained by adding 2-handles to a product 
 $\Sigma\times I$ ($I=[0,1]$) along $\Sigma\times\{0\}$, where $\Sigma$ is a closed and orientable surface, and 
 then capping off resulting 2-sphere boundary components by 3-balls.  The surface $\Sigma\times\{1\}$ 
 is denoted by $\partial_+U$ and $\partial U \ssm \partial_{+}U$ is denoted by $\partial_{-}U$.  If 
 $\partial_{-}U = \emptyset $ then $U$ is a handlebody.   One can also view a compression body 
 with $\partial_-U\ne\emptyset$ as a manifold obtained by adding 1-handles on the same side of 
 $\partial_-U\times I$.

 A {\it Heegaard splitting} of a $3$-manifold $Y$, denoted by $(U,W)$, is a decomposition of $Y$ into two 
 compression bodies $U$ and $W$ along a closed surface $\Sigma=\partial_+U=\partial_+W$, and $\Sigma$ 
 is called a {\it Heegaard surface}.
 \end{definition}

Let $K \subset M$  be  a knot in $M$.  Throughout the paper, we will use  $(U,W)$ to denote a Heegaard 
splitting for $M \ssm \NN(K)$ with $\partial_{-} U = \partial (M \ssm \NN(K))$ being the boundary torus, and 
we use $\Sigma$ to denote the Heegaard surface.  Note that a trivial Dehn filling along $\partial_-U$ extends 
$U$ to a handlebody $V$, and hence $\Sigma$ is also a Heegaard surface of $M$. Throughout the paper, 
we use $(V,W)$ to denote this Heegaard splitting of $M$, with $\Sigma=\partial V=\partial W$, $K\subset V$ 
and $U=V \ssm \NN(K)$.

For any surface $S$ properly embedded in a compression body $U$, we use $\partial_{-}S$ to denote 
$S \cap \partial_{-}U $  and   $ \partial_{+}S$ to denote  $S \cap \partial_{+}U  \subset  \Sigma$. 

Throughout this paper $P$ will denote a planar surface which is not a disk in $U$. We require further that 
$P$  has a single boundary component, denoted by $\partial_{+}P$ on $\partial_{+}U$  and  $N \geq 1$ 
boundary  components on  $\partial_{-} U$.  

An annulus $A$ properly embedded in the compression body $U$ is called a {\it vertical annulus} if $A$ 
is incompressible and has one boundary component in  $\partial_{+}U$ and the other boundary component 
in $\partial_{-}U$.  Note that there are many vertical annuli in $U$ even with the same $\partial_-A$ curve, 
by taking band sums of one such $A$ with  compressing disks in $U$.

 \begin{definition}\label{def:PD} 
 Let $(U,W)$ be a Heegaard splitting of $M\ssm\NN(K)$ as above.  Let $P$ be a planar surface 
 properly embedded in $U$ and $D$ a compressing disk in $W$.  Suppose $\partial_+P$ is a 
 single curve transversely intersecting $\partial D$ in $\Sigma$.  We say that $(P,D)$ 
 is a {\it $(\mathcal{P},\mathcal{D})$-pair} with respect to the Heegaard splitting with slope 
 $r$ if $\partial D\cap\partial_+P$ is a single point and $\partial_-P$ consists of essential curves 
 of slope $r$ in the boundary torus.  If $P$ is an annulus then we say that $(P,D)$ is an  
 {\it $(\mathcal{A}, \mathcal{D})$-pair}.
\end{definition}  

In this paper we assume that some  nontrivial Dehn surgery on $K$ yields a lens space. 

\begin{lemma}\label{lem:integer slope}
Let $K$ be a knot in $M$ with irreducible knot exterior, where $M$ is either $S^3$ or $(S^2 \times S^1) \# L(r,s)$. Suppose $K$ 
admits a non-trivial Dehn surgery resulting in a lens space. Then either $K$ is doubly primitive, or the surgery slope is an integer slope.
\end{lemma}

\begin{proof}
If  $M = S^3$,  then since $K$ admits lens space surgery, by the Cyclic Surgery Theorem (see \cite{CGLS}), 
either the surgery slope $r$ is an integer or $K$ is a torus knot. Since a torus knot is doubly primitive 
(see ~\cite{Be1}), the lemma holds. 

Suppose therefore that $M$ is  $S^2\times S^1$ or  $(S^2 \times S^1) \# L(r,s)$. 
Since $M \ssm \NN(K)$ is irreducible, by \cite[Corollary 1.4]{BoZa}, either the 
surgery slope is an integer slope or $M\ssm \NN(K)$ is a simple Seifert fibered space.  If $M\ssm \NN(K)$ 
is a simple Seifert fibered space, then it contains an essential annulus $A$ which divides $M\ssm \NN(K)$ 
into two solid tori. Since $M$ is  $(S^2 \times S^1) \# L(r, s)$ and since $M \ssm \NN(K)$ is irreducible, 
$M\ssm \NN(K)$ admits an essential planar surface with $\infty$-slope, i.e.~$S^2  \ssm  \NN(K)$, where 
$S^2$ is an essential $2$-sphere in  $(S^2 \times S^1) \# L(r, s)$.  By \cite{GL2}, the boundary slope of 
$A$ and the boundary slope of the planar surface $S^2  \ssm  \NN(K)$ has geometric intersection one.  Hence 
the boundary slope of $A$ is an integer slope.  Since $A$ divides $M\ssm \NN(K)$ 
into two solid tori and $\partial A$ has integer slope, $K$ must lie on a torus which divides $M$ into two solid tori, which implies that $M=S^2\times S^1$.  Similar to the argument about torus knots (see ~\cite{Be1}), $K$ is doubly 
primitive with respect to a genus two  Heegaard surface of $M$.
\end{proof}

\begin{proposition} \label{pro:DPEquivalentToAD} 
Let $K\subset M$ be a tunnel number one knot with irreducible knot exterior, where 
$M$ is $S^3$ or $(S^2 \times S^1) \# L(r,s)$. Suppose $K$ admits a lens space surgery,  Then, if there is an 
$(\mathcal{A},\mathcal{D})$-pair with respect to a genus two Heegaard splitting of $M \ssm\NN(K)$ 
 where $\partial_{-} A$ has the surgery slope, then $K$ is doubly primitive.
\end{proposition}

\begin{proof} 
Suppose  $(A, D)$ is such a pair, i.e. $A \subset U$ is a vertical annulus meeting  the 
disk $D \subset W$ in a single point in $\partial_{+}A$. If we perform Dehn filling  along 
the slope $r$ of $\partial_{-}A$, $A$ extends to a disk and the $(A,D)$-pair extends to a 
stabilizing pair of disks in the manifold $M(r)$  obtained by Dehn filling along the slope $r$. 
Since the Heegaard splitting has genus two, this implies that $M(r)$ has Heegaard genus 
at most one and hence is a lens space.  

As in Lemma~\ref{lem:integer slope}, this means that the slope $r$ of $\partial_-A$ is an integer or $K$ 
must be doubly primitive.  Suppose  $r$ is an integer, then $\partial_-A$ and hence $\partial_+A$ are 
isotopic to the knot $K$ in $M$. So we may isotope $K$ to the curve $\partial_+A\subset\Sigma$.  
Since $\partial_+A$ transversely intersects $\partial D$ in one point, $K$ is primitive in $W$. Since 
$U$ is a compression body, $U$ contains a vertical annulus $A'$ such that $\partial_-A'$ has
 meridional slope.  Hence $\partial_+A'$ bounds a disk $D'$ in the handlebody $V=U \cup \NN(K)$ 
 intersecting $K$ in one point. So $K$ is primitive in $V$ as well. 
\end{proof} 

\begin{remark}\label{rem:NotAllSurfaces} 
A tunnel number one knot may have multiple different unknotting tunnels. Thus, 
a priori, it is possible that a knot $K$ is not doubly primitive with respect to 
some genus two Heegaard surface of $S^3 \ssm \NN(K)$ but  is doubly primitive 
with respect to a different genus two Heegaard surface of $S^3 \ssm \NN(K)$.  

In fact, the $(-2, 3, 7)$-pretzel knot $K$ is a tunnel number one knot with four unknotting tunnels  
(see \cite{Heath-Song}). John Berge pointed out that the compression body in one  of these 
genus two Heegaard  splittings can be obtained by adding a $2$-handle to the boundary of a genus 
 two handlebody $H$  along a simple closed curve $\rho \subset \partial H$ which represents the
 word $x^3 y^2 x^2 y^2 x^3 y^3$ in  $\pi_1(H) = F(x,y)$, the free group on two generators. In this case, 
 there is only one simple closed  curve disjoint from $\rho$ which represents  a primitive in $H$, namely 
 the curve $x y$, which is the  meridian of $H(\rho) = S^3 \ssm \NN(K)$.  So clearly $K$ cannot be doubly 
 primitive with respect to this Heegaard splitting. However the $(-2, 3, 7)$-pretzel knot is on the Berge list. 
 
 This example implies that the proof of the main theorem, Theorem \ref{thm:MainTheorem}, will have to
 deal with this possible situation. This will be done in the next section.
\end{remark}


\section{A notion of complexity for a $(\mathcal{P},\mathcal{D})$-pair}\label{sec:complexity}

Let $K \subset M$ be a tunnel number one knot and $(U, W)$ a genus two Heegaard splitting of 
$M \ssm \NN(K)$, where $W$ is a handlebody, and let $\Sigma$ be the Heegaard surface.  
The surface $\Sigma$ also determines a Heegaard splitting $(V,W)$ of $M$ with $V=U\cup\NN(K)$. 
Suppose the Dehn filling along a slope $r$ yields a lens space.   Note that $\Sigma$ is 
also a genus two Heegaard surface of the lens space.  By Bonahon-Otal ~\cite{BO}, the genus 
two Heegaard splitting of the lens space is stabilized by  a stabilizing pair of disks $(D, E)$. 
Assume $D \subset W$ and $E \subset V$. Hence  $E \ssm \NN(K) \subset V \ssm \NN(K) = U$   is  a 
planar surface $P \subset U$ which intersects $D\subset W$ in a single point. Thus, stabilizing pair of 
disks gives rise to a $(\mathcal{P},\mathcal{D})$-pair $(P,D)$, where $\partial_-P$ has slope $r$. 
By Lemma~\ref{lem:integer slope}, we may suppose that $r$ is an integer slope.  If $P$ is compressible in $U$, then by compressing $P$, we obtain a new planar surface $P'$ with $\partial_+P'=\partial P_+$. 
So, after finitely many compressions, we obtain a new $(\mathcal{P},\mathcal{D})$-pair with incompressible planar surface. 
Assume, throughout the 
paper that  $P$ is incompressible in $U$.  By  Proposition~\ref{pro:DPEquivalentToAD}, we shall assume $P$ is not an annulus (as otherwise 
we are done). 

We now proceed to define a notion of complexity for $(\mathcal{P},\mathcal{D})$-pairs. As the 
Heegaard splitting is of genus two, there is a unique, up to isotopy, nonseparating disk in  $U$ which is 
denoted by $C$. There is a vertical annulus $A$ in $U$ such that $\partial_-A$ has the
$\infty$-slope (i.e.~the meridional slope of $K$ in $M$).  Note that $A$ is not unique, in  fact, one can construct infinitely many such annuli 
by band summing $A$ with $C$ in different ways.  Suppose $C\cap A=\emptyset$. 
Let $\gamma=\partial C$ and $\alpha=\partial_+A$.  Note that $\gamma$ and $\alpha$ lie on the 
Heegaard surface $\Sigma$, and  both $\gamma$ and $\alpha$ bound nonseparating disks in the 
handlebody $V=U\cup \NN(K)$.

 \begin{definition}\label{def:complexity} 
Given a $(\mathcal{P},\mathcal{D})$-pair $(P, D)$  in a Heegaard splitting of $M\ssm\NN(K)$ and 
a collection of disjoint simple closed curves $\{\alpha_1, \dots, \alpha_n\} \subset \Sigma$ in general 
position with respect to $(\partial P_{+}, \partial D)$, define the complexity $c_0(P, D, \alpha_1, \dots, \alpha_n)$ 
to be the number of segments of $\{\alpha_1 \cup \dots \cup \alpha_n\} \ssm \partial_{+} P$ which intersect $D$.

In the current situation, for a given $(P,D)$, consider the complexity $c_0(P, D, \alpha, \gamma)$, where 
$\gamma = \partial C$ and $\alpha = \partial_+ A$ and where $A$ is a vertical annulus in $U$ such that 
$\partial_-A$ has the $\infty$-slope.  Further, define $c_1(D, \alpha, \gamma)$ to be the number of 
intersection points of $\partial D$ and  $\alpha \cup \gamma$.  Now, define  the pair
$$c(P, D, \alpha, \gamma) = (c_0(P, D, \alpha, \gamma), c_1(D, \alpha, \gamma))$$ 
with the lexicographical order to be the complexity of $(P, D)$.
\end{definition}

 As mentioned above, the general idea of the proof is to show one can reduce the complexity of the diagram. 
 The following proposition shows that if one reaches a diagram with small complexity then $K$ is a Berge 
 knot. It is thus a main ingredeant  of  the proof. This proposition also deals with the situation described in 
 Remark \ref{rem:NotAllSurfaces}, where a knot space has  multiple  Heegaard splittings in which  the knot  
 $K$ is not doubly primitive with respect to some and is doubly primitive with respect to others.
 
\begin{proposition}\label{pro:Tao}
Let $(P,D)$ be a $(\mathcal{P},\mathcal{D})$-pair and  $\alpha = \partial_+ A$ as above. Suppose  $|\,P\cap A\,|$ is minimal up to isotopy on $P$ and $A$. If $c_0(P, D, \alpha) \leq 1$ then $K$ is 
doubly primitive with respect to a (possibly different) genus two Heegaard splitting of $M$. 
\end{proposition}

\begin{proof}
We have assumed at the beginning of section~\ref{sec:complexity} that $P$ is incompressible. This implies $P$ does not have a $\partial$-compressing disk at the 
boundary torus $\partial \NN(K)$. To see this, if $P$ has such a $\partial$-compressing disk, after boundary compressing 
$P$ along this disk, we obtain a new surface $P'$ with a boundary curve bounding a disk on the boundary torus
$\partial \NN(K)$. If this curve bounds a disk in $P'$ as well, then since $M\ssm\NN(K)$ is irreducible, $P$ must be a 
boundary parallel annulus. If this curve is essential in $P'$, then the disk in $\partial \NN(K)$ gives rise to a 
compressing disk for $P$.  Both cases contradict the choice of $P$.

 By Lemma~\ref{lem:integer slope}, we may assume that the  curves in $\partial_-P$ have integer slope. 
 Hence each component of $\partial_-P$ intersects $\partial_-A$ in exactly one  point. Since $P$ is 
 incompressible and $|\,P\cap A\,|$ is minimal we may assume that, perhaps after isotopy,  $P\cap A$ 
 contains no closed curve that is trivial in $A$. Since $P$ does not have $\partial$-compressing disk at 
 the boundary torus, after isotopy, $P\cap A$ contains no arc with both endpoints in $\partial_-A$.  Further, 
 each component of $A\cap P$ must be either an essential arc in $A$ or an arc in $A$ with both endpoints
  in $\alpha=\partial_+A$. 
 
Let $k_1,\dots, k_N$ be arcs in  $A\cap P$ that are essential in $A$.  Thus $k_1,\dots, k_N$ divide $A$ into $N$ 
rectangles $R_1,\dots R_N$.   By Proposition~\ref{pro:DPEquivalentToAD}, we may assume $P$ is not an annulus 
and hence $N\ge 2$.

Suppose $c_0(P, D, \alpha) \leq 1$. Then the intersection $\alpha\cap\partial D$, if nonempty, must lie in a single 
arc of $\alpha\ssm P$ and thus in the boundary of a single rectangle $R_i$, and without loss of generality 
we may assume that $\alpha\cap\partial D\subset\partial R_1$.  Set $\widehat{R}=R_2\cup\cdots\cup R_N$. 
If $P\cap\widehat{R}$ contains arcs that are $\partial$-parallel in $A$, we can perform $\partial$-compressions 
on $P$ along these arcs.  Since $\widehat{R}\cap D=\emptyset$, such $\partial$-compressions do not create 
new intersection points with $D$.  Hence, after the $\partial$-compressions, one of the resulting planar surfaces 
still  intersects $D$ in a single point.  Thus we  may assume that no arc in $P\cap\widehat{R}$ is $\partial$-parallel 
in $A$, in other words, $P\cap\widehat{R}=\cup_{i=1}^N k_i$.

Let $K^*$ denote the dual knot to $K$ in the corresponding lens space $L$, where  
$L\ssm\NN(K^*)\cong M\ssm\NN(K)$.  We view $P$, $D$ and $A$ in the exterior of $K^*$.  The planar surface 
$P$ can be viewed as an $N$-punctured disk in $L\ssm \NN(K^*)$. The arcs $k_1,\dots, k_N$, when viewed 
in $P$, are arcs connecting the $N$ components of $\partial_-P$ to $\partial_+P$.  Let $\mu_i$ be the 
component of $\partial_-P$ that contains an endpoint of $k_i$ ($i=1,\dots,N$), and let $N_{\mu_i}$ be the 
closure of a small neighborhood of $\mu_i\cup k_i$ in $P$. So $N_{\mu_i}$ is an annulus of which $\mu_i$ 
is a boundary component.  We may view the annulus $N_{\mu_i}$ as a sub-surface of $P$ cut off by an arc $d_i$, 
where $d_i\subset\partial N_{\mu_i}$ and $d_i$ is properly embedded in $P$, see Figure~\ref{Finter}(a).  
The arc $d_i$ is the arc in $P$ that goes around $\mu_i\cup k_i$ as indicated in Figure~\ref{Finter}(a).

\begin{figure} 
\begin{overpic}[width=10cm]{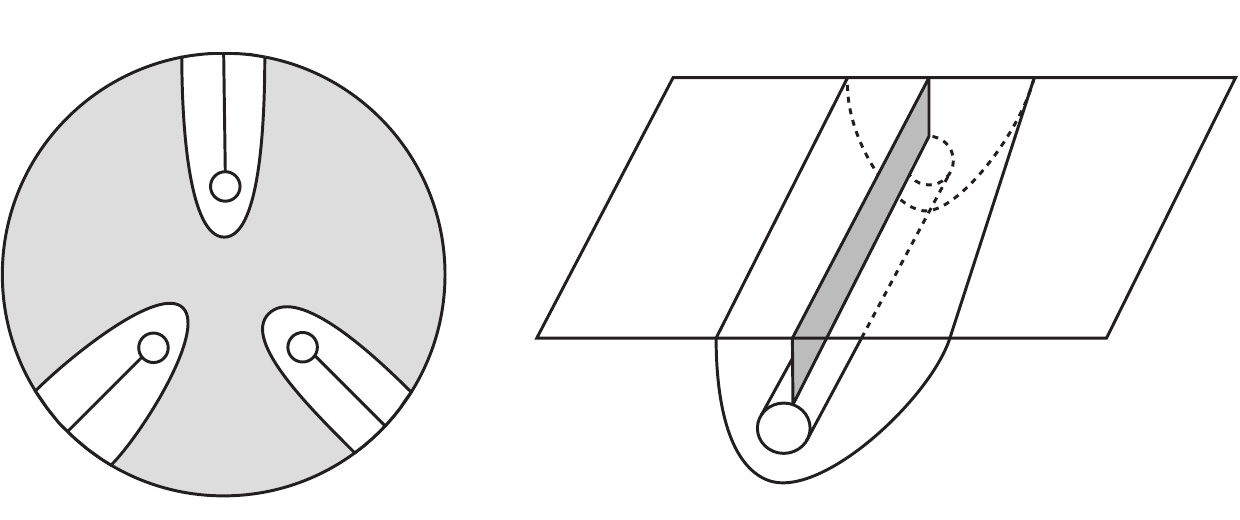}
\put(5,28){$\Delta_P$}
\put(18.3,34){{\footnotesize $k_i$}}
\put(21.5,29){$d_i$}
\put(31,32){$P$}
\put(47,17){$\Sigma$}
\put(49,9){$U$}
\put(75,31.5){{\footnotesize $\widehat{R}$}}
\put(75.5,9){{\small $\widehat{R}_d$}}
\put(15,-4){(a)}
\put(68,-4){(b)}
\end{overpic}
\vskip5pt
\caption{The surfaces $\Delta_P$ and $\widehat{R}_d$ used to modify the Heegaard splitting.}\label{Finter}
\end{figure}

Since $k_1,\dots, k_N$ are vertical arcs of the rectangle $\widehat{R}$ there is a rectangle $\widehat{R}_d=I\times I$
``following"  $\widehat{R}$ so that each $d_i$ is a vertical arc in $\widehat{R}_d$ of the form of $\{x\}\times I$ 
and $\widehat{R}_d\cap \Sigma=I\times\partial I$, see  Figure~\ref{Finter}(b) (the shaded region denotes $\widehat{R}$).
The rectangle  $\widehat{R}_d$ can be obtained by first taking two parallel copies of $\widehat{R}$ and then connecting 
them around $K^*$. Note that $\{d_1,\dots, d_N\}$  is a collection of disjoint vertical arcs in $\widehat{R}_d$.
Since $\widehat{R}$ is disjoint from $D$ we have $\widehat{R}_d\cap D=\emptyset$.

Next we  modify the Heegaard splitting $(U, W)$ of $M\ssm\mathcal{N}(K)=L\ssm\mathcal{N}(K^*)$.  
Let $U'=U\ssm\NN(\widehat{R}_d)$.  Note that $\widehat{R}_d$ deformation retracts to $d_i$ for 
any $i$. So $U'\cong U\ssm\NN(d_i)$ and one may view $\mathcal{N}(\widehat{R}_d)$ as a 
``fat" tunnel in $U$.  As $U$ is a compression body, there is a vertical annulus $A_{\mu_i}\subset U$ 
between any $\mu_i$ and a curve in $\Sigma$.  
Notice that, after isotopy, $A_{\mu_i}\cap A$ is a single arc $k$ which is vertical in both $A_{\mu_i}$ and $A$. Hence $k_i$ is isotopic to $k$ in the annulus $A$. As $k\subset A_{\mu_i}$, we can isotope $\mu_i\cup k_i$ into $A_{\mu_i}$. 
Hence the annulus $N_{\mu_i}$ can be isotoped 
to a subsurface of  $A_{\mu_i}$. In particular $d_i$ can be isotoped into $A_{\mu_i}$.  

As $d_i$ has  both end points on $\Sigma$, $d_i$ is $\partial$-parallel in  $A_{\mu_i}$ and hence 
$\partial$-parallel  in $U$. Thus drilling out a tunnel in $U$ with the arc $d_i$ as its core results 
in a compression body isotopic to $U'=U\ssm\mathcal{N}(\widehat{R}_d)$. The complement 
of $U'$ is a genus-three handlebody $W'$.  Thus  $(U', W')$ is a stabilized genus-three Heegaard splitting of   
$M\ssm\mathcal{N}(K)=L\ssm\mathcal{N}(K^*)$.  

Note that  $N_{\mu_i}$ is a vertical annulus in $U'$. Let $\Delta\subset W'$ be a cocore 
(meridional) disk of  the (fat) tunnel $\mathcal{N}(\widehat{R}_d)$ that we drill out from $U$. 
So $\Delta\cap N_{\mu_i}$ is a  single point and $(N_{\mu_i}, \Delta) $ is an 
$(\mathcal{A}, \mathcal{D})$-pair.  

Let $\Delta_P$ be the closure of $P\ssm \cup_{i=1}^N N_{\mu_i}$, see the shaded region in 
Figure~\ref{Finter}(a). Since $\widehat{R}$ 
connects all the components of $\partial_-P$, $\Delta_P$ is a disk which can be viewed as a properly 
embedded disk in $U'$. Since $\widehat{R}_d$ is disjoint from $D$,  the operation of drilling this tunnel 
does not affect $D$. Hence we may view $D$ as a disk in $W'$ and $D\cap\Delta_P=D\cap P$ is a single point. 
Thus $(D,\Delta_P)$ is a stabilizing pair of disks for the genus-$3$ splitting $(U',W')$. By cutting $W'$ along $D$ and adding $\NN(D)$ as a $2$-handle 
to $U'$ along $\partial D$, we destabilize the splitting $(U', W')$ into a new genus two Heegaard splitting 
of $M \ssm \NN(K)$.  As $D$ is disjoint from both $N_{\mu_i}$ and $\Delta$, the pair $(N_{\mu_i}, \Delta)$ 
described above remains an $(\mathcal{A}, \mathcal{D})$-pair in the resulting genus two Heegaard splitting. 
By Proposition~\ref{pro:DPEquivalentToAD}, $K$ is doubly primitive with respect to  this new genus-two 
Heegaard splitting.
\end{proof}

\vskip5pt

By Definition~\ref{def:complexity}, $c_0(P, D,\alpha)\le c_0(P, D,\alpha, \gamma)$. Hence $K$ is doubly primitive if  $c_0(P, D,\alpha, \gamma)\le 1$. 
In light of Proposition \ref{pro:Tao} we we adopt the following:

\vskip5pt

\begin{assumption}\label{ass:minimality}
Fix a $(\mathcal{P}, \mathcal{D})$-pair $(P, D)$ and suppose $P$ is incompressible in $U$. 
By Proposition~\ref{pro:DPEquivalentToAD}, we may assume $P$ is not an annulus. 
Throughout the paper, we assume $P$, $D$ and $\alpha$ are chosen so that 
$c(P, D, \alpha, \gamma)$ is minimal among all the $(\mathcal{P},\mathcal{D})$-pairs and 
all such annuli $A$ ($\alpha=\partial_+A$). 
\end{assumption}

We finish this section with the following lemma.

\begin{lemma}\label{lem:PCIntersection} Let $P$ and $\gamma$ be as in Assumption~\ref{ass:minimality}. Then $\partial_+P\cap\gamma\ne\emptyset$. 
\end{lemma}
\begin{proof}
Suppose on the contrary that $\partial_+P\cap\gamma=\emptyset$. 
Recall that $\gamma=\partial C$ where $C$ is a compressing disk in $U$.  Since $P$ is incompressible in $U$ and since $\partial_+P\cap\gamma=\emptyset$, we can perform an isotopy so that $P\cap C=\emptyset$.  

If we cut $U$ open along the disk $C$, then the resulting manifold is 
homeomorphic to $T^2 \times I$ and $P\subset T^2\times I$.  Note that if $P$ is $\partial$-compressible 
in $T^2\times I$, then $P$ is compressible (one can construct a compressing disk using two parallel 
copies of a $\partial$-compressing disk).  Since $P$ is assumed to be incompressible, $P$ must be both 
incompressible and $\partial$-incompressible and hence $\pi_1$-injective in $T^2\times I$.  This means that 
$P$ must be an annulus, contradicting Assumption~\ref{ass:minimality}.  
\end{proof}

\vskip10pt


\section{Meridional disks and Whitehead graphs}\label{sec:MeridionalDisksandWhiteheadGraphs}

In this section we define Whitehead graphs and \emph{waves} which are the main tools in the proof 
of the theorem. 

Let  $(V, W)$ be a  genus two Heegaard splitting of a closed $3$-manifold $M$ with a Heegaard surface 
$\Sigma$.  Let $\widetilde{V}=\{V_1, V_2\}$ and  $\widetilde{W}=\{W_1,W_2\}$  be {\it complete meridian sets} 
respectively. That is, each is a pair of essential disks so  that $V\ssm \widetilde{V}$ and $W \ssm \widetilde{W}$ 
are $3$-balls.  Let $v_i=\partial V_i$ and $w_i=\partial W_i$ ($i=1,2$). To simplify notation, we also call $v_i$ 
and $w_i$ meridians of $V$ and $W$ respectively.  If we cut $\Sigma$ open along $v_1\cup v_2$, we get a 
$4$-punctured  sphere whose boundary consists  of four components  $\{v^+_1, v^-_1, v^+_2, v^-_2\}$, and 
$w_1\cup w_2$ is cut into a collection of properly embedded curves in the  $4$-punctured  sphere.

Let $\widehat{V} = \{ v_1, v_2\}$ and $\widehat{W} = \{ w_1, w_2\}$. 
A Heegaard diagram $(\Sigma, \widehat{V}, \widehat{W})$ for the $3$-manifold $M$ is said 
to be {\it normal} if no domain in  $\Sigma\ssm  \{\widehat{V} \cup \widehat{W}\}$ is a bigon, see \cite{Och} for more detailed definitions. 
We can always isotope the curves $v_1$, $v_2$, $w_1$ $w_2$ in $\Sigma$ to eliminate all the bigons and obtain a normal Heegaard diagram. 
Throughout the paper, we assume all the Heegaard diagrams are normal.

\begin{definition}[Genus two Whitehead graphs] \label{def:WhiteheadGraphs} The graph 
$\Gamma(\{v_1, v_2\})$ (or $\Gamma(\widehat{V})$) obtained, as above, by setting  $\{v^+_1, v^-_1, v^+_2, v^-_2\}$ to be
vertices and the arcs of $(w_1\cup w_2) \ssm (v_1\cup  v_2)$ to be  edges is called the 
{\it Whitehead graph}  corresponding to  $\{v_1, v_2\}$. As $\Sigma\ssm \widehat{V}$ is a 4-hole sphere, $\Gamma(\{v_1, v_2\})$ is a planar graph in the sphere. 
Similarly we have the Whitehead  graph 
$\Gamma(\{w_1, w_2\})$ (or $\Gamma(\widehat{W})$) corresponding  to $\{w_1, w_2 \}$.
\end{definition}

Ochiai  states, in ~\cite [Theorem 1] {Och}, that any Whitehead graph of a genus two Heegaard 
diagram of a $3$-manifold is isomorphic as a planar graph to one of the three graphs $(i), (ii)$  
or  $(iii)$ in Figure \ref{graphs} below, where the integers $a,\ b, c$ and $d$ represent the 
number, which can be zero, of  parallel arcs corresponding to each graph edge.  Note that 
it is easy to understand this theorem and Figure~\ref{graphs} using the hyperelliptic involution $\pi$: 
Each nonseparating simple closed curve in a genus two surface is invariant (up to isotopy) under 
 $\pi$ with orientation reversed. So the involution interchanges $v_i^+$ and 
$v_i^-$ ($i=1,2$) and induces an involution on the Heegaard diagram. For example, as $w_i$ 
is invariant under the involution, the numbers of arcs connecting $v_1^+$ to $v_2^+$ and $v_1^-$ 
to $v_2^-$ are the same (indicated by $a$ in Figure~\ref{graphs}).  Furthermore, this theorem is 
about the graphs only. If one interchanges $v_1$ and $v_2$ or switches the $\pm$-signs of $v_i$, 
then the diagram is of the same type.  Note that if $w_i$ is parallel to $v_i$ (i.e.~$M$ has an 
$S^2\times S^1$ summand) then there are no $w_i$ arcs in $\Gamma(\{v_1,v_2\})$.

\begin{figure}[h!]
\subfigure[]{%
\begin{overpic}[width=4.5cm]{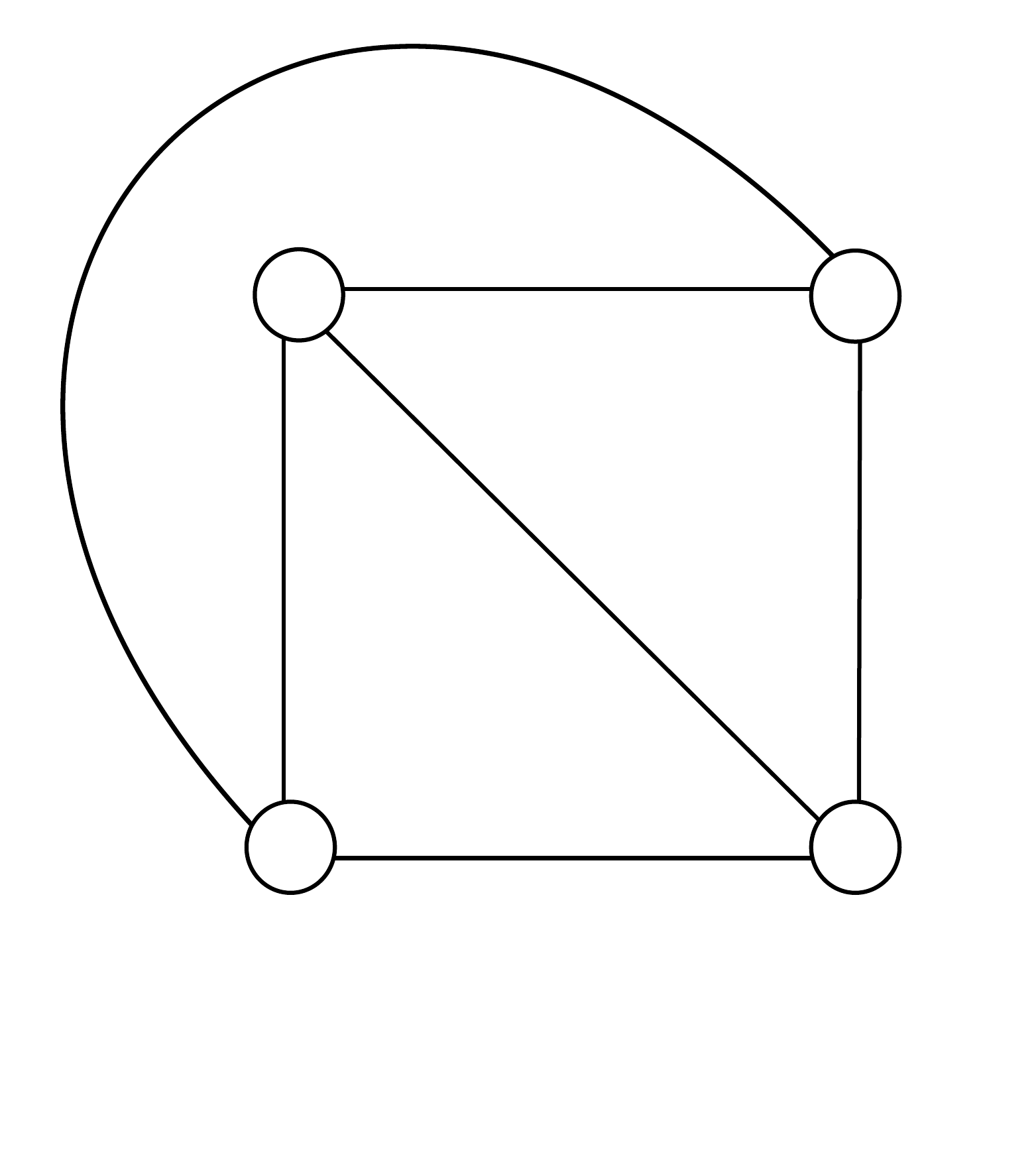}
\put(12,68){$v_1^+$}
\put(12,18){$v_1^-$}
\put(76,69){$v_2^+$}
\put(76,19){$v_2^-$}
\put(45,77){$a$}
\put(45,21){$a$}
\put(6,84){$b$}
\put(48,52){$b$}
\put(18,47){$c$}
\put(75,47){$d$}
\end{overpic}}
\subfigure[]{%
\begin{overpic}[width=4.5cm]{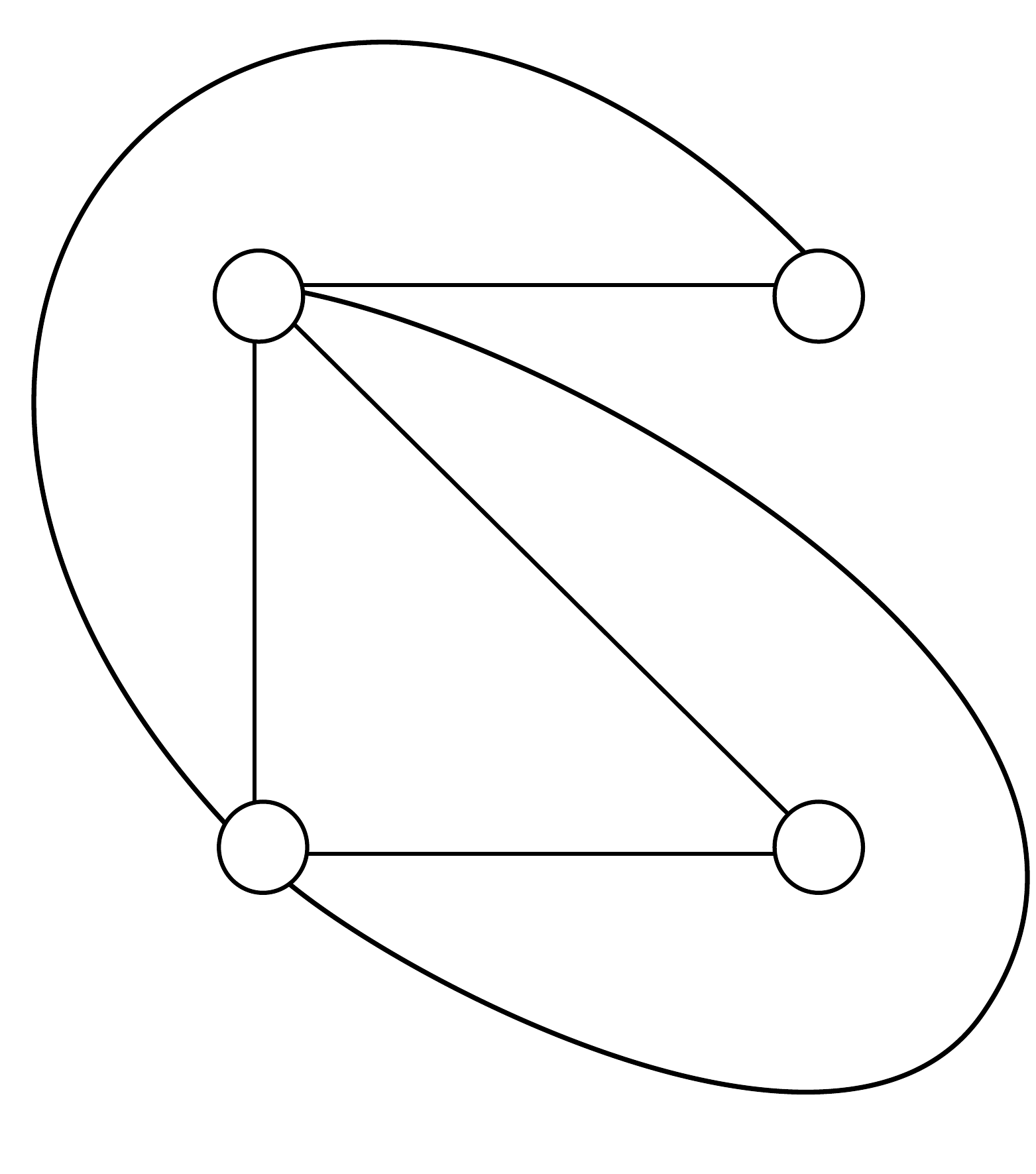}
\put(9,74){$v_1^+$}
\put(9,20){$v_1^-$}
\put(74,70){$v_2^+$}
\put(74,20){$v_2^-$}
\put(45,77){$a$}
\put(45,20){$a$}
\put(7,88){$b$}
\put(46,53){$b$}
\put(16,50){$c$}
\put(74,50){$d$}
\end{overpic}}
\subfigure[]{%
\begin{overpic}[width=4.5cm]{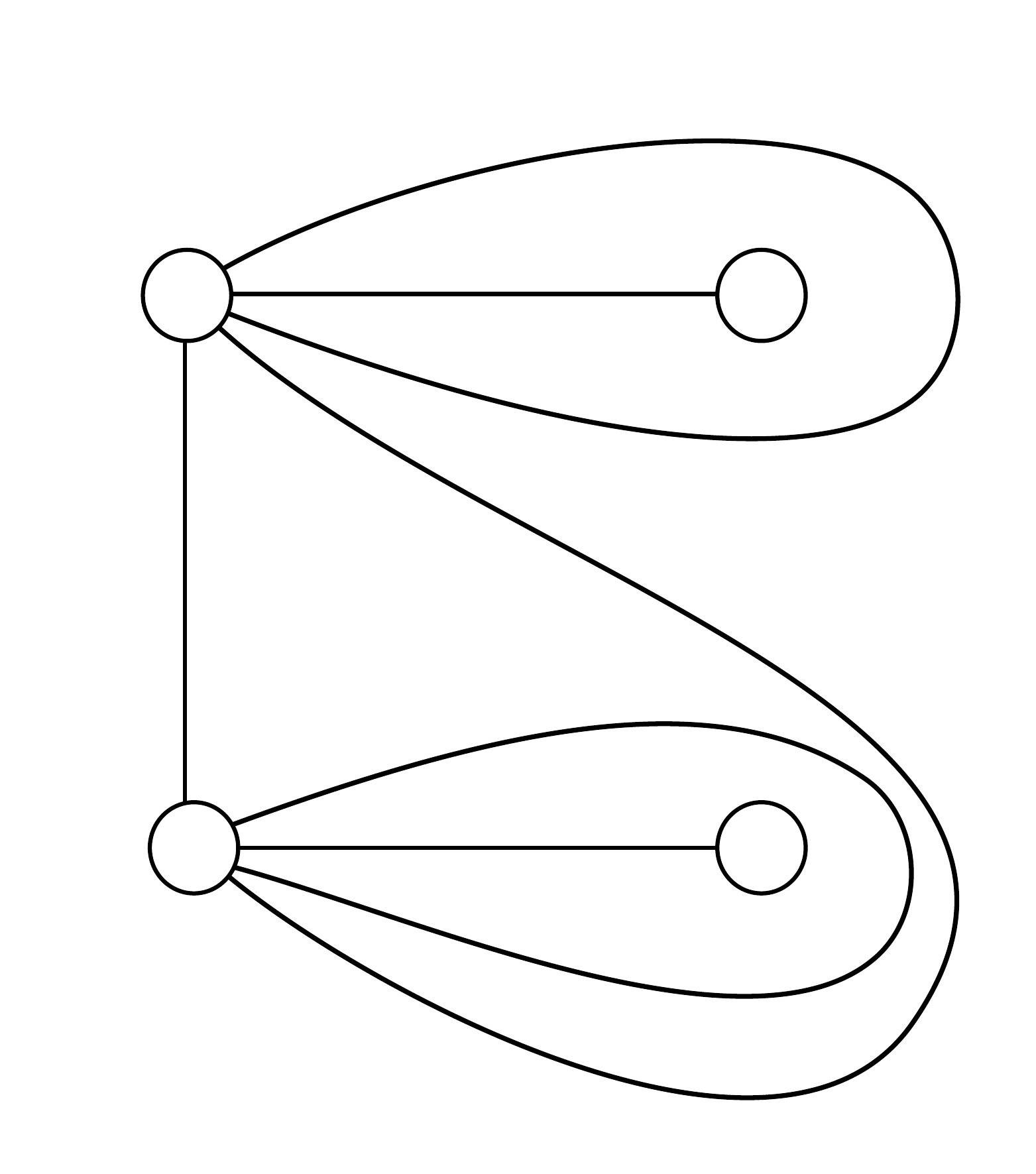}
\put(3,74){$v_1^+$}
\put(3,25){$v_1^-$}
\put(69,72){$v_2^+$}
\put(68,22){$v_2^-$}
\put(45,76){$a$}
\put(45,29){$a$}
\put(45,88){$b$}
\put(46,39){$b$}
\put(17,50){$c$}
\put(56,50){$d$}
\end{overpic}}
\caption{The three graphs given in \cite{Och}. }
\label{graphs}
\end{figure}

\begin{definition}\label{def:Shortcuts}\hfill

 \noindent (1)  A {\it standard wave $\zeta$ } with respect to a closed curve $\kappa$ on an orientable surface  $\Sigma$ 
 is a simple arc so that $\zeta \cap \kappa = \partial \zeta$,  the arc $\zeta$ intersects $\kappa$ from the 
 same side and $\zeta$ is not homotopic into a subarc of $\kappa$ rel $\partial \zeta$. 
\vskip 3pt 
 \noindent (2) When $\Sigma$ is the boundary of a genus two handlebody $V$ with a fixed  meridian system
 $\{v_1, v_2\}$, a standard wave $\zeta$ with respect to $v_i$ is called an \emph{s-wave}  if $\zeta$ is disjoint 
 from $v_j$ ($j\neq i$) and $\Sigma\ssm(v_i\cup\zeta)$ is connected.
  \vskip 3pt 
 \noindent (3) In a Heegaard diagram $(\Sigma, \widehat{V}, \widehat{W})$ of some $3$-manifold $M$, 
 where $\widehat{V}$ and $\widehat{W}$ are complete sets of meridians, an $s$-wave $\eta \subset \Sigma$  
 with respect to  a meridian $v_i \subset \widehat{V}$ is called a {\it wave} with respect to $v_i$ if 
 $\eta\cap \widehat{W} =\emptyset$ or $\eta$ is a subarc of $\widehat{W}$.
\end{definition}

For example, in the Whitehead graph depicted in Figure~\ref{graphs}(iii), the edges labeled $b$ 
correspond to subarcs of $w_1\cup w_2$ that are waves with respect to $v_1$.

\begin{remark}\label{rem:duralwave}
	Note that the configurations in Figure~\ref{graphs} are symmetric under the hyperelliptic involution $\pi$. 
	Therefore, if there is a wave $\eta$ with endpoints in $v_i^+$, there must be a wave $\eta'=\pi(\eta)$ with 
	endpoints in $v_i^-$, for example, see the two edges labeled $b$ in Figure~\ref{graphs}(iii). Moreover, $\eta$ and $\eta'$ are attached to opposite sides of $v_i$ and hence are distinct waves. 
\end{remark}

\begin{definition}\label{def:wavemove}
Let $\widehat{V}=\{v_1,v_2\}$ be as above. 
In case there exists an $s$-wave $\eta$ with respect to $v_i$, each of the three boundary components of 
$\NN(v_i\cup\eta)$ bounds a disk in the handlebody $V$. Since $V$ is of genus two, we necessarily 
have two of the boundary components isotopic to $v_1$ and $v_2$  respectively, while  the third 
component is the boundary of a new nonseparating disk. We call the operation of replacing the disk 
$V_i$ in the complete meridian set $\widetilde{V}$ by the new disk the {\it wave move} of $v_i$ along $\eta$. 
In this paper we often view the wave move as a two step operation: The first step is a surgery on $v_i$ which 
connects the endpoints of  $v_i\ssm\NN(\partial\eta)$ using two parallel copies of $\eta$, thus obtaining 
two curves. Note that  one of these curves  is parallel to $v_j$ ($j\neq i$). In the second step we delete the 
curve parallel to $v_j$ ($j\neq i$) and the remaining  curve is the boundary of the new disk  obtained 
from the wave move.
\end{definition}

\begin{remark} Note that $\NN(v_i\cup\eta)$ in Definition~\ref{def:wavemove} is a pair of pants in $\Sigma$ with two of its boundary 
components isotopic to $v_1$ and $v_2$ respectively.  So the new curve obtained by the wave move 
can also be obtained by a band sum of $v_1$ and $v_2$ along an arc in $\NN(v_i\cup\eta)$.
\end{remark}

\begin{definition}\label{def:standard}
Meridian sets $\widehat{V} = \{ v_1, v_2\}$ and $\widehat{W} = \{ w_1, w_2\}$ for a genus two Heegaard splitting 
of $S^3$ are {\it standard} if $|v_1\cap w_1|=|v_2\cap w_2|=1$ and $v_1\cap w_2=v_2\cap w_1=\emptyset$.
Meridian sets $\widehat{V} = \{ v_1, v_2\}$ and $\widehat{W} = \{ w_1, w_2\}$ for a genus two Heegaard splitting 
of $(S^1 \times S^2) \# L(r ,s)$ are {\it standard} if both conditions below are satisfied:
\vskip5pt
\noindent (1) The curves $v_1$ and $w_1$ are parallel and disjoint from $v_2\cup w_2$, and
\vskip5pt
\noindent (2) There is a separating essential simple closed curve in $\Sigma$ disjoint from 
$v_1\cup v_2\cup w_1\cup w_2$. 
\vskip5pt
\noindent
Note that a separating curve in (2) bounds disks in both handlebodies and the two 
disks form a separating essential 2-sphere in $M$. 
A genus two Heegaard diagram of $S^3$ or 
$(S^1 \times S^2) \# L(r, s)$ is {\it standard} if the corresponding meridian sets are standard.  
\end{definition}

We are now ready to state the following theorems which play a crucial role in the argument. 
The main theorem of Homma, Ochiai and Takahashi in \cite{HOT} states:

\begin{theorem}[Homma, Ochiai, Takahashi \cite{HOT}] \label{thm:HOT} Any genus two Heegaard diagram 
for $S^3$  either is standard or contains a wave.  
\end{theorem}

For $M = (S^1 \times S^2) \# L(r, s)$ we have the following theorem of Negami and Okita: 

\begin{theorem}[Negami-Okita  \cite{NeOk}]\label{thm:NeOk} Any genus two Heegaard diagram of  
$(S^1 \times S^2) \# L(r, s)$ is either standard or contains a wave.
\end{theorem}

\begin{lemma}\label{lem:+ and -} Suppose  $\widetilde{V}=\{V_1, V_2\}$ and $\widetilde{W}=\{W_1,W_2\}$  
are complete meridian sets  for a genus two Heegaard splitting $(V, W)$ for a $3$-manifold $M$ 
and the Whitehead graph of $\Gamma(\widehat{V})$ has edges connecting the vertices as 
in one of the following cases:
\begin{enumerate}
 \item  there are edges connecting $v_1^+$ to $v_1^-$ and  $v_2^+$ to $v_2^-$, or 
 \item there are edges connecting both $v_2^+$ and $v_2^-$  to the same vertex $v_1^+$ or $v_1^-$, or
 \item there are edges connecting both $v_1^+$ and $v_1^-$ to the same vertex $v_2^+$ or $v_2^-$. 
  \end{enumerate}

Then  the Heegaard diagram $(\widehat{V}, \widehat{W})$ contains no waves with respect to $\widehat{V}$.
\end{lemma} 

\begin{proof} Since this is a Heegaard  diagram of a $3$-manifold,   Theorem 1 of ~\cite{Och}  applies and 
the Whitehead graph is as in Figure ~\ref{graphs}.  For case (1), the Heegaard diagram must be of type (i) 
in  Figure~\ref{graphs} with $c>0$ and $d>0$.  Note that a wave with respect to $v_j$ must separate $v_i^+$ 
from $v_i^-$ ($i\ne j$), see the edge labeled $b$ in  Figure~\ref{graphs}(iii) for an example of a wave. 
So an edge connecting $v_i^+$ to $v_i^-$ prevents any wave with respect to $v_j$, 
$i \neq j, \, i,j \in \{1,2\}$. Hence Case (1) holds.

Case (2) and (3) are symmetric.  For case (2), the Heegaard diagram must be of type (i) or (ii) in  Figure~\ref{graphs} 
with $a > 0$ and $b > 0$.  Similar to case (1), since the edges labeled $a$ and $b$ in Figure~\ref{graphs} (i) and (ii) 
connect $v_1^+$, $v_1^-$, $v_2^+$ and $v_2^-$ together, there cannot be any wave  with respect to  $v_1$   
and  $v_2$.  
\end{proof}

\begin{definition}\label{def: blockingedges} 
Throughout the paper, for any different vertices $v_i^{\pm}$ and $v_j^{\pm}$, we use $[v_i^\pm, v_j^\pm]$ to denote an edge of the Whitehead graph $\Gamma(\widehat{V})$ connecting $v_i^\pm$ to $v_j^\pm$.  If $[v_i^\pm, v_j^\pm]$ is an edge in either Case $(1)$, $(2)$  or  $(3)$  of 
Lemma~\ref{lem:+ and -}, we say this edge is a $[v_i^\pm, v_j^\pm]$ {\it blocking-edges} with respect to $V$. 
Lemma~\ref{lem:+ and -} says that,  
if the Whitehead graph $\Gamma(\widehat{V})$ has both types of blocking-edges in either Case (1), (2), or (3) of 
Lemma~\ref{lem:+ and -}, then the Heegaard diagram $(\widehat{V}, \widehat{W})$ contains no waves with respect to $\widehat{V}$. 
The phrase ``blocking-edge" is used in the paper whenever we are searching for edges in Case (1), (2), or (3) of Lemma~\ref{lem:+ and -}. 
We would like to emphasize that one type of blocking-edge alone does not ``block" waves, and we need both types in Case (1), (2), or (3) of Lemma~\ref{lem:+ and -}.  
The corresponding statement of Lemma~\ref{lem:+ and -} also holds for the Whitehead graph $\Gamma(\widehat{W})$ and we can similarly define {\it blocking-edges} with respect to $W$.
\end{definition}

For  example, if two adjacent intersection points of $w_i\cap(v_1\cup v_2)$ along $w_i$ both belong to $v_1$ 
and have the same sign of intersection, then the subarc of $w_i$ between these two points is a 
$[v_1^+, v_1^-]$ blocking-edge.  Moreover, if there is a similar subarc of $w_j$ with endpoints in $v_2$ and having the same sign of intersection, then we also have a $[v_2^+, v_2^-]$ blocking-edge. 
By Case (1) of Lemma~\ref{lem:+ and -}, this means that the Heegaard diagram has no wave with respect to $V$.

We finish this chapter with the following lemma on the Whitehead graphs.

\begin{lemma}\label{lem:dual}
Let $\widehat{V}=\{v_1, v_2\}$ and $\widehat{W}=\{w_1, w_2\}$ be complete sets of meridians. 
Suppose $v_1\cup v_2$ contains a $[w_1^+, w_1^-]$ blocking-edge and a $[w_2^+, w_2^-]$ 
blocking-edge. Suppose there is a wave with respect to $\{v_1, v_2\}$ and suppose 
$\Gamma(\widehat{W})$ is a connected graph. 
\begin{enumerate}
\item If no subarc of $w_1\cup w_2$ is a wave with respect to $\{v_1, v_2\}$, then the Whitehead graphs 
$\Gamma(\widehat{W})$ and $\Gamma(\widehat{V})$ are as shown in Figure~\ref{dual1} where all the labels are 
non-zero.
\item If a subarc of $w_1$ or $w_2$ is a wave with respect to $\{v_1, v_2\}$, then the Whitehead graphs 
$\Gamma(\widehat{W})$ and $\Gamma(\widehat{V})$ are as shown in Figure~\ref{dual2} where all the labels 
are non-zero.
\end{enumerate}
\end{lemma} 
\begin{proof}
Since $v_1\cup v_2$ contains a $[w_1^+, w_1^-]$ blocking-edge and a $[w_2^+, w_2^-]$ blocking-edge, 
the Whitehead graph $\Gamma(\widehat{W})$ must be of type (i) in Figure~\ref{graphs} with $c\ge 1$ and 
$d\ge 1$.  If $a=b=0$ in Figure~\ref{graphs}(i), then the Whitehead graph $\Gamma(\widehat{W})$ is disconnected, 
a contradiction to the hypothesis.  Hence at least one of $a$ and $b$ in Figure~\ref{graphs}(i) is non-zero.  
If $b=0$ and $a\ne 0$ in Figure~\ref{graphs}(i), then $\Gamma(\widehat{W})$ is exactly the left picture of Figure~\ref{dual1}. If $a=0$ and $b\ne 0$ in Figure~\ref{graphs}(i), then by interchanging the locations of $v_2^+$ and $v_2^-$ in Figure~\ref{graphs}(i),  $\Gamma(\widehat{W})$ is isotopic to the left picture of Figure~\ref{dual1}. If neither 
$a$ or $b$ is zero in Figure~\ref{graphs}(i),  then $\Gamma(\widehat{W})$ is the left picture of Figure~\ref{dual2}.

If $\Gamma(\widehat{W})$ is as in the left picture of Figure~\ref{dual1}, then 
$\Sigma\ssm(v_1\cup v_2\cup w_1\cup w_2)$ consists of two octagons and a collection of quadrilaterals 
(i.e.~the regions between parallel arcs of the same edge type).  If  $\Gamma(\widehat{W})$ is as in the left 
picture of Figure~\ref{dual2}, then $\Sigma\ssm(v_1\cup v_2\cup w_1\cup w_2)$ consists of four hexagons 
and a collection of quadrilaterals. 

\begin{figure}[ht]
	\begin{overpic}[width=7cm]{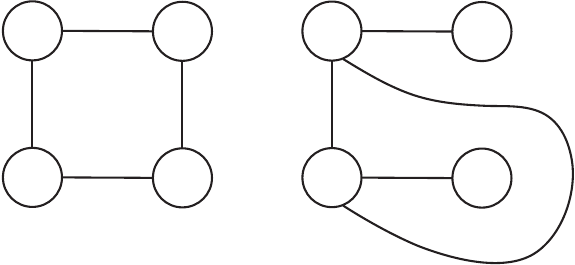}
		\put(14,-6){$\Gamma(\widehat{W})$}
		\put(67,-6){$\Gamma(\widehat{V})$}
		\put(2,39){$w_1^+$}
		\put(2,13){$w_1^-$}
		\put(55,39){$v_1^+$}
		\put(55,13){$v_1^-$}
		\put(18,41.5){$a$}
		\put(18,11.5){$a$}
		\put(2,26){$c$}
		\put(32,26){$d$}
		\put(69,41.5){$a'$}
		\put(69,11){$a'$}
		\put(53,26){$c'$}
		\put(84,28){$d'$}
\put(28,39){$w_2^\pm$}
\put(28,13){$w_2^\mp$}
\put(80,39){$v_2^\pm$}
\put(80,13){$v_2^\mp$}
	\end{overpic}
\vspace{8pt}
	\caption{Dual Whitehead graphs: case (1)}\label{dual1}
\end{figure}

Next, consider the Whitehead graph  $\Gamma(\widehat{V})$. Since there is a wave with respect to 
$\{v_1, v_2\}$ and since the diagram is connected, the Whitehead graph  $\Gamma(\widehat{V})$ must 
be of type (iii) in Figure~\ref{graphs} with $a\ne 0$ and at least one of $c$ and $d$ non-zero. Without 
loss of generality, suppose $c\ne 0$. 

\begin{figure}[ht]
	\begin{overpic}[width=8cm]{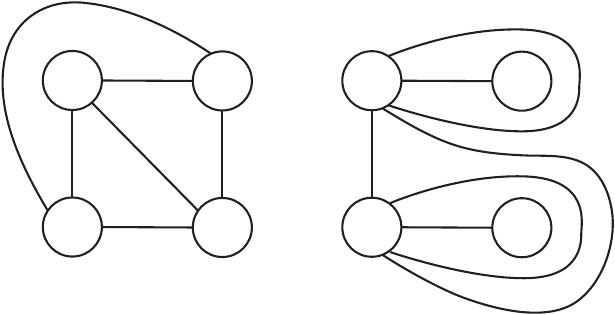}
		\put(18,-6){$\Gamma(\widehat{W})$}
		\put(69,-6){$\Gamma(\widehat{V})$}
		\put(9,37){$w_1^+$}
		\put(9,13){$w_1^-$}
		\put(57,37){$v_1^+$}
		\put(57,13){$v_1^-$}
		\put(22,39){$a$}
		\put(22,10.5){$a$}
		\put(23,28){$b$}
		\put(-2,41){$b$}
		\put(8,26){$c$}
		\put(36.5,26){$d$}
		\put(72,38.5){$a'$}
		\put(73,15){$a'$}
		\put(71,45){$b'$}
		\put(70,21){$b'$}
		\put(56.5,25){$c'$}
		\put(94,25){$d'$}
		\put(32.5,37){$w_2^\pm$}
		\put(32.5,13){$w_2^\mp$}
		\put(81.5,37){$v_2^\pm$}
		\put(81.5,13){$v_2^\mp$}
	\end{overpic}
	\vspace{6pt}
	\caption{Dual Whitehead graphs: case (2)}\label{dual2}
\end{figure}

If $b=d=0$ in Figure~\ref{graphs}(iii), then the picture of Figure~\ref{graphs}(iii) indicates that 
$\Sigma\ssm(v_1\cup v_2\cup w_1\cup w_2)$ consists of one 12-gon and a collection of 
quadrilaterals. This configuration is incompatible with the two possible configurations of 
$\Gamma(\widehat{W})$ discussed above. So this cannot happen and at least one 
of $b$ and $d$ in Figure~\ref{graphs}(iii) is non-zero.

If $b\ne 0$ and $d=0$ in Figure~\ref{graphs}(iii), then $\Sigma\ssm(v_1\cup v_2\cup w_1\cup w_2)$ 
consists of one octagon, two hexagons and a collection of quadrilaterals.  This configuration is also 
incompatible with the two possible configurations of $\Gamma(\widehat{W})$ discussed above. 
So this cannot happen either. 

If $b=0$ and $d\ne 0$ in Figure~\ref{graphs}(iii), then $\Gamma(\widehat{V})$ is as in the right 
picture of Figure~\ref{dual1}, in which case $\Sigma\ssm(v_1\cup v_2\cup w_1\cup w_2)$ consists of
two octagons and a collection of quadrilaterals. By the description of $\Gamma(\widehat{W})$ above, 
this means that $\Gamma(\widehat{W})$ must be as in the left picture of Figure~\ref{dual1}. 
Note that an edge in $\Gamma(\widehat{V})$ is a wave if and only if $b\ne 0$ in Figure~\ref{graphs}(iii). 
In this configuration, since no edge in $\Gamma(\widehat{V})$ is a wave, 
no subarc of $w_1$ and $w_2$ is a wave with respect to $\{v_1, v_2\}$. This is Case (1) of the lemma. 

If $b\ne 0$ and $d\ne 0$ in Figure~\ref{graphs}(iii), then $\Gamma(\widehat{V})$ is as in the right 
picture of Figure~\ref{dual2}, in which case $\Sigma\ssm(v_1\cup v_2\cup w_1\cup w_2)$ consists of four 
hexagons and a collection of quadrilaterals. By the discussion of $\Gamma(\widehat{W})$ above, 
this means that $\Gamma(\widehat{W})$ must be as in the left picture of Figure~\ref{dual2}. Moreover, 
since $b\ne 0$ in Figure~\ref{graphs}(iii), $b'\ne 0$ in Figure~\ref{dual2} and hence a subarc of 
$w_1$ or $w_2$ (i.e.~edge marked $b'$ in Figure~\ref{dual2}) is a wave with respect to $\{v_1, v_2\}$. 
This is Case (2) of the lemma. 
\end{proof}



\chapter{The planar surface and the curves $\alpha$, $\gamma$, 
$\delta$ and $\varepsilon$}\label{cpt:PlanarSurfacesGenusTwo}

\section{Fixing the decomposition of $\widehat{\Sigma}$} \label{sec:NoCircularArcs}

	
As in Section~\ref{sec:complexity}, we have a genus two Heegaard splitting $(U,W)$ of $M \ssm \NN(K)$ along the Heegaard surface $\Sigma$. The Heegaard splitting extends to a Heegaard splitting $(V,W)$ of $M$ with $V=U\cup\NN(K)$.  
Given a $(\mathcal{P},\mathcal{D})$-pair $(P,D)$ as in Section~\ref{sec:complexity}, the Dehn filling along the 
slope of $\partial_-P$ yields a lens space $L$.  The union of $U$ and the surgery solid torus is a genus two 
handlebody  $V'$ in $L$, and $(V',W)$ is a genus two Heegaard splitting of the lens space $L$.  Moreover 
$(P,D)$ extends to a stabilizing pair of disks for the Heegaard splitting $(V',W)$ of $L$.  This implies that there is a
nonseparating disk $E$ in $W$ disjoint from $P\cup D$.  So $\{D, E\} $ form a complete meridian set for $W$.  
Let $\delta=\partial D$ and $\varepsilon=\partial E$.  Hence $\varepsilon\cap\partial_+P=\emptyset$ and 
$\delta\cap\partial_+P$ is a single point. In Section~\ref{sec:complexity}, we denote by $\gamma$ the curve 
which bounds the nonseparating disk $C \subset U$ and by  $\alpha=\partial_+A$ where $A$ is a vertical 
annulus and the slope of $\partial_-A$ is the $\infty$-slope.  Thus both $\gamma$ and $\alpha$ bound disks 
in the handlebody $V$ in $M$ and $\{\gamma,\alpha\}$ is a complete set of meridians for $V$. 
In this paper, we will study the Heegaard diagram of $M$ formed by $\widehat{V}=\{\gamma,\alpha\}$ 
and $\widehat{W}=\{\delta,\varepsilon\}$.

\begin{figure}[ht]
\begin{overpic}[width=5cm]{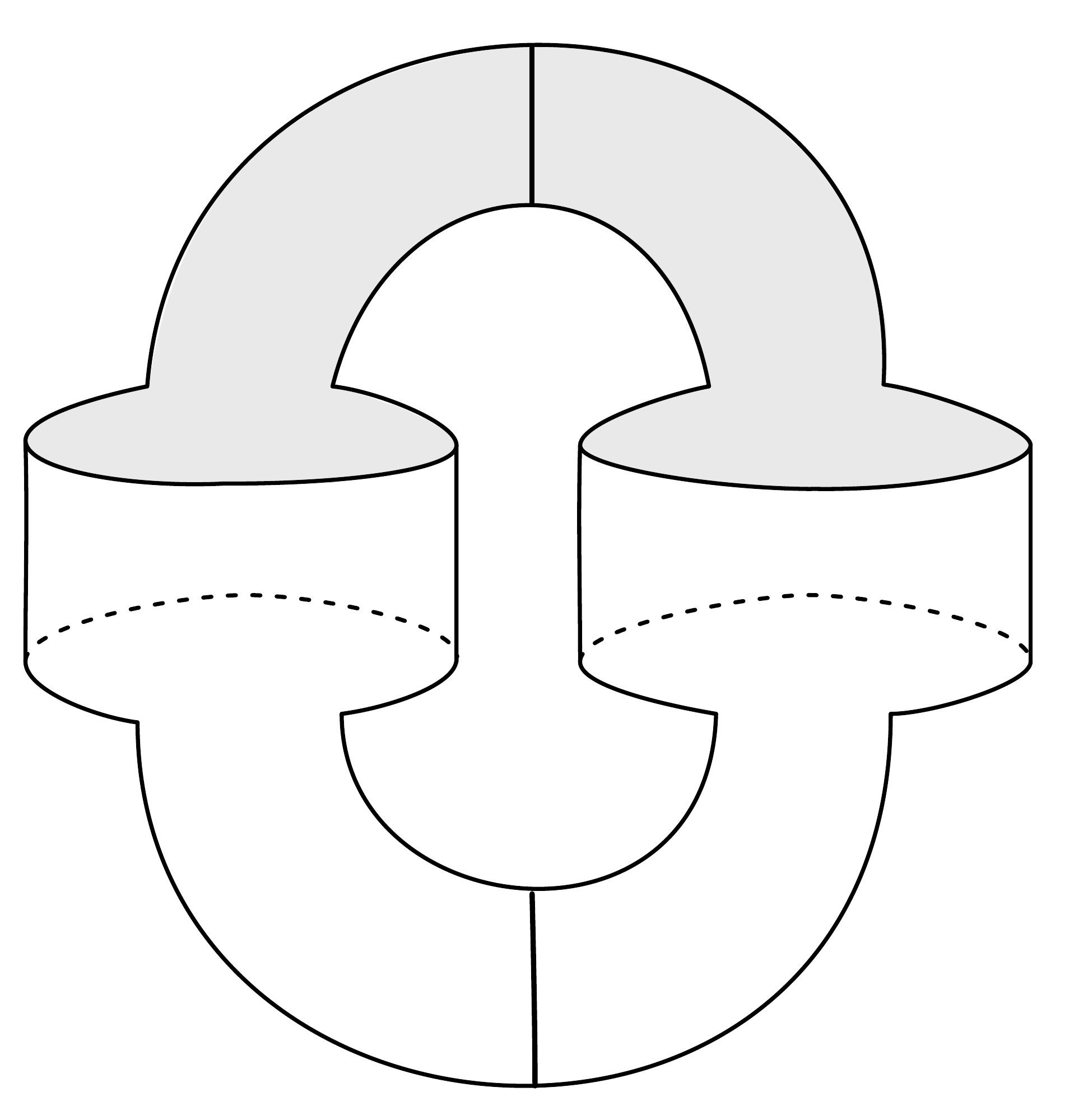}
\put(75,85){$\partial_+P^u$}
\put(74.5,10){$\partial_+P^d$}
\put(42,90){$\gamma$}
\put(42,10){$\gamma$}
\put(25,80){$\mathcal{R}^u$}
\put(25,15){$\mathcal{R}^d$}
\put(-7,50){$\mathcal{A}_l$}
\put(95,50){$\mathcal{A}_r$}
 \end{overpic}
\caption{The decomposition of $\widehat{\Sigma}$ into two annuli and two rectangles}\label{circular}
\end{figure}

Let $\widehat{\Sigma}$ be the surface obtained by cutting $\Sigma$ open along $\partial_+P$, 
i.e., $\widehat{\Sigma} = \Sigma \ssm \NN(\partial_+P)$.   Throughout the paper, we often 
view $\widehat{\Sigma}$ as the closure of $\Sigma\setminus\partial_+P$ under path metric. 
As $\partial_+P$ is nonseparating,  $\widehat{\Sigma}$ is a twice-punctured torus. 
Fix a rectangle-annulus decomposition 
of this twice-punctured  torus as a union of two annuli $\mathcal{A}_r$  and  $\mathcal{A}_l$ connected  
by two rectangles $\mathcal{R}^u$ and $\mathcal{R}^d$, as in Figure \ref{circular} 
(see Lemma~\ref{lem:GammaRectangles} below for detailed requirements). The two boundary 
components of this twice punctured torus are denoted by $\partial_+P^u$ and  $\partial_+P^d$  
(the upper and lower components as depicted in Figure~\ref{circular}). The core curves of the 
annuli are denoted by $\mathfrak{a}_r$ and   $\mathfrak{a}_l$ respectively. 

Arcs in the rectangles $\mathcal{R}^u$ and $\mathcal{R}^d$ that are properly embedded and 
non-$\partial$-parallel in $\widehat{\Sigma}$ will be called {\it cocore arcs in the rectangles} 
(e.g.~the arcs marked $\gamma$ in Figure~\ref{circular} represent components of $\gamma\cap\widehat{\Sigma}$ that are cocore arcs of rectangles $\mathcal{R}^u$ and $\mathcal{R}^d$).  An arc in $\mathcal{R}^u$ or 
$\mathcal{R}^d$ connecting $\mathcal{A}_l$ to $\mathcal{A}_r$ will be called a {\it core arc} of 
the rectangles.  An arc properly embedded in $\mathcal{A}_l$ or $\mathcal{A}_r$ connecting
$\partial_+P^u$ to  $\partial_+P^d$ will be called a cocore of the annulus.  

\vskip5pt

 In the remaining part of Section \ref{sec:NoCircularArcs} and Section \ref{sec:Delta} we prove, 
 using the above definition of complexity $c(P, D, \alpha,\gamma)$, a series of technical lemmas which are dedicated to a 
 careful analysis of possible intersections of the curves $\alpha, \gamma, \delta$ and $ \varepsilon$ 
 with the surface $\widehat{\Sigma}$.

\begin{lemma}\label{lem:GammaRectangles}
The rectangle-annulus decomposition of $\widehat{\Sigma}$ in Figure \ref{circular} can be constructed 
so that each segment of $(\alpha\cup\gamma)\cap\widehat{\Sigma}$ is a cocore arc of a rectangle 
or an annulus.  Furthermore, in each of the rectangles $\mathcal{R}^u$ and $\mathcal{R}^d$, there 
is at least one  $\gamma$ segment which is a cocore arc for the rectangle.
\end{lemma} 

\begin{proof}
Consider the disk $C\subset U$ ($\gamma=\partial C$). Assume $|P\cap C|$ is minimal up to isotopy. 
Since $P$ is incompressible,  $P\cap C$ contains no closed curve. Since $P$ is not an annulus by assumption, $P \cap C \neq \emptyset$ by 
Lemma ~\ref{lem:PCIntersection}.  Let  $\rho$ be an arc in $P \cap C$ that 
is outermost in $C$.  The arc $\rho$ cuts off a subdisk $\Delta_{\rho}$ of $C$ with $\Delta_{\rho}\cap P=\rho$.  
As $|P\cap C|$ is minimal, $\Delta_{\rho}$ is a $\partial$-compressing disk for $P$.  
Let $\kappa_1=\partial\Delta_\rho\ssm \Int(\rho)$.  So $\kappa_1\subset\gamma=\partial C$.
As $M$ is orientable and since the two endpoints of $\kappa_1$ are connected by the arc $\rho$ in $C$, 
$\partial\kappa_1$ is a pair of  points of $\gamma\cap\partial_+P$ with opposite signs of intersection. 
Thus, we may view $\kappa_1$ as a properly embedded arc in $\widehat{\Sigma}$ with both endpoints of 
$\kappa_1$ on the same boundary curve, say $\partial_+P^d$.  As $\gamma$ is a closed curve we are 
guaranteed to have another arc $\kappa_2$  of  $\gamma \cap\widehat{\Sigma}$ in $\widehat{\Sigma}$ on 
the other side of $P$ with both endpoints in $\partial_+P^u$ (one can also see the existence of $\kappa_2$ 
using the involution $\pi$).

For each $i=1,2$, consider a regular neighborhood $\NN(\kappa_i)$ of $\kappa_i$ in 
$\widehat{\Sigma}$.  We isotope all the arcs of $(\alpha\cup\gamma)\cap\widehat{\Sigma}$ that are 
parallel to $\kappa_i$ into $\NN(\kappa_i)$ ($i=1,2$), and set $\mathcal{R}^d=\NN(\kappa_1)$ and 
$\mathcal{R}^u=\NN(\kappa_2)$.  We may set all the arcs of $(\alpha\cup\gamma)\cap\widehat{\Sigma}$ in 
$\mathcal{R}^d$ and $\mathcal{R}^u$ as cocore arcs of the rectangles.

As $\NN(\kappa_1)$ and $\NN(\kappa_2)$ are disjoint in the twice-punctured torus $\widehat{\Sigma}$, 
the complement of $\mathcal{R}^d \cup \mathcal{R}^u$ is a pair of annuli, which we denote by 
$\mathcal{A}_l$ and $\mathcal{A}_r$.  Moreover, each boundary curve of $\mathcal{A}_l$ and $\mathcal{A}_r$ 
contains exactly one subarc of $\partial\mathcal{R}^d$ and $\partial\mathcal{R}^u$, as depicted in 
Figure~\ref{circular}.  Since all the arcs of $(\alpha\cup\gamma)\cap\widehat{\Sigma}$ that are 
parallel to $\kappa_i$ are in $\mathcal{R}^d \cup \mathcal{R}^u$, this implies that the remaining 
arcs of $(\alpha\cup\gamma)\cap\widehat{\Sigma}$ are nonintersecting arcs in $\mathcal{A}_l$ and $\mathcal{A}_r$ 
with endpoints in different components of $\partial\mathcal{A}_l$ and $\partial\mathcal{A}_r$.  Thus we 
may set all the arcs of $(\alpha\cup\gamma)\cap\widehat{\Sigma}$ in  $\mathcal{A}_l$ and $\mathcal{A}_r$ 
as cocore arcs of the annuli.  This gives a required decomposition.
\end{proof}

\begin{remark}\label{rem:core}
As in the proof of Lemma~\ref{lem:GammaRectangles}, if we boundary compress $P$ along $\Delta_{\rho}$, 
we obtain two planar surfaces $P_l$ and $P_r$.  By the construction in the proof, $\partial_+P_l$ and 
$\partial_+P_r$ are isotopic to the two core curves $\frak{a}_l$ and $\frak{a}_r$ of the annuli $\mathcal{A}_l$ 
and $\mathcal{A}_r$ respectively.  In particular $\partial_+P_l$ and $\partial_+P_r$ are  nonseparating 
curves in $\Sigma$.  If $P_l$ is a disk, then since $C$ is the only nonseparating essential disk in the compression 
body $U$ (up to isotopy), $P_l$ must be parallel to $C$, which contradicts the assumption that $\Delta_{\rho}$ 
is a $\partial$-compressing disk for $P$.  Thus neither $P_l$ nor $P_r$ is a disk. 
\end{remark}

\begin{remark}\label{rem: epsilonanddelta}\hfill

\noindent (1) Note that as $\varepsilon \cap \partial_{+}P = \emptyset$ the curve  $\varepsilon$ is 
contained in the  interior  of $\widehat{\Sigma}$. 
 \vskip5pt
\noindent (2) Since the hyperelliptic involution $\pi$ leaves each simple closed curve in $\Sigma$ 
 invariant up to isotopy,  we may suppose $\delta$ and $\partial_+P$ are invariant under $\pi$.  Thus 
 $X=\delta\cap\partial_+P$ is a fixed point of $\pi$.  Moreover we may suppose the two annuli $\mathcal{A}_l$ 
 and $\mathcal{A}_r$ are invariant under $\pi$ and hence  $\pi$ interchanges the two rectangles $\mathcal{R}^d$ 
 and $\mathcal{R}^u$.

\end{remark}

\begin{lemma}\label{lem:TwoDeltaArcs}
Each of the intersections $\delta\cap  \mathcal{R}^u$ and $\delta\cap\mathcal{R}^d$ contains at least 
two arcs that are core arcs of the corresponding rectangle.
\end{lemma}

\begin{proof} 
Without loss of generality, suppose that at most one component of $\delta\cap  \mathcal{R}^d$ is a 
core arc of $\mathcal{R}^d$. Note that a component of $\delta\cap  \mathcal{R}^d$ is either a core 
arc of $\mathcal{R}^d$ or an arc connecting a boundary edge of $\mathcal{R}^d$ to  the intersection 
point $X=\delta\cap\partial_+P$.  This means that, after an isotopy in $\mathcal{R}^d$, a cocore arc 
of $\mathcal{R}^d$ intersects $\delta$ at most once.  In the proof of Lemma~\ref{lem:GammaRectangles}, 
the cocore arc $\kappa_1$ of $\mathcal{R}^d$ is a boundary arc of a $\partial$-compressing disk $\Delta_\rho$ 
for $P$. So we may isotope $\kappa_1$ in $\mathcal{R}^d$ so  that $\kappa_1$ intersects $\delta$ at most once.  
Since $P\cap\delta=X$ is a single point, the assumption on $\kappa_1\cap\delta$ implies that the surface 
obtained by the $\partial$-compression along $\Delta_\rho$ intersects $\delta $ in at most three points: the 
point $X$ and the two copies of  $\kappa_1\cap\delta$ created by the boundary compression.

As in Remark~\ref{rem:core}, the two planar surfaces obtained by $\partial$-compressing $P$ 
along $\Delta_\rho$ are $P_l$ and $P_r$.  So $(P_l\cup P_r)\cap\delta$ consists of at most 3 points. Moreover,  $(P_l\cup P_r)\cap\delta$ consists of exactly 3 points if and only if  $\kappa_1\cap\delta\ne\emptyset$, in which case $P_l\cap\delta\ne\emptyset$ and $P_r\cap\delta\ne\emptyset$.
This implies that either $P_l\cap\delta$ or $P_r\cap\delta$ is a single point, say, $P_l\cap\delta$. Thus $(P_l, D)$ 
forms a new $(\mathcal{P},\mathcal{D})$-pair. 

By Remark~\ref{rem:core}, $\partial_+P_l$ is isotopic to  the core curve $\frak{a}_l$ of the annulus 
$\mathcal{A}_l$.  We may isotope $\partial_+P_l$ to the upper boundary component of $\mathcal{A}_l$. 
After this isotopy,  any intersection point of $\alpha\cup\gamma$ with $\partial_+ P_l$ is an 
intersection point with $\partial_+P$.  Therefore, any segment of $(\alpha\cup\gamma)\ssm\partial_+ P_l$ 
is composed of a nonzero number of segments of  $\alpha\cup\gamma$ in $\widehat{\Sigma}$.  It follows that  
$c(P_l,D,\alpha,\gamma) < c(P,D,\alpha,\gamma)$, contradicting Assumption~\ref{ass:minimality}.
\end{proof}

\begin{lemma}\label{lem:FirstWave}
Let $P$, $D$ and $\widehat{V}=\{\gamma,\alpha\}$  be as above.  Then there is no $s$-wave with 
respect to $\alpha$ that is disjoint from $\delta$.
\end{lemma}
\begin{proof}
Suppose, in contradiction, that there is such an  $s$-wave $\eta$ with  $\partial\eta\subset\alpha$ 
and $\eta\cap\delta=\emptyset$.  

Let $\alpha_1$ and $\alpha_2$ be the two components of $\alpha\ssm\NN(\partial\eta)$.  We define 
$c_0(\alpha_i)$ ($i=1,2$) to be the number of segments of $\alpha_i\ssm\partial_+P$ that intersect $\delta$.  
This implies that $c_0(P, D,\alpha)+2\ge c_0(\alpha_1)+c_0(\alpha_2)$ (see Definition~\ref{def:complexity}) 
because a segment of $\alpha\ssm\partial_+P$ may be split into two segments by $\partial\eta$,

Recall that the first step of a wave move along $\eta$  is to connect the two endpoints of $\alpha_1$ 
and the two endpoints of $\alpha_2$ by a pair of arcs parallel to $\eta$.  Denote the two resulting 
closed curves  by $\alpha_1'$ and $\alpha_2'$ respectively.   Since $\eta\cap\delta=\emptyset$ 
for $i \in \{1,2\}$ we have   $c_0(\alpha_i)\ge c_0(P,D,\alpha_i')$. Therefore, 
$c_0(P, D,\alpha)+2\ge c_0(P,D,\alpha_1')+c_0(P,D,\alpha_2')$.

Because $\Sigma$ is of genus two, one of the two closed curves $\alpha_1'$ and $\alpha_2'$, obtained after 
the first step of the wave move, is isotopic to $\gamma$. Say it is $\alpha_2'$. By Lemma~\ref{lem:TwoDeltaArcs} 
and the construction of $\mathcal{R}^u$ and $\mathcal{R}^d$, $c_0(P,D,\alpha_2') = c_0(P,D,\gamma)\ge 2$.  
This plus the inequality above implies that $c_0(P, D,\alpha)\ge c_0(P,D,\alpha_1')$.  Recall that in the case 
of a genus two surface, a curve obtained by a wave move can also be obtained by a band sum. So $\alpha_1'$ 
can be obtained by a band sum of $\alpha$ and $\gamma$.  Since $\alpha=\partial_+A$ and $\gamma=\partial C$,
$\alpha_1'$ is a boundary curve of an annulus obtained by a band sum of $A$ and $C$.  Moreover,
$|\,\alpha\cap\delta\,|=|\,\alpha_1'\cap\delta|\,+|\,\alpha_2'\cap\delta\,|$, and by Lemma~\ref{lem:TwoDeltaArcs}, 
$\gamma\cap\delta = \alpha_2'\cap\delta \ne\emptyset$. Thus $|\alpha_1'\cap\delta| < |\alpha\cap\delta|$ 
and  $c(P,D,\alpha_1',\gamma) < c(P,D,\alpha,\gamma)$, contradicting Assumption~\ref{ass:minimality}.
\end{proof}

\begin{lemma} \label{lem:NoAWaves} 
Any $s$-wave with respect to $\widehat{V}=\{\gamma,\alpha\}$ must intersect $\partial_{+}P$.
\end{lemma}

\begin{proof}
As in the proof of Lemma~\ref{lem:PCIntersection}, since $P$ is incompressible and $U\ssm\NN(C)\cong T^2\times I$, 
each component of $P\ssm\NN(C)$ is either a disk or a vertical annulus in $U\ssm\NN(C)$.  Let $\rho'$ be an  
arc in the intersection $P\cap C$ which is outermost in $P$.  So $\rho'$ cuts off a subsurface $Q'$ of $P$ 
with $Q'\cap C=\rho'$.  Since $|P\cap C|$ is assumed to be minimal up to isotopy, $Q'$ cannot be a disk. 
Hence $Q'$ is a vertical  annulus in  $U\ssm\NN(C)$.  Now consider the annulus $A$ with $\partial_+A=\alpha$.  
As $C\cap A=\emptyset$ and since $\partial_-P$ has integer slope, after isotopy, $Q'\cap A$ is a single vertical arc,  
and thus $\partial_+Q' \cap  \alpha$ is a single point.  Note that $\partial_+Q'$ is the union of $\rho'$ and a subarc 
$\kappa'$ of $\partial_+P$, see Figure~\ref{fig:genus2}. So $\kappa'$ is an arc  properly embedded 
in $\Sigma\ssm\NN(\gamma)$ intersecting $\alpha$ exactly once and having both endpoints in the 
same boundary component of $\Sigma\ssm\NN(\gamma)$. 

\begin{figure}[ht]
	\begin{overpic}[width=6cm]{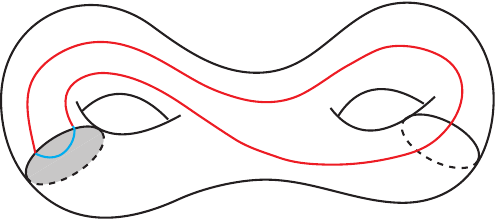}
		\put(15.5,13.5){{\tiny $C$}}
		\put(20,9){$\gamma$}
		\put(91,12){$\alpha$}
		\put(8.5,9){{\tiny $\rho'$}}
		\put(55,18){$\kappa'$}
	\end{overpic}
	\caption{Subarcs of $\partial_+P$ forming blocking-edges for $s$-wave of $\{\gamma,\alpha\}$}
	\label{fig:genus2}
\end{figure}

Now cut $\Sigma$ open along $\gamma$ and $\alpha$ and obtain a four-holed sphere. Denote the boundary 
curves of the four-holed sphere by $\gamma^+$, $\gamma^-$, $\alpha^+$ and $\alpha^-$. The restriction of 
$\partial_+P$ to this four-holed sphere is a collection of arcs. It follows that the arc $\kappa'$ is cut into two subarcs 
in the four-holed sphere connecting $\alpha^+$ and $\alpha^-$ to the same curve, say $\gamma^+$.  By the 
symmetry induced from the involution $\pi$,  the curve $\partial_+P$ also contains subarcs in the four-holed 
sphere connecting both $\alpha^+$ and $\alpha^-$ to $\gamma^-$.   Similar to the proof of cases (2) and (3) 
of Lemma~\ref{lem:+ and -}, these two subarcs of $\partial_+P$ ``block'' any $s$-wave with respect to 
$\widehat{V}=\{\gamma,\alpha\}$.   In other words,  any $s$-wave with respect to $\widehat{V}=\{\gamma,\alpha\}$ 
must intersect $\partial_{+}P$. 
\end{proof}

\begin{corollary}\label{cor:NoWaves}
The curves $\delta$ and $\varepsilon$ do not contain subarcs which are  waves with respect to  
$\alpha$ or $\gamma$. 
\end{corollary}

\begin{proof}
As above we denote the boundary curves of the four-holed sphere $\Sigma\ssm \NN(\alpha \cup \gamma)$ 
by $\gamma^+$, $\gamma^-$, $\alpha^+$ and $\alpha^-$.  If  a subarc $\eta$ of $\delta\cup\varepsilon$ 
is a wave  with $\partial\eta$ in $\gamma^\pm$ (or $\alpha^\pm$), then by the symmetry from the hyperelliptic 
involution, another subarc $\eta'$ of $\delta\cup\varepsilon$ is a wave with $\partial\eta'$ in $\gamma^\mp$ 
(or $\alpha^\mp$).  By Lemma~\ref{lem:NoAWaves}, both $\eta$ and $\eta'$ must intersect $\partial_+P$, which 
means that $\delta\cup\varepsilon$ has at least two intersection points with $\partial_+P$.  This is a contradiction 
because $\varepsilon\cap P=\emptyset$ and $|\delta\cap P|=1$.
\end{proof}

Let  $\{\delta,\varepsilon'\}$ be a complete set of meridians for $W$, where $\varepsilon'$ may 
not be the same as $\varepsilon$. Fix an orientation for each of $\alpha$, $\gamma$, $\delta$ and 
$\varepsilon'$. These orientations induce a $\pm$-sign for each point of intersection of the curves.

\begin{proposition}\label{pro:Paths}
Let $\alpha$, $\gamma$, $\delta$, $\varepsilon'$ be as above. Consider the Heegaard diagram given by  
$\{\delta, \varepsilon'\}$ and $\{\alpha, \gamma\}$.  Suppose the Heegaard diagram is not a standard 
Heegaard diagram of $S^3$ or  $(S^2\times S^1)\# L(r, s)$. Then
\begin{enumerate}
\item  Suppose there is no wave with respect to $\{\delta, \varepsilon'\}$, then all  the intersection points, of $\delta$ 
with each curve of $\{\alpha, \gamma\}$,  have the same sign.  Moreover, if no subarc of $\varepsilon'$ is a wave 
with respect to  $\{\alpha, \gamma\}$, then the intersection points of  $\varepsilon'$ with each curve of 
$\{\alpha, \gamma\}$ all have the same sign. 
\vskip4pt
\item If $\varepsilon' = \varepsilon$, then 
$\delta$ and $\varepsilon$ admit orientations so that the intersection points of 
$\delta\cup\varepsilon$ with each component of $(\alpha\cup\gamma)\ssm \partial_+P$ all have the same sign.
\end{enumerate}
\end{proposition}

\begin{proof}  
Assume in contradiction that the statement of part (1) of the Proposition is false and that $\delta$ intersects 
a curve $v_1\in\{\alpha,\gamma\}$ with different signs. In particular, $\delta$ has a subarc connecting 
one side of $v_1$, say the plus side, to the same side.  Denote an innermost such subarc by $\tau$, so 
$\tau\cap v_1=\partial\tau$ and $\tau$ connects the plus side of $v_1$ to the same side. If $\tau$ does not 
have other intersection points with the other curve $v_2\in\{\alpha,\gamma\}$ in its interior, then $\tau$ is 
a wave with respect to $v_1\in\{\alpha,\gamma\}$, contradicting Corollary~\ref{cor:NoWaves}. Hence $\tau$ 
must intersect  the other curve $v_2\in\{\alpha,\gamma\}$  in its interior. These intersection points of $\tau$ 
with $v_2$ must have the same sign (otherwise $\tau$ contains a subarc that is a wave with respect to 
$v_2\in\{\alpha,\gamma\}$  contradicting Corollary~\ref{cor:NoWaves}). 

Denote the two sides of $v_i$ by $v_i^\pm$. The arc $\tau$ was chosen to connect $v_1^+$ to $v_1^+$. 
The curve $v_2$ divides $\tau$ into a collection of subarcs. Since all the points of $\tau\cap v_2$ have the 
same sign, the two subarcs of $\tau$ containing $\partial\tau$ must be arcs connecting $v_1^+$ to $v_2^+$ 
and $v_1^+$ to  $v_2^-$ respectively, in other words, the two arcs are $[v_1^+,v_2^+]$ and $[v_1^+, v_2^-]$ 
blocking-edges. Hence, by Lemma~\ref{lem:+ and -}, there can be no wave with respect to 
$\widehat{V}=\{\alpha,\gamma\}$.  Since, by the hypothesis, there is no wave with respect to $\{\delta, \varepsilon'\}$, 
there is no wave in the Heegaard diagram given by $\{\delta, \varepsilon'\}$ and $\{\alpha, \gamma\}$, 
contradicting Theorems ~\ref{thm:HOT} and ~\ref{thm:NeOk}.  Thus all  the intersection points, of $\delta$ with 
each curve of $\{\alpha, \gamma\}$, must have the same sign. 

If no subarc of $\varepsilon'$ is a wave with respect to  $\{\alpha,\gamma\}$, then the same proof also 
works for $\varepsilon'$ and hence the intersection points of  $\varepsilon'$ with each curve of 
$\{\alpha, \gamma\}$ all have the same sign.  This proves part (1) of the proposition.
 
To prove part (2), we first show that the intersection points of each curve of $\{\delta, \varepsilon\}$ 
with any component of $(\alpha\cup\gamma)\ssm \partial_+P$ all have the same sign:  Suppose, 
to the contrary,  that a curve $w_1\in\{\delta, \varepsilon\}$ has opposite signs at intersection points  
with a segment  $h$ of $(\alpha\cup\gamma)\ssm \partial_+P$. Then a subarc of $h$ connects  
 one side of $w_1$, say the plus side, to the same side.  Let $\tau'$ be an innermost such subarc.  
 So $\tau'\cap w_1=\partial\tau'$.   
 
 First suppose that $\tau'$ does not intersect the other curve 
 $w_2\in\{\delta, \varepsilon\}$ in its interior, then $\tau'$ is a wave with respect to  
 $w_1\in\{\delta, \varepsilon\}$ and is disjoint from  $\partial_+P$.   
If $w_1=\delta$, since $\tau'\cap\partial_+P=\emptyset$, a wave move on $\delta$ along $\tau'$ yields a new 
meridional curve $\delta'$ of $W$ that intersects $\partial_+P$ in a single point. Let $D'$ be the disk in $W$ 
bounded by $\delta'$. So $(P, D')$ is a $(\mathcal{P},\mathcal{D})$-pair.  However, it is clear from the wave 
move that the complexity $c(P,D',\alpha,\gamma) < c(P,D,\alpha,\gamma)$. This contradicts 
Assumption~\ref{ass:minimality}.  If $w_1=\varepsilon$ then a wave move on $\varepsilon$ along $\tau'$ 
yields a meridional curve of $W$ that is disjoint from $\delta\cup\partial_+P$ and not parallel to $\varepsilon$, 
which is impossible in the genus two surface $\Sigma$.   

Thus $\tau'$ cannot be a wave and hence $\tau'$ must contain intersection points with the other curve 
$w_2\in\{\delta, \varepsilon\}$  in its interior.  By applying the argument used on the arc $\tau$ as above to $\tau'$,  
 we conclude that the arc $\tau'$ contains two subarcs that are $[w_1^+, w_2^+]$ and $[w_1^+, w_2^-]$ 
 blocking-edges. Thus, by Lemma~\ref{lem:+ and -},  there can be no wave with respect to $\{\delta, \varepsilon\}$.  
 By Corollary~\ref{cor:NoWaves}, no subarc of $\varepsilon$ is a wave with respect to $\{\alpha,\gamma\}$.  Hence, 
 it follows from part (1) of this proposition that the intersection points of $w_1$ with each curve of $\{\alpha,\gamma\}$ 
 all have the same sign, a contradiction to our assumption on $h$ at the beginning. 
 This proves that the intersection points of each curve of $\{\delta, \varepsilon\}$ 
 with any component of $(\alpha\cup\gamma)\ssm \partial_+P$ all have the same sign. 

This conclusion implies that part (2) of the proposition holds if no component of  $(\alpha\cup\gamma)\ssm \partial_+P$ 
 intersects both $\delta$ and $\varepsilon$.  Assume therefore, in contradiction, that $h$ is a component of  
 $(\alpha\cup\gamma)\ssm \partial_+P$ which intersects both $\delta$ and $\varepsilon$.  By the conclusion 
 above, we may fix orientations for $\delta$ and $\varepsilon$ so that the intersection points of $\delta\cup\varepsilon$ 
 with $h$ all have the same sign.  So part (2) of the proposition is true unless there is another component $h'$ of 
 $(\alpha\cup\gamma)\ssm \partial_+P$ such that points in $h'\cap\delta$ and  $h'\cap\varepsilon$ have opposite signs. 
 Suppose there is such an arc $h'$. A subarc of $h'$ has one endpoint in $\delta$ and the other endpoint in $\varepsilon$.  
 As they have opposite signs, this subarc is either a $[\delta^+,\varepsilon^+]$ or a $[\delta^-,\varepsilon^-]$ edge.  
 Similarly, since the intersection points of $\delta\cup\varepsilon$ with $h$ all have the same sign, a subarc of $h$ 
 is a $[\delta^-,\varepsilon^+]$ or $[\delta^+,\varepsilon^-]$ edge. Thus $h\cup h'$ contains two subarcs that either 
connect both $\delta^+$ and $\delta^-$ to the same $\varepsilon^\pm$ or connect both $\varepsilon^+$ and 
$\varepsilon^-$ to the same $\delta^\pm$.  It now follows from (2) and (3) of Lemma~\ref{lem:+ and -} that 
the Heegaard diagram contains no wave with respect to $\{\delta, \varepsilon\}$. 

By Corollary~\ref{cor:NoWaves} no subarc of $\varepsilon$ is a wave.  So
 part (1) of this proposition implies that $h$ and $h'$ must belong to different curves of $\{\alpha,\gamma\}$.  
Without loss of generality,  suppose $h\subset\alpha$, $h'\subset\gamma$. Moreover, by part (1) of the 
proposition, we may choose the orientation of $\delta$ and $\varepsilon$ so that: 
\begin{enumerate}[(a)]
\item points in $\alpha\cap(\delta\cup\varepsilon)$ all have positive signs, 

\item points in $\gamma\cap\delta$ have positive signs, and

\item points in $\gamma\cap\varepsilon$ have negative signs.
\end{enumerate}  

Similar to the argument on $h$ and $h'$ above, this means that a subarc of $\delta$ is an $[\alpha^+,\gamma^-]$ 
or $[\alpha^-,\gamma^+]$ edge and a subarc of $\varepsilon$ is an  $[\alpha^+,\gamma^+]$ or $[\alpha^-,\gamma^-]$ edge.  
Again by (2) and (3) of Lemma~\ref{lem:+ and -}, the Heegaard diagram contains no wave with respect to $\{\alpha, \gamma\}$.   
 Hence the Heegaard diagram has no wave, contradicting Theorems ~\ref{thm:HOT} and ~\ref{thm:NeOk}.
\end{proof}

\begin{lemma}\label{lem:Compatibility}  
Let $\alpha$, $\gamma$, $\delta$, $\varepsilon'$ be as in Proposition~\ref{pro:Paths}.  
Consider the Heegaard diagram formed by $\widehat{V} = \{\gamma,\alpha\}$ and 
$\widehat{W} = \{\delta,\varepsilon'\}$. Suppose the Heegaard diagram is not a standard 
Heegaard diagram of $S^3$ or  $(S^2\times S^1)\# L(r, s)$.  Suppose that $\alpha\cup\gamma$ 
contains both $[\delta^+ ,\delta^-]$ and  $[\varepsilon'^+ ,\varepsilon'^-]$ blocking-edges. Furthermore,  
suppose that no subarc of $\varepsilon'$ is a wave with respect to $\widehat{V}=\{\alpha,\gamma\}$. Then 
\begin{enumerate}
\item $\Gamma(\widehat{W})$ is of type (i) in Figure~\ref{graphs}, where $c\ne 0$, $d\ne 0$ and
 exactly one of $a$, $b$ is $0$. In other words, $\Gamma(\widehat{W})$ is as in the left 
 picture of Figure~\ref{dual1}. See Figure~\ref{fig:octagon}(a) for an example of such Heegaard diagram. 
\vskip5pt
\item Suppose that there are more than one $[\delta^+ ,\delta^-]$-edges in the Heegaard diagram. Then 
there is a wave with respect to $\{\alpha,\gamma\}$ connecting a $[\delta^+ ,\delta^-]$-edge to an
$[\varepsilon'^+ ,\varepsilon'^-]$-edge, see the dashed arc in Figure~\ref{fig:octagon}(b). 
\end{enumerate}
\end{lemma}

\begin{proof} 
First note that since $\Gamma(\widehat{W})$ contains both  $[\delta^+ ,\delta^-]$ and  $[\varepsilon'^+ ,\varepsilon'^-]$ 
blocking-edges, it follows from case (1) of Lemma~\ref{lem:+ and -} that there is no wave with respect to 
$\{\delta,\varepsilon'\}$. Moreover, it implies  that $\Gamma(\widehat{W})$ must be of type (i) in Figure~\ref{graphs} 
with $c\ne 0$ and $d\ne 0$. 

If both $a$ and $b$ are $0$, then the Whitehead graph is disconnected and there is a circle $\Psi\subset\Sigma$ 
which is disjoint from $\widehat{V}\cup\widehat{W}$ and which separates the $[\delta^+ ,\delta^-]$-edges 
from the $[\varepsilon'^+ ,\varepsilon'^-]$-edges.  Since $\Psi$ is a separating circle in $\Sigma$ disjoint 
from $\widehat{V}\cup\widehat{W}$, $\Psi$ bounds separating disks in both 
handlebodies $V$ and $W$ and the two disks form a separating $S^2$ in $M$. This implies that 
$\Gamma(\widehat{W})$ is  a standard Heegaard diagram of a connected sum of two lens spaces $L(c,s) \# L(d,t)$.
Since $c>0$ and $d>0$, none of $L(c,s)$ and $L(d,t)$ is $S^2\times S^1$.  As $\Gamma(\widehat{W})$ is not the
standard Heegaard diagram of $S^3$ and since $M$ is either $S^3$ or a connected sum of $S^2\times S^1$ and 
a lens space, this is a contradiction. 

The argument above implies that the Whitehead graph is connected. Moreover, since there is no wave 
with respect to $\{\delta,\varepsilon'\}$, by Theorems~\ref{thm:HOT} and \ref{thm:NeOk}, there must 
be a wave with respect to $\{\alpha,\gamma\}$. 

By Corollary~\ref{cor:NoWaves}, no subarc of $\delta$ is a wave with respect to $\{\alpha,\gamma\}$. Thus, 
by the hypothesis that no subarc of $\varepsilon'$ is a wave with respect to $\{\alpha,\gamma\}$, the conditions 
in  part (1) of Lemma~\ref{lem:dual} are satisfied. Hence, $\Gamma(\widehat{W})$ must be as in the left picture of 
Figure~\ref{dual1}. This proves part (1) of the lemma.

Next, we prove part (2) of the lemma. Since no subarc of $\varepsilon'$ is a wave, by part (1) of 
Proposition~\ref{pro:Paths}, we may assign orientations for $\delta$, $\varepsilon'$,  
and $\gamma$ so that the intersection points of $\gamma$ with $\delta\cup\varepsilon'$ all have the same sign. 
In particular, the Heegaard diagram is as in Figure~\ref{fig:octagon}(a). Furthermore, 
the configuration of Figure~\ref{fig:octagon}(a) implies that we may choose the orientation for $\alpha$ so 
that all the intersection points have the same sign. 

\begin{figure}[ht]
\begin{overpic}[width=10cm]{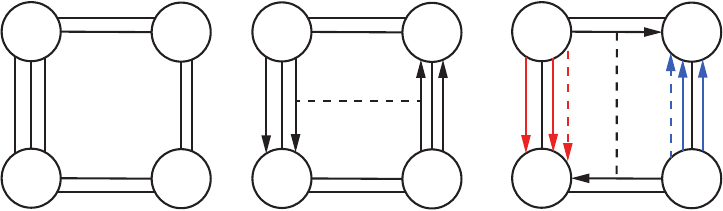}
\put(14,-4){(a)}
\put(48,-4){(b)}
\put(84,-4){(c)}
\put(48,16.5){$\eta$}
\put(86,15){$\eta$}
\put(3,3){$\delta^-$}
\put(37,3){$\delta^-$}
\put(73.5,3){$\delta^-$}
\put(23,3){$\varepsilon'^+$}
\put(58,3){$\varepsilon'^+$}
\put(94,3){$\varepsilon'^+$}
\put(3,24){$\delta^+$}
\put(37,24){$\delta^+$}
\put(73.5,24){$\delta^+$}
\put(23,24){$\varepsilon'^-$}
\put(58,24){$\varepsilon'^-$}
\put(94,24){$\varepsilon'^-$}
\end{overpic}
\vspace{10pt}
\caption{Heegaard diagrams and waves}
\label{fig:octagon}
\end{figure}

As illustrated in Figure~\ref{fig:octagon}(a), the four curves $\delta$, $\varepsilon'$, $\alpha$ and $\gamma$ 
cut the Heegaard surface into two octagons and a collections of quadrilaterals.  All the $[\delta^+ ,\delta^-]$-edges 
are parallel and all the $[\varepsilon'^+ ,\varepsilon'^-]$-edges are parallel.  We say that a $[\delta^+ ,\delta^-]$-edge 
or an $[\varepsilon'^+ ,\varepsilon'^-]$-edge is outermost if it is a boundary edge of an octagon. The hyperelliptic
involution $\pi$ interchanges the two octagons and the two outermost $[\delta^+ ,\delta^-]$-edges.  Furthermore, 
if there is more than one $[\varepsilon'^+ ,\varepsilon'^-]$-edge, then $\pi$ also interchanges the two outermost 
$[\varepsilon'^+ ,\varepsilon'^-]$-edges, and if there is a single $[\varepsilon'^+ ,\varepsilon'^-]$-edge, then $\pi$ 
leaves this edge invariant.

Since all the intersection points have the same sign, if the outermost $[\delta^+ ,\delta^-]$-edges and 
$[\varepsilon'^+ ,\varepsilon'^-]$-edges belong to the same curve $\alpha$ or $\gamma$, then the dashed arc 
$\eta$ in Figure~\ref{fig:octagon}(b) is a wave and part (2) of the lemma holds.  Next, suppose that they 
belong to different curves. Without loss of generality, suppose the outermost $[\delta^+ ,\delta^-]$-edges 
belong to $\gamma$ and the outermost $[\varepsilon'^+ ,\varepsilon'^-]$-edges belong to $\alpha$.

By Theorems~\ref{thm:HOT} and \ref{thm:NeOk}, there must be a wave $\eta$ with respect to 
$\{\alpha, \gamma\}$.  Since no subarc of $\delta$ or $\varepsilon'$ is a wave, a wave  with 
respect to $\{\alpha, \gamma\}$ must be in an octagon connecting two opposite boundary edges. 
Thus, either there is a wave as in Figure~\ref{fig:octagon}(b) and part (2) of the lemma holds, 
or the wave $\eta$ is as shown in Figure~\ref{fig:octagon}(c). 

Suppose the wave $\eta$ is as shown in Figure~\ref{fig:octagon}(c).
By Lemma~\ref{lem:FirstWave}, $\eta$ must be a wave with respect to $\gamma$.  Perform the wave 
move along $\eta$ in two steps: The first step is a surgery on $\gamma$ along $\eta$, resulting a new 
meridian $\gamma'$ and a curve $\alpha'$ parallel to $\alpha$, see the dashed arcs in 
Figure\ref{fig:octagon}(c). The second step is to delete $\alpha'$ and obtain a new set of meridians 
$\{\alpha,\gamma'\}$. 

Since $\alpha'$ is parallel to $\alpha$ and since the outermost $[\delta^+,\delta^-]$-edges before 
the wave move belong to $\gamma$, the newly created $[\varepsilon'^+ ,\varepsilon'^-]$-edge by 
the first step (see the blue dashed arc in Figure\ref{fig:octagon}(c)) must belong to $\alpha'$, and 
the newly created $[\delta^+,\delta^-]$-edge (see the red dashed arc in Figure\ref{fig:octagon}(c)) 
must belong to $\gamma'$. Moreover, this newly created $[\delta^+,\delta^-]$-edge is an outermost 
$[\delta^+,\delta^-]$-edge. As every essential curve in $\Sigma$ is invariant under the involution 
$\pi$, the other outermost $[\delta^+,\delta^-]$-edge must also belong to $\gamma'$. Since 
the second step removes $\alpha'$ and leave $\gamma'$ unchanged, after the wave move, 
we must have two outermost $[\delta^+,\delta^-]$-edges, both belong to $\gamma'$. Furthermore, 
since $\alpha$ is unchanged by this wave move, the outermost $[\varepsilon'^+ ,\varepsilon'^-]$-edges 
still belong to $\alpha$, and if there is more than one  $[\varepsilon'^+ ,\varepsilon'^-]$-edge before 
the wave move, there is more than one $[\varepsilon'^+ ,\varepsilon'^-]$-edge in the new Heegaard diagram.
Therefore, we can conclude: 

\begin{enumerate}
\item In the new Heegaard diagram after the wave move, the outermost $[\delta^+ ,\delta^-]$-edges 
belong to $\gamma'$ and the outermost $[\varepsilon'^+ ,\varepsilon'^-]$-edges still belong to $\alpha$.

\vskip5pt

\item Since there is more than one $[\delta^+ ,\delta^-]$-edge in the diagram before the wave move, 
there is more than one $[\delta^+ ,\delta^-]$-edge in new Heegaard diagram. Moreover, if there is 
more than one $[\varepsilon'^+ ,\varepsilon'^-]$-edge before the wave move, there is more than 
one $[\varepsilon'^+ ,\varepsilon'^-]$-edge in the new Heegaard diagram.
\end{enumerate}

Since there is a $[\delta^+ ,\delta^-]$-edge and an $[\varepsilon'^+ ,\varepsilon'^-]$-edge in the new 
Heegaard diagram, there is no wave with respect to $\{\delta,\varepsilon'\}$ by 
Lemma~\ref{lem:+ and -}.  Moreover, $\gamma'$ has an induced orientation from $\gamma$ and 
hence the intersection points of $\delta$ and $\varepsilon'$ with $\gamma'$ all have the same sign. 
So no subarcs of $\delta$ or $\varepsilon'$ is a wave with respect to $\{\alpha,\gamma'\}$. 
The argument above implies that a wave with respect to $\{\alpha,\gamma'\}$ must be in an 
octagon similar to the arc $\eta$ in Figure~\ref{fig:octagon}(c).  Thus we can repeat the wave 
moves described above until the new Heegaard diagram is of type (i) in Figure~\ref{graphs} with $a=b=0$, 
which means that $M$ is a connected sum of two lens spaces $L(c,s) \# L(d,t)$ and none of $L(c,s)$ and 
$L(d,t)$ is $S^2\times S^1$.  

Furthermore, the argument above implies that there is always more than one $[\delta^+ ,\delta^-]$-edge 
after each wave move, that is $c\ge 2$. Hence $M$ cannot be $S^3$.  This contradicts the 
hypothesis that  $M$ is either $S^3$ or a connected sum of $S^2\times S^1$ with a lens space.
\end{proof}

\section{ Paths, Junctions, Traintracks  and  the curves $\delta, \varepsilon$}\label{sec:Delta}
 
As  $\widehat{\Sigma}$  is endowed with the structure from Lemma ~\ref{lem:GammaRectangles},  
it has a decomposition as two annuli and two rectangles  
$\mathcal{A}_r  \cup  \mathcal{R}^u \cup \mathcal{A}_l  \cup  \mathcal{R}^d$.  As explained above, 
we assume that the hyperelliptic involution $\pi$ maps $\partial_+P$ to itself. Hence $\pi$ induces 
an involution on $\widehat{\Sigma}$ interchanging the two boundary curves $\partial_+P^u$ and 
$\partial_+P^d$.

\begin{remark}\label{rem:Symmetries} Note the following: \\
\noindent   (A) The surface, as in Figure \ref{circular}, has the following three symmetries:

\begin{enumerate}
\item  The hyperelliptic involution: a  $180^\circ$-rotation about the horizontal axis within the projection 
plane in Figure~\ref{circular},  puncturing each annulus in two points. It interchanges the rectangles, 
and  keeps each annulus in place.
\vskip4pt
\item A reflection along a vertical plane cutting through the ``middle'' of the figure. The vertical plane 
is disjoint from the two annuli and intersects each rectangle in a cocore arc.  The reflection interchanges 
the two annuli $\mathcal{A}_r$ and $ \mathcal{A}_l$. 
\vskip4pt
\item A reflection along a cylinder which meets each of $ \mathcal{R}^u$ and  $\mathcal{R}^d$ in a core 
arc and each of $\mathcal{A}_r$ and $ \mathcal{A}_l$ in two cocore arcs. The restriction of this reflection 
on each rectangle  is a reflection along its core arc, and the restriction on each annulus is a reflection of 
the annulus along a vertical plane that cut the annulus vertically into two halves.  This symmetry 
flips the orientations of the core curves of the two annuli.
\end{enumerate}
\vskip4pt
\noindent (B) The second and the third reflections commute with the hyperelliptic involution.
\end{remark}

As it is assumed that $\pi$ leaves $\delta=\partial D$ and $\partial_+P$  invariant, $\pi$  fixes 
$X = \delta\cap\partial_+P$. The image of the point $X$ in $\partial_+P^d$ and $\partial_+P^u$ 
is denoted by $X^d$  and $X^u$ respectively.  So the restriction of 
$\delta$ to $\widehat{\Sigma}$ is an arc connecting $X^d$ to $X^u$, and the restriction of $\pi$ on 
$\widehat{\Sigma}$ interchanges $X^u$ and $X^d$.   

First consider the possibility that $X^d\in\partial\mathcal{R}^d$. Hence, by the involution,
$X^u\in\partial\mathcal{R}^u$. In this case, let $k^d$ be the arc component of 
$\delta\cap \mathcal{R}^d$ that contains $X^d$.  If $k^d\cap(\alpha\cup\gamma)=\emptyset$, 
then by the construction of the rectangles, we may ``shrink'' the rectangle $\mathcal{R}^d$, by an 
isotopy so that both $k^d$ and $X^d$ are no longer in  $\mathcal{R}^d$.  Thus, we may assume that if
 $X^d\in\partial\mathcal{R}^d$, then $k^d\cap(\alpha\cup\gamma)\ne\emptyset$.  Note that we can 
 perform this  shrinking operation symmetrically on $\mathcal{R}^u$ so that  the symmetry from  
 $\pi$ is preserved.
 
\begin{lemma}\label{lem:NoCircular}
No component of $(\delta\cup\varepsilon)\cap\mathcal{A}_l$ and  $(\delta\cup\varepsilon)\cap\mathcal{A}_r$ 
is an arc with both endpoints in the same rectangle $\mathcal{R}^d$ or $\mathcal{R}^u$.
\end{lemma}

\begin{proof}
Assume to the contrary that the statement of the lemma is false. Without loss of generality, suppose 
$(\delta\cup\varepsilon)\cap\mathcal{A}_l$ has a component $\rho$ with $\partial\rho\subset\mathcal{R}^d$.  
Note that $\rho$ is a topologically $\partial$-parallel arc in the annulus $\mathcal{A}_l$.  Let $k_1$ and $k_2$ 
be the components of either $\delta\cap \mathcal{R}^d$ or $\varepsilon\cap\mathcal{R}^d$ that contain the two endpoints of $\rho$. So $k_i$ ($i=1,2$) is either a core arc of $\mathcal{R}^d$ or the component $k^d$ of $\delta\cap \mathcal{R}^d$ that contains $X^d$. 
By the assumption before Lemma~\ref{lem:NoCircular}, $k^d\cap(\alpha\cup\gamma)\ne\emptyset$. Since a core arc of $\mathcal{R}^d$ intersects every component of $(\alpha\cup\gamma)\cap \mathcal{R}^d$,  $k_1$ and $k_2$ must both
intersect a component of $(\alpha\cup\gamma)\cap\mathcal{R}^d$. This implies that the two intersection points of $k_1\cup\rho\cup k_2$ with 
a component of $(\alpha\cup\gamma)\cap\mathcal{R}^d$ have opposite signs. This contradicts part (2) of 
Proposition~\ref{pro:Paths}.
\end{proof}

Choose product  structures $S^1\times I$ on each of the annuli $\mathcal{A}_l$ and $\mathcal{A}_r$  and 
$I\times I$ on each of the rectangles $\mathcal{R}^u$ and $\mathcal{R}^d$, so that  each component of 
$(\alpha\cup\gamma)\cap\widehat{\Sigma}$ is of the form $\{x\}\times I$ in $\mathcal{A}_l$, $\mathcal{A}_r$, 
$\mathcal{R}^u$ or $\mathcal{R}^d$.  We call an arc of the form $\{x\}\times I$ a \emph{vertical} arc in the 
annuli or rectangles.  After performing an isotopy on $\delta$, we may assume each arc in $\delta\cap \mathcal{A}_l$ and $\delta\cap \mathcal{A}_r$ is transverse to all the vertical arcs of $\mathcal{A}_l$, $\mathcal{A}_r$.

\begin{figure}[!ht]
	\begin{overpic}[width=7.5cm]{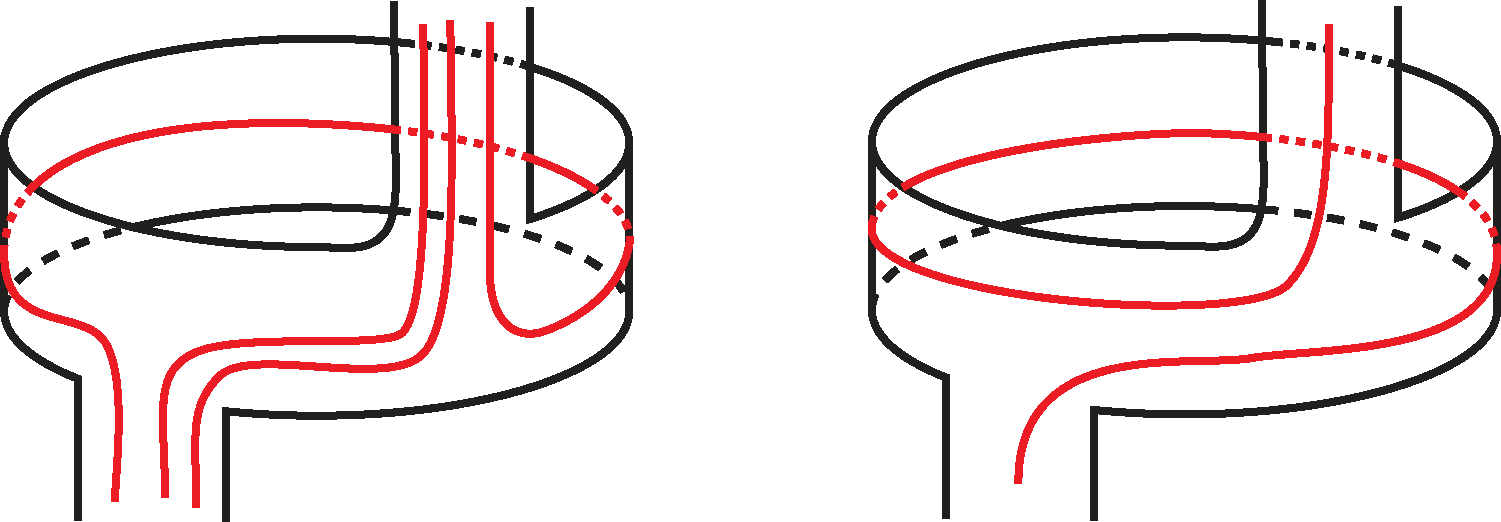}
		\put(29,34){$\delta$}
		\put(88,34){$\delta$}
		\put(-3,-6){$\delta$ takes 2 short paths}
		\put(59,-6){$\delta$ takes a long path}
	\end{overpic}
	\vspace{10pt}
	\caption{$\delta$-paths}
	\label{fig:path}
\end{figure}

\begin{definition}\label{def:Path} {\bf paths, long and short}:
Let $l$ be a subarc of $\delta$ that is properly embedded in $\mathcal{A}_r$ or $\mathcal{A}_l$.
By Lemma~\ref{lem:NoCircular}, $l$ is an arc connecting $\mathcal{R}^d$ to $\mathcal{R}^u$.  By the assumption before Definition~\ref{def:Path}, $l$ is transverse to all the vertical arcs. We call $l$ a {\it regular arc} if one of the following two conditions is satisfied 
\begin{enumerate}
	\item either $l$ is intersects a vertical arc of $\mathcal{A}_r$ or $\mathcal{A}_l$ more than once (in other words, the curve $l$ wraps around the annulus more than once), or
	\item $l$ shares endpoints with two components of $\delta\cap(\mathcal{R}^d\cup\mathcal{R}^u)$ that are core arcs of $\mathcal{R}^d$ and $\mathcal{R}^u$ respectively. 
\end{enumerate} 
The purpose of Condition (2) is to rule out a special case that $l$ is attached to the component of $\delta\cap \mathcal{R}^u$ or $\delta\cap \mathcal{R}^d$ which contains $X^u$ or $X^u$ (e.g.~the arc $k^d$ described before Lemma~\ref{lem:NoCircular}).  
A {\it path} is an isotopy class of regular arcs which setwise preserves 
$\mathcal{R}^d$ and $\mathcal{R}^u$. Two regular arcs of $\delta$ are said to \emph{belong to the same path} if the two arcs are isotopic via an isotopy that setwise preserves $\mathcal{R}^d$ and $\mathcal{R}^u$.  We say that 
$\delta$  {\it takes a path} if it contains a subarc that belongs to that path.  A path is 
called a {\it long path} in $\mathcal{A}_r$ or  $\mathcal{A}_l$ if a regular arc in that path intersects a vertical arc of $\mathcal{A}_r$ or $\mathcal{A}_l$ more than once (i.e.~it satisfies Condition (1) above).  Otherwise it will be called a {\it short} path.  Note that $\delta$ 
can take at most two different paths in one annulus and that there may be multiple arcs in the same 
path, see Figure~\ref{fig:path} for a picture.  We would like to emphasize that, in a short path, both ends of a regular arc must attach to core arcs of $\mathcal{R}^d$ and $\mathcal{R}^u$ in $\delta$. 
\end{definition}

The decomposition of $\widehat{\Sigma}$ divides $\delta$ into a collection of regular arcs and core arcs of $\mathcal{R}^d$ and $\mathcal{R}^u$, as well as two special arcs connecting $X^d$ and $X^u$ to $\partial \mathcal{R}^d$ and $\partial \mathcal{R}^u$.  The two special arcs join at $X$ to form a subarc of $\delta$ containing $X$, which we denote by $\rho_x$. 
By Definition~\ref{def:Path}, each end of $\rho_x$ is attached to either a regular arc in a long path or a core arc in $\mathcal{R}^d$ and $\mathcal{R}^u$ in $\delta$.

\begin{figure}[ht]
	\begin{overpic}[width=6cm]{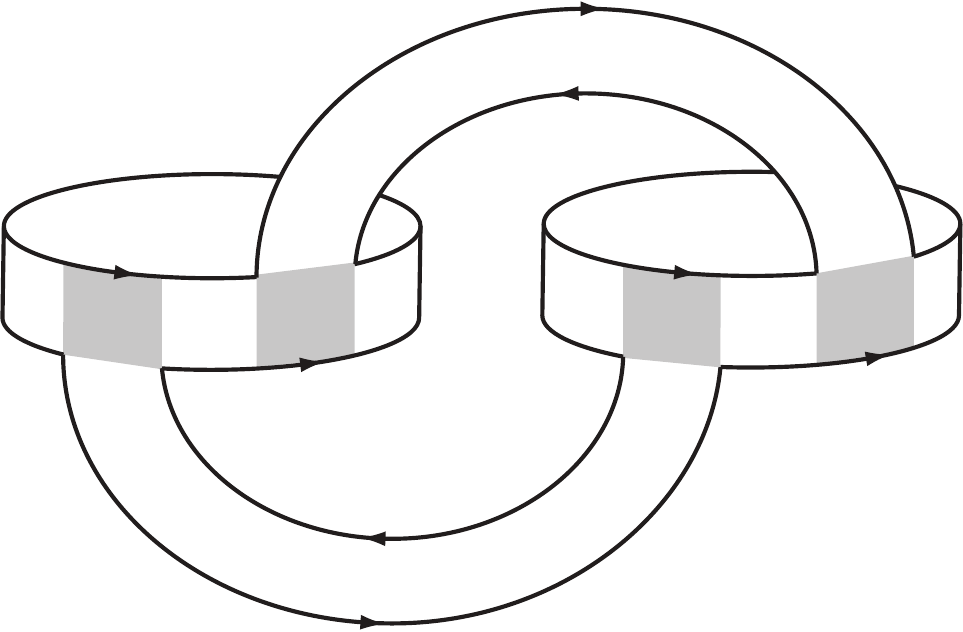}
		\put(7.5,30){{\footnotesize $\mathcal{J}_l^d$}}
		\put(28,31.5){{\footnotesize $\mathcal{J}_l^u$}}
		\put(65.5,30){{\footnotesize $\mathcal{J}_r^d$}}
		\put(86.5,31.5){{\footnotesize $\mathcal{J}_r^u$}}
		\put(58,57.5){{\footnotesize $\mathcal{R}^u$}}
		\put(35,2.5){{\footnotesize $\mathcal{R}^d$}}
		\put(-9,34){$\mathcal{A}_l$}
		\put(101,34){$\mathcal{A}_r$}
	\end{overpic}
	\caption{The junctions in $\widehat{\Sigma}$}\label{fig:junctions}
\end{figure}

\begin{definition}\label{def:junctions} {\bf Junctions}:
	The intersections $\partial\mathcal{R}^u \cap \partial \mathcal{A}_l$  and  
	$\partial\mathcal{R}^d \cap \partial \mathcal{A}_l$ are subarcs of  different components 
	of $ \partial \mathcal{A}_l$,  see  Figure~\ref{circular}. Choose the product structure on 
	$\mathcal{A}_l$ so that no vertical arc  has one endpoint in  $\partial\mathcal{R}^u$ and 
	the other end point in  $\partial\mathcal{R}^d$.
	
	As illustrated in Figure~\ref{fig:junctions}, let $\mathcal{J}_l^u$ and $\mathcal{J}_l^d$ be the unions 
	of all vertical arcs of $\mathcal{A}_l$ with  an endpoint in $\partial\mathcal{R}^u$ and $\partial\mathcal{R}^d$ 
	respectively. Similarly let  $\mathcal{J}_r^u$ and $\mathcal{J}_r^d$ be the unions of vertical arcs of 
	$\mathcal{A}_r$ with an endpoint in $\partial\mathcal{R}^u$ and $\partial\mathcal{R}^d$ 
	respectively. With this choice $\mathcal{J}_l^u$,  $\mathcal{J}_l^d$, $\mathcal{J}_r^u$ and $\mathcal{J}_r^d$ 
	are disjoint squares in $\mathcal{A}_l$ and $\mathcal{A}_r$, see Figure~\ref{fig:junctions}. Call these  four squares 
	{\it junctions}.
\end{definition}

The purpose of introducing the structure of junctions in Definition~\ref{def:junctions} is to construct train tracks which are our main tools in analyzing curves.

\begin{definition}\label{def:traintracks} {\bf Train tracks}:  The four junctions give $\widehat{\Sigma}$ a 
decomposition into ten rectangles: $\mathcal{R}^u$, $\mathcal{R}^d$, $\mathcal{J}_l^u$,  $\mathcal{J}_l^d$, 
$\mathcal{J}_r^u$, $\mathcal{J}_r^d$, $\mathcal{A}_l\ssm (\mathcal{J}_l^u\cup\mathcal{J}_l^d)$ and 
 $\mathcal{A}_r\ssm (\mathcal{J}_r^u\cup\mathcal{J}_r^d)$.  When  all the arcs of $\delta\setminus\rho_x$ in each of these 
 rectangles with endpoints in the same pair of boundary edges are collapsed into a single arc, we obtain a train track.
Now, we consider how the special arc $\rho_x$ is affected by the train-track construction. 
If $\delta$ does not take a long path, then the train-track construction pinches the two endpoints of $\rho_x$ onto the train track, and the construction does not affect $\Int(\rho_x)$, see Figure~\ref{fig:pinch-to-traintrack} for a local picture. 
If $\delta$ takes a long path, then two subarcs of $\rho_x$ at its ends may be affected by the train-track construction. 
For example, in the left picture of Figure~\ref{fig:longisotopy}, $\rho_x$ contains arcs in the middle of the shaded region where parallel arcs on both sides of $\rho_x$ are collapsed into a single arc in the train track. 
If this happens, we collapse these two subarcs of $\rho_x$ (which are at the two ends of $\rho_x$) onto the train track (i.e.~collapse all the arcs in shaded region of Figure~\ref{fig:longisotopy} into a single segment in the train track). 
In either case, this train track, together with $\rho_x$, forms a slightly larger train track, which we denote by  $\tau_D$, see Figure~\ref{fig:pinch-to-traintrack} and the left picture of Figure~\ref{fig:longisotopy}. 
The special arc $\rho_x$ becomes a special segment of $\tau_D$, which we still denote by $\rho_x$.

\begin{figure}[!ht]
	\begin{overpic}[width=8cm]{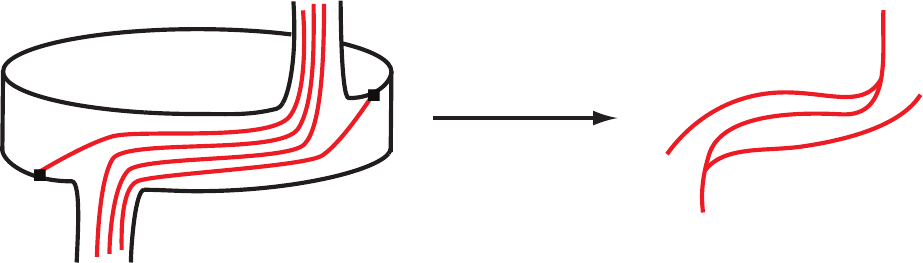}
		\put(82,8){$\tau_D$}
		\put(99.5,20){$\rho_x$}
		\put(69,9){$\rho_x$}
	\end{overpic}
	\caption{Pinch $\delta$ into a traintrack $\tau_D$}
	\label{fig:pinch-to-traintrack}
\end{figure}
 
 The switches (or cusp points) of the train track divide the train track into a collection of arcs which are 
 called the {\it segments} of the train track.  The {\it weight} of $\delta$ at a segment of $\tau_D$ is the 
 number of $\delta$-arcs which collapse onto this segment. (In other  words, the number of times that 
 $\delta$ passes this segment.) 
 It follows from our construction that the arc $\rho_x$ is a segment of $\tau_D$ that contains the 
 point $X$, see Figure~\ref{fig:pinch-to-traintrack}. Clearly the weight of $\delta$ at the segment $\rho_x$ 
 is one by the construction.   Moreover, by part (2) of Proposition~\ref{pro:Paths}, each segment of 
 $\tau_D$ has an induced orientation from $\delta$ and the orientations of the segments are compatible 
 at each cusp.
 \end{definition}

\begin{figure}[!ht]
\begin{overpic}[width=9cm]{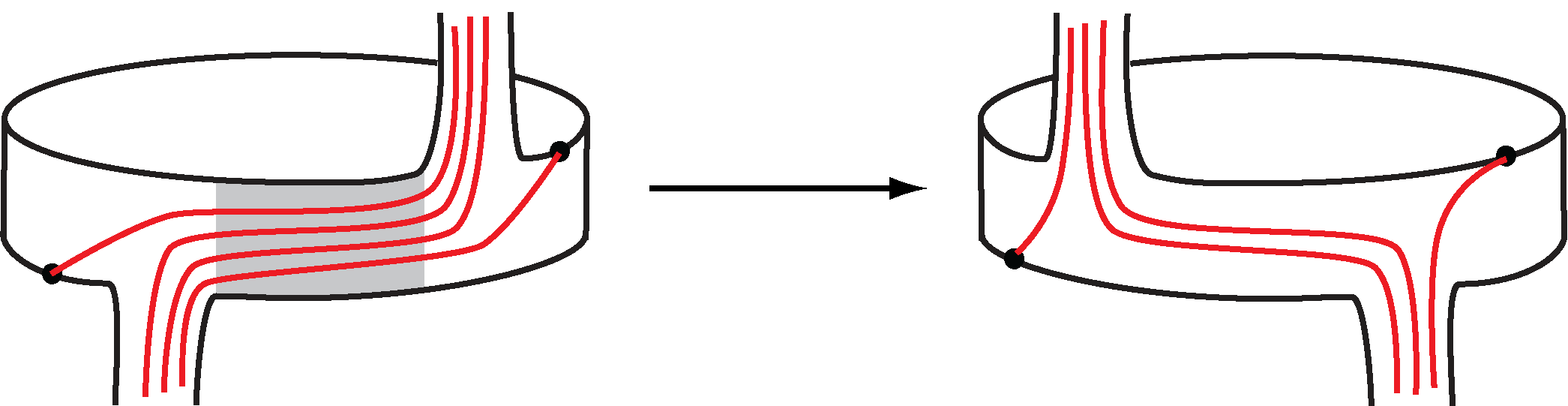}
\put(43,10){isotopy}
\end{overpic}
\caption{Isotopy on an annulus}
\label{fig:isotopy}
\end{figure}

\begin{remark}
In the decomposition of $\widehat{\Sigma} = \Sigma \ssm \NN(\partial_+P)$, the product structures on the 
two rectangles $\mathcal{R}^u$ and $\mathcal{R}^d$ depend on their intersections with $\gamma\cup\alpha$, 
see the $\gamma$ arcs in Figure~\ref{circular}. So the product structures of $\mathcal{R}^u$ and $\mathcal{R}^d$ 
are ``basically'' fixed with respect to $\delta$, and there is no ambiguity in the traintrack construction.  
	
The product structures on $\mathcal{A}_l$ and $\mathcal{A}_r$ are somewhat flexible on the 
regions which do not intersect $\gamma\cup\alpha$. For example, if $\gamma\cup\alpha$ does 
not intersect the shaded region in Figure~\ref{fig:isotopy}, then we can perform an isotopy 
on the annulus as in in Figure~\ref{fig:isotopy} that changes the product structure.  The effect 
of this isotopy on the traintrack $\tau_D$ is performing a splitting move on the train track. 
\end{remark}

\begin{definition}\label{def:splitting-arc}
Consider $\tau_D$ as a traintrack in the Heegaard surface $\Sigma$ and  a small neighborhood of 
$\tau_D$ in $\Sigma$, denoted by $\mathcal{N}(\tau_D)$. We can view $\mathcal{N}(\tau_D)$ as an $I$-bundle 
over the traintrack $\tau_D$. The curve $\delta \subset\mathcal{N}(\tau_D)$  is transverse to the $I$-fibers of 
$\mathcal{N}(\tau_D)$. Similar to a traintrack, $\mathcal{N}(\tau_D)$ has a collection of cusps at its boundary. 
There is a collection of arcs $s_1,\dots,s_k$ in $\mathcal{N}(\tau_D)$  which are disjoint from $\delta$ and are 
transverse to the $I$-fibers, and which connect the cusps of $\mathcal{N}(\tau_D)$ in pairs. If one cuts open 
$\mathcal{N}(\tau_D)$ along $s_1,\dots,s_k$  the resulting space is an $I$-bundle over $\delta$. We call each 
$s_i$ a {\it splitting arc}.
\end{definition}



\chapter{Proof  of the main theorem}\label{cpt:ObtainingTheContradiction}

In this chapter we give the proof of Theorem~\ref{thm:MainTheorem} using the terms defined and facts proved
in the previous chapters. 

We begin with a description of the arc $\delta \ssm \partial_+P \subset \widehat{\Sigma}$ connecting $X^u$ to 
$X^d$. The possible locations for $X^u$ and $X^d$ in $\partial_+P^u$ and $\partial_+P^d$, respectively, 
are not determined by the construction as the precise gluing map between $\partial_+P^u$  and $\partial_+P^d$ 
is unknown. Nonetheless, since $X$ is assumed to be a fixed point of the hyperelliptic involution $\pi$, the points 
$X^u$ and $X^d$ are symmetric under $\pi$.   

\vskip5pt
It follows from Lemma~\ref{lem:TwoDeltaArcs} that $\delta$ must take at least one path in each annulus $\mathcal{A}_l$, $\mathcal{A}_r$. 
 The following  are the four possible  configurations in $\widehat{\Sigma}$  for this arc:

\vskip5pt
\noindent \underline{Configuration 1}: $\delta$ takes one short path in $\mathcal{A}_l$ and one short path in $\mathcal{A}_r$.
\vskip5pt
\noindent \underline{Configuration 2}: $\delta$ takes one short path in $\mathcal{A}_l$ or $\mathcal{A}_r$ and two short paths 
in the other annulus.
\vskip5pt
\noindent \underline{Configuration 3}: $\delta$ takes two short paths in $\mathcal{A}_l$ and two short paths in $\mathcal{A}_r$.
\vskip5pt
\noindent  \underline{Configuration 4}: $\delta$ takes a long path in one (or both) of the annuli.
\vskip7pt
A key difference between these four configurations can be seen in the train track $\tau_D$ (constructed in 
Definition \ref{def:traintracks}).   

\vskip5pt

\begin{figure}[!ht]
\begin{overpic}[width=3.5cm]{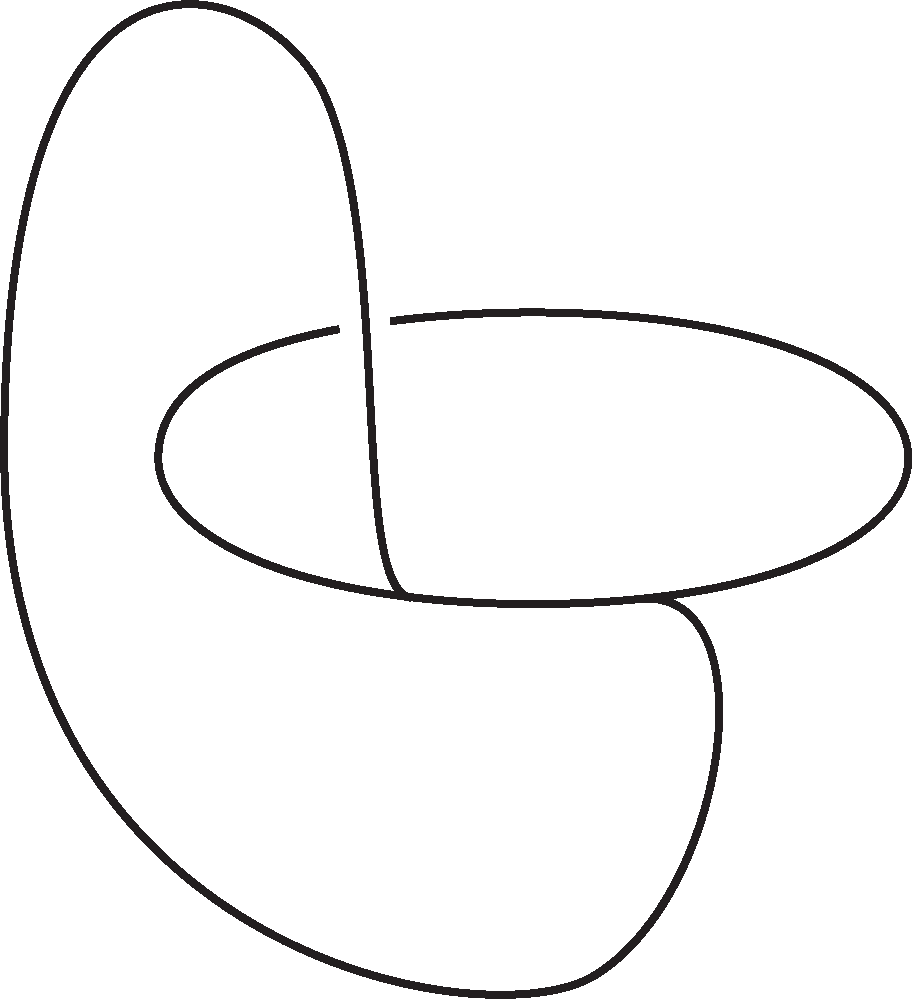}
\put(41,33){$x$-arc}
\put(-25,50){$y$-arc}
\put(47,71){$y$-arc}
\end{overpic}
\caption{A train track whose complement is a once-punctured torus}
\label{fig:traintrack}
\end{figure}

In  Configuration 1 the two short paths and a core arc from each rectangle  $\mathcal{R}^u$ and $\mathcal{R}^d$ 
are pinched into a circle, and the train track $\tau_D$ is the union of this circle and the special segment $\rho_x$ which 
contains the point $X = \delta \cap  \partial_+P$, as illustrated in the construction of the train track $\tau_D$ in 
Figure~\ref{fig:pinch-to-traintrack}. Thus $\tau_D$ in the first configuration must be as shown in  Figure~\ref{fig:traintrack}. 
The $x$- and $y$-arcs marked in Figure~\ref{fig:traintrack} are the three segments of $\tau_D$, where the $x$-arc is the 
segment with cusp directions at its endpoints pointing into this arc. The special segment $\rho_x$ is one of the $y$-arcs in
Figure~\ref{fig:traintrack}.  Note that in this configuration, $\Sigma\ssm\tau_D$ is a once-punctured torus with two cusps at 
its boundary.  

In Configuration 2, see Figure~\ref{fig:OneTwo1}, 
the train track $\tau_D$ can be obtained by adding an additional segment to the train track in Figure~\ref{fig:traintrack}.  
Hence $\tau_D$ in this configuration has four cusps. Since $\partial_+P \cap \tau_D$ is just a single 
point $X$, the train track $\tau_D$ is non-separating in $\Sigma$. The complement $\Sigma\ssm\tau_D$ in 
Configuration 2 can be obtained by cutting the once-punctured torus in Configuration 1 along an essential 
arc, which means that $\Sigma\ssm\tau_D$ in Configuration 2 is a topological annulus.  Moreover, since 
every essential curve in a genus-2 surface is invariant by the hyperelliptic involution, after isotopy, we 
may assume $\tau_D$ is invariant by the involution and this implies that  the annulus $\Sigma\ssm\tau_D$ 
has two cusps at each boundary curve.  

Similarly, in Configuration 3, since $\tau_D\cap\partial_+P$ is a single point,   
$\Sigma\ssm\tau_D$ is connected. As $\tau_D$ has six cusps in this configuration, 
for Euler characteristic reason, $\Sigma\ssm\tau_D$ is a disk with six cusps at its boundary. 

In Configuration 4, $\tau_D$ takes a long path in one annulus, say $\mathcal{A}_r$. Then
$\tau_D$ has two possible configurations depending on whether $\delta$ takes  one short path, 
two short paths, or another long path in $\mathcal{A}_l$. If $\tau_D$ takes one short path in $\mathcal{A}_l$, 
then $\tau_D$ has four cusps, and similar to Configuration 2, the surface $\Sigma\ssm\NN(\tau_D)$ is an 
annulus with two cusps on each boundary curve. If $\tau_D$ takes two short paths or another long path in 
$\mathcal{A}_l$, then $\tau_D$ has six cusps, and similar to Configuration 3, $\Sigma\ssm\tau_D$
is a disk with six cusps at its boundary. 

We discuss  each of these configurations in separate sections. The  sections below deal with the 
configurations above in order  2, 3, 4 and 1 for convenience of the proof.


\section{The curve $\delta$ takes three short paths}\label{sec:OneArTwoAl}

\begin{figure}[ht]
\begin{overpic}[width=7cm]{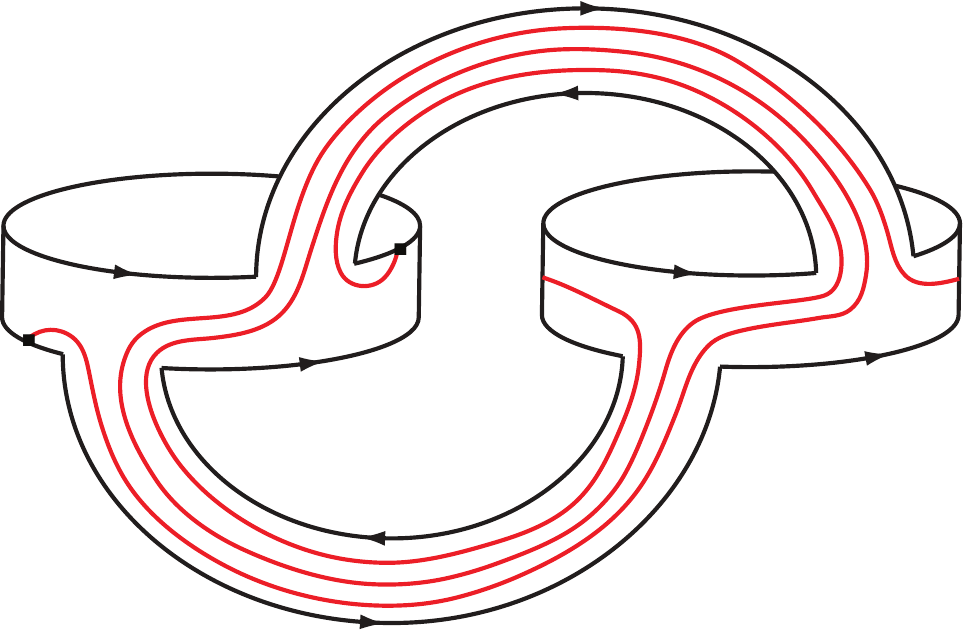}
\put(32,33){$\delta$}
\put(-3,24.5){$X^d$}
\end{overpic}
\caption{A configuration for $\delta$ that takes one short path in $\mathcal{A}_l$ and two short paths 
in $\mathcal{A}_r$}
\label{fig:OneTwo1}
\end{figure}

In this section, Configuration 2 is considered, where $\delta$ takes two short paths in one of the annuli,
 and a single short path in the other annulus. An example of this situation is given in Figure~\ref{fig:OneTwo1}. 

\begin{proposition}\label{pro:OneRTwoL}
If  $\delta$ takes two short paths in one annulus and a single short path in another annulus, then
 $K$  is doubly primitive.
\end{proposition}

\begin{proof}
After applying Symmetry (2) in Remark~\ref{rem:Symmetries}
if necessary,  assume without loss of generality that $\delta$ takes  two short paths in $\mathcal{A}_r$ and 
one short path in $\mathcal{A}_l$.   Since Symmetry (3) in Remark~\ref{rem:Symmetries} interchanges the two 
short paths in $\mathcal{A}_l$,  assume  further that $\delta$ takes a fixed short path in $\mathcal{A}_l$, 
see Figure~\ref{fig:OneTwo1}.

As explained at the beginning of the chapter, $\Sigma\ssm\NN(\tau_D)$ is an annulus  denoted by $A_D$.  
The annulus $A_D$ has two cusps in each boundary component  of $A_D$ and the cusps correspond to the 
cusps of the train track $\tau_D$.

By  construction, the train track $\tau_D$ has a special segment $\rho_x$ which contains the intersection 
point $X$ and the weight of $\delta$ at $\rho_x$ is one. The two ends of $\rho_x$ are two cusps of $\tau_D$, 
which are called the $\rho$-cusps.  The other two cusps of $\tau_D$ are at the two junctions of $\mathcal{A}_r$ 
where the two short paths meet, see Figure~\ref{fig:OneTwo1}. We call these two cusps the $j$-cusps.  As 
$\delta$ is invariant under  $\pi$ and $X$ is a fixed point of $\pi$, we may assume  $\rho_x$ is invariant under 
$\pi$ and $\pi$ interchanges the two $\rho$-cusps and the two $j$-cusps.

A neighborhood $\NN(\tau_D)$ of $\tau_D$ contains two splitting arcs as in Definition~\ref{def:splitting-arc}.
Denoted them by $s_\rho$ and $s_j$. These splitting arcs connect the cusps of $\NN(\tau_D)$ in pairs such that 
$s_\rho$ and $s_j$ are disjoint from $\delta$, and  if one splits $\NN(\tau_D)$ along $s_\rho$ and $s_j$, the 
resulting surface $\NN(\tau_D)\ssm (s_\rho\cup s_j)$ is a product neighborhood of $\delta$.  

We claim that one splitting arc, say $s_\rho$,  connects the two $\rho$-cusps and the other splitting arc $s_j$ 
connects the two $j$-cusps. To see this, consider the modified train track $\tau_D^-$ obtained by removing 
$\rho_x$ from $\tau_D$. As in the discussion at the beginning of the chapter, the train track $\tau_D^-$ is 
as in Figure~\ref{fig:traintrack} and $\NN(\tau_D^-)$ is a once-punctured torus.  We view $\NN(\tau_D^-)$ as 
an $I$-bundle over  $\tau_D^-$ and view $s_\rho$ and $s_j$ as arcs in $\NN(\tau_D^-)$ transverse to the 
$I$-fibers. 

Assume, in contradiction, $s_\rho$ connects a $\rho$-cusp to a $j$-cusp. Then $s_\rho$ connects a cusp 
of $\NN(\tau_D^-)$ to a smooth boundary point of $\NN(\tau_D^-)$. Since $s_\rho$ is transverse to the $I$-fibers, 
this implies that $s_\rho$ is nontrivial in the once-punctured torus $\NN(\tau_D^-)$. Similarly $s_j$ is also nontrivial 
$\NN(\tau_D^-)$. 

The symmetry induced by $\pi$ implies that the two endpoints of $\rho_x$ (i.e., the two 
$\rho$-cusps) and the two $j$-cusps alternate along the boundary curve of $\NN(\tau_D^-)$, see 
Figure~\ref{fig:OneTwo1}. 

If $s_\rho$ connects a $\rho$-cusp to a $j$-cusp, then the endpoints $\partial s_\rho$ and $\partial s_j$ 
do not alternate along $\partial\NN(\tau_D^-)$. However, the endpoints of any pair of disjoint non-parallel 
essential arcs in a once-punctured torus do alternate along the boundary curve. Thus $s_\rho$ and $s_j$ 
must be parallel in the once-puncture torus $\NN(\tau_D^-)$. 
Hence, $\NN(\tau_D^-)\ssm (s_\rho\cup s_j)$ has two components, one of which is an annulus. 
Since $s_\rho$ and $s_j$ are transverse to the 
$I$-fibers of $\NN(\tau_D^-)$, this annular component of $\NN(\tau_D^-)\ssm (s_\rho\cup s_j)$ must 
be an $I$-bundle over a simple closed curve with each $I$-fiber a subarc of an $I$-fiber of $\NN(\tau_D^-)$. 
This is  impossible because $s_\rho$ and $s_j$ are the splitting arcs of $\NN(\tau_D)$ and 
$\NN(\tau_D)\ssm (s_\rho\cup s_j)$ is an $I$-bundle over $\delta$, i.e., we have a contradiction. 
Thus, one splitting arc, say,  $s_\rho$ connects the two $\rho$-cusps and the other $s_j$ connects 
the two $j$-cusps.

\begin{claim}\label{cl:splitting-arc}
Each of the splitting arcs $s_\rho$ and $s_j$ passes through both $\mathcal{R}^u$ and 
$\mathcal{R}^d$. In other words, each contains subarcs that are core arcs of $\mathcal{R}^u$ 
and $\mathcal{R}^d$. In particular, each splitting arc intersects $\gamma$ at least twice.
\end{claim}

\begin{proof}[Proof of Claim~\ref{cl:splitting-arc}]
Since the splitting arc $s_j$ connects the two $j$-cusps, as illustrated in Figure~\ref{fig:OneTwo1}, $s_j$ 
must pass through both rectangles $\mathcal{R}^d$ and $\mathcal{R}^u$. 
Now we consider $s_\rho$.   

Recall the construction of the special arc $\rho_x$. 
As discussed before Definition~\ref{def:junctions}, in general, each end of $\rho_x$ is attached either to a long path or to a component of $\delta\cap(\mathcal{R}^d\cup\mathcal{R}^u)$. 
Since there is no long path in the configuration of Proposition~\ref{pro:OneRTwoL}, the two ends of $\rho_x$ are attached to two components of $\delta\cap(\mathcal{R}^d\cup\mathcal{R}^u)$ that are core arcs of $\mathcal{R}^d$ and $\mathcal{R}^u$ respectively.

In the construction of the train track $\tau_D$, all the core arcs of $\delta\cap\mathcal{R}^d$ (and $\delta\cap\mathcal{R}^u$) are pinched onto the same segment of $\tau_D$. As illustrated in Figure~\ref{fig:pinch-to-traintrack}, this means that, similar to the two $j$-cusps, the cusp directions at the two ends of $\rho_x$ must point into $\mathcal{R}^d$ and $\mathcal{R}^u$ respectively. Thus, $s_\rho$ passes through both $\mathcal{R}^d$ and $\mathcal{R}^u$.
\end{proof}

 As $\varepsilon\cap\delta=\emptyset$, the intersection $\varepsilon\cap\NN(\tau_D)$ is a collection of 
 arcs  through the cusps of $\NN(\tau_D)$. By part (2) of Proposition~\ref{pro:Paths}, we can choose an 
 orientation on $\varepsilon$ so that the induced orientation on all the arcs in  $\varepsilon\cap\NN(\tau_D)$  
 is compatible with the orientation of $\delta$ near the cusps of $\NN(\tau_D)$, see the arrows in 
 Figure~\ref{fig:Annulus1}.  

\begin{figure}[ht]
\begin{overpic}[width=5.5cm]{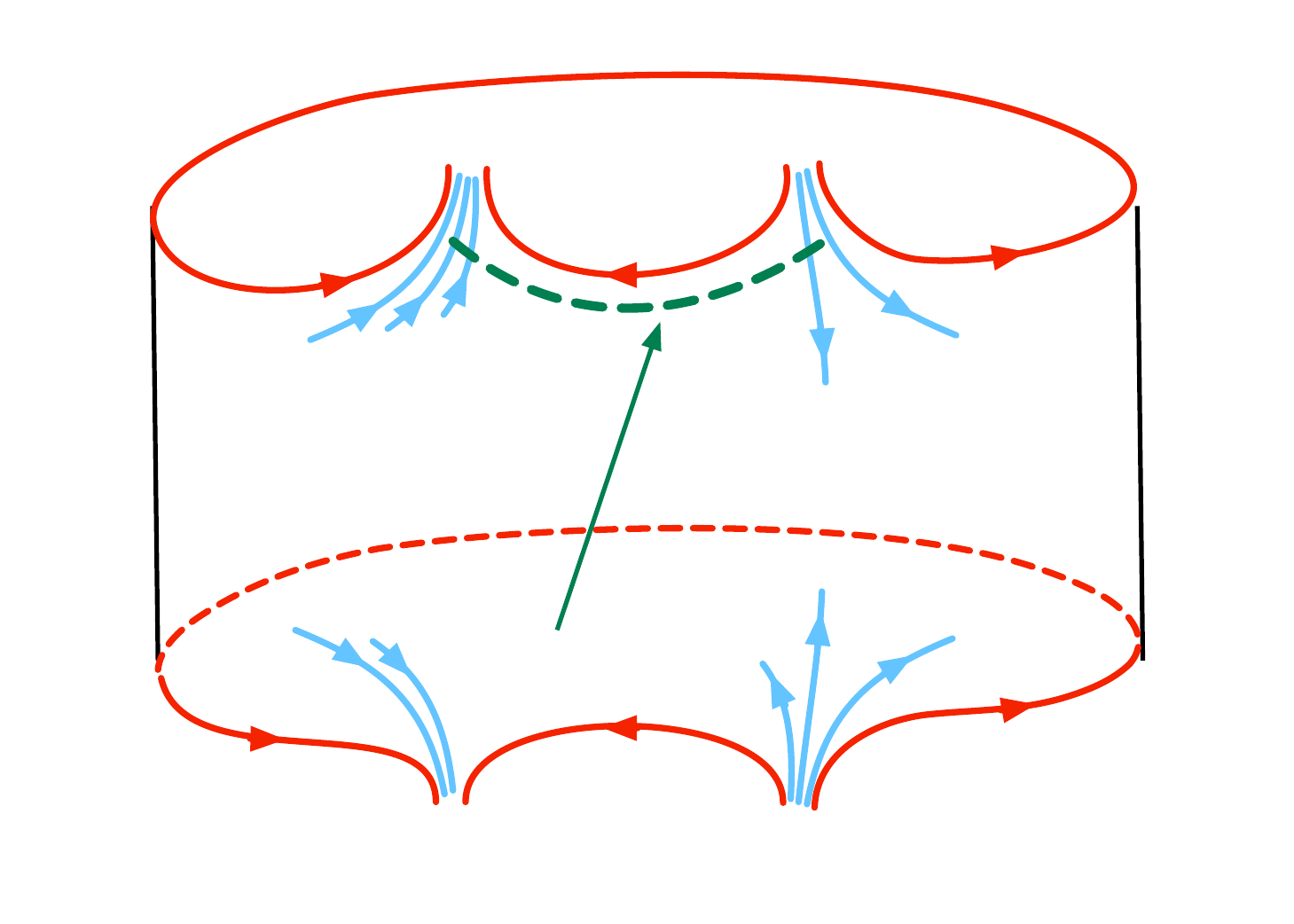}
\put(40,18.5){$\rho$}
\put(29,41){$\varepsilon$}
\put(68,42){$\varepsilon$}
\put(47,59){$\delta$}
\end{overpic}
\caption{The annulus $A_D=\Sigma\ssm\tau_D$ and the $s$-wave $\rho$ with respect to $\varepsilon$.}
\label{fig:Annulus1}
\end{figure}

If $\varepsilon$ goes through all the cusps of $\NN(\tau_D)$, then as shown in Figure~\ref{fig:Annulus1},
 there must be an $s$-wave $\rho$ with respect to $\varepsilon$ which is parallel to a $\delta$-arc. Perform 
 a sequence of wave moves  on $\varepsilon$ along $s$-waves similar to $\rho$ (as in Figure~\ref{fig:Annulus1}) 
 so that the resulting curve $\varepsilon'$ goes into at most one cusp in each boundary component of 
 $A_D=\Sigma\ssm\NN(\tau_D)$.  This implies that the arcs in  $\varepsilon'\cap\NN(\tau_D)$, if there 
 are any, are all parallel to either $s_\rho$ or $s_j$.

\begin{claim}\label{cl:MoreThanOneArc}
The intersection $\varepsilon'\cap\NN(\tau_D)$ contains at least two arcs.
\end{claim}

\begin{proof}[Proof of  Claim~\ref{cl:MoreThanOneArc}]

 The proof is by contradiction and is divided into two cases: 
 \begin{enumerate}[(i)]
 \item  $|\,\varepsilon'\cap\NN(\tau_D)\,| = 0$, i.e.~$\varepsilon'\cap\NN(\tau_D)=\emptyset$,  or
 \vskip4pt
  \item  $|\,\varepsilon'\cap\NN(\tau_D)| = 1$. 
\end{enumerate}
 
\vskip5pt
\noindent{\it Case (i)}:  $\varepsilon'\cap\NN(\tau_D)=\emptyset$, i.e., $\varepsilon'$ is isotopic to 
the core curve of the annulus $A_D$.  
\vskip5pt

Since  $\partial_+P\cap\rho_x = \partial_+P\cap\delta$ is the single point $X$, the intersection
 $\partial_+P\cap A_D$ is an essential arc in $A_D$.  As $\varepsilon'$ is isotopic to the 
core curve of $A_D$, the intersection $\varepsilon'\cap\partial_+P$ is a single point.   Let $E'$ 
denote the disk in the handlebody $W$ which is bounded by $\varepsilon'$. Hence $(P, E')$ is 
a $(\mathcal{P},\mathcal{D})$-pair.  
  
\vskip5pt

Showing that  the $(\mathcal{P},\mathcal{D})$-pair $(P, E')$ has a smaller complexity than $(P, D)$ will  rule 
out Case (i). To compute the complexity, consider the arcs in $(\alpha\cup\gamma)\cap\widehat{\Sigma}$, 
or equivalently arcs in $(\alpha\cup\gamma)\ssm\partial_+P$.  There are two possible sub-cases: 

The first subcase  is that when a component of $(\alpha\cup\gamma)\ssm\partial_+P$ intersects $\varepsilon'$ 
then it also intersects $\delta$.  In this subcase,  $c_0(P,E',\alpha,\gamma)\le c_0(P,D,\alpha,\gamma)$ by 
the definition of the  complexity.

The second subcase is that there is a component of $(\alpha\cup\gamma)\ssm\partial_+P$ 
which intersects $\varepsilon'$ but does not intersect $\delta$. 
Let $\mathcal{C}$ denote the union of the components of 
$(\alpha\cup\gamma)\cap A_D$ which are not $\partial$-parallel in  $A_D$. So, in Case (i), we have 
$|(\alpha\cup\gamma)\cap\varepsilon'|=|\mathcal{C}|$.  

The arcs in $\mathcal{C}$ can be viewed as vertical arcs of the annulus $A_D$ and  $\partial_+P\cap A_D$ 
can be  viewed as a spiral in $A_D$. In the second subcase, since $A_D=\Sigma\ssm\tau_D$, a component of 
$(\alpha\cup\gamma)\ssm\partial_+P$ which intersects $\varepsilon'$ but does not intersect $\delta$ must be 
contained in the interior of $A_D$. This implies that the arc in $\mathcal{C}$ which contains such a component of 
$(\alpha\cup\gamma)\ssm\partial_+P$ must meet $\partial_+P$ more than once. 
Hence the curve $\partial_+P\cap A_D$ wraps around $A_D$ more than once, which means 
that each arc in $\mathcal{C}$ must intersect $\partial_+P$ at least once.  Thus, for each 
component $\sigma$ of $\mathcal{C}$, we can pick one intersection point of $\sigma\cap \partial_+P$, and so we have a total of $|\mathcal{C}|$ such intersection points. These $|\mathcal{C}|$ intersection 
points divide $\alpha\cup\gamma$ into $|\mathcal{C}|$ arcs, each of which meets $\delta$.  Since 
these $|\mathcal{C}|$ points are in $\partial_+P$, each of the these $|\mathcal{C}|$ arcs contains 
at least one component of  $(\alpha\cup\gamma)\ssm\partial_+P$ which meets $\delta$. By our 
definition of complexity, this implies that $c_0(P,D,\alpha,\gamma)\ge |\mathcal{C}|$. 
Moreover, by the definition of the complexity,   
$$c_0(P,E',\alpha,\gamma)\le |(\alpha\cup\gamma)\cap\varepsilon'|.$$  
As $|(\alpha\cup\gamma)\cap\varepsilon'|=|\mathcal{C}|$, we have 
$c_0(P,E',\alpha,\gamma)\le c_0(P,D,\alpha,\gamma)$ in both subcases. 

\vskip5pt

Since each arc in $\mathcal{C}$ intersects $\varepsilon'$ in one point,  Lemma~\ref{lem:TwoDeltaArcs} 
implies that  in fact 
$$|(\alpha\cup\gamma)\cap\varepsilon'| < |(\alpha\cup\gamma)\cap\delta|$$ 
and thus  
$$c(P,E',\alpha,\gamma)< c(P,D,\alpha,\gamma).$$

\vskip5pt
 
This contradicts  Assumption~\ref{ass:minimality}  on the choice of 
the pair $(P, D)$ and therefore  $\varepsilon'\cap\NN(\tau_D)\ne\emptyset$.

\vskip8pt
\noindent{\it Case (ii)}:    $\varepsilon'\cap\NN(\tau_D)$ is a single arc.
\vskip5pt

Consider the train track $\tau_D^-$ obtained by removing the segment $\rho_x$ from $\tau_D$.  Note that 
$\NN(\tau_D^-)$ is a once-punctured torus and $\delta\cap\NN(\tau_D^-)$ is a single arc in $\NN(\tau_D^-)$. 
Since $\varepsilon' \cap \rho_x = \emptyset$ we have $\varepsilon'\cap\NN(\tau_D)=\varepsilon'\cap\NN(\tau_D^-)$. 
Hence  $\varepsilon'\cap\NN(\tau_D^-)$  is also a single arc in $\NN(\tau_D^-)$.  Further, note that the arc 
$\varepsilon'\cap\NN(\tau_D)$ is parallel to one of the two splitting arcs $s_\rho$ or $s_j$.  Similar to the 
argument before Claim~\ref{cl:splitting-arc},  the endpoints $\partial s_\rho$ (or $\partial s_j$) and the 
endpoints $\partial\rho_x$ alternate along the boundary of $\NN(\tau_D^-)$.  Hence the endpoints of the arc 
$\varepsilon'\cap\NN(\tau_D^-)$ and the endpoints of the arc $\delta\cap\NN(\tau_D^-)$ alternate along  the 
boundary of $\NN(\tau_D^-)$. As  $\NN(\tau_D^-)$ is a once-punctured torus, this implies that $\delta\cap\NN(\tau_D^-)$ 
and $\varepsilon'\cap\NN(\tau_D^-)$ are non-parallel essential arcs in $\NN(\tau_D^-)$. Thus we may assume that 
$\delta\cap\NN(\tau_D^-)$ is of slope $1/0$ and  $\varepsilon'\cap\NN(\tau_D^-)$ is of slope $0/1$ with respect to a 
basis of the relative first  homology of the  once-punctured torus $\NN(\tau_D^-)$.

Since $\delta$ takes two short paths in $\mathcal{A}_r$, we can isotope the core curve $\mathfrak{a}_r$ 
of $\mathcal{A}_r$ into $\NN(\tau_D^-)$.  So there is an arc of slope $p/q$ in $\NN(\tau_D^-)$ interesting 
$\mathfrak{a}_r$ in a single point.  Let $E'$ and $D$ be the disks bounded by $\varepsilon'$ and $\delta$ in 
$W$ respectively. We can take $p$ parallel copies of $D$ and $q$ parallel copies of $E'$, and perform a 
sequence of band sums of these disks along the boundary of $\NN(\tau_D^-)$ to obtain a disk $D'$ in $W$ 
so that $\partial D'\cap \NN(\tau_D^-)$ is an arc of slope $p/q$. Hence $\partial D'$  intersects $\mathfrak{a}_r$ 
in a single point.  Recall that $\mathfrak{a}_r$ is a boundary curve of a planar surface $P_r$, 
i.e.,~$\mathfrak{a}_r=\partial_+P_r$, where $P_r$ is obtained by a single $\partial$-compression on $P$, see 
Remark~\ref{rem:core}. 
Thus the pair $(P_r, D')$ is a new $(\mathcal{P},\mathcal{D})$-pair.

Now compute the complexity of this new  $(P_r, D')$ pair. The curve $\mathfrak{a}_r$ can be viewed as the 
union of two segments of $\tau_D$.   Since $\delta$ intersects every $\alpha$- and $\gamma$-arc in 
$\mathcal{A}_r$ and intersects at least two other $\gamma$-arcs in the rectangles $\mathcal{R}^d$ 
and $\mathcal{R}^u$, we have $$|(\alpha\cup\gamma)\cap\mathfrak{a}_r| < c_0(P,D,\alpha,\gamma).$$ 
By the definition of the complexity, 
$$c_0(P_r,D',\alpha,\gamma)\le |(\alpha\cup\gamma)\cap\partial_+P_r|=|(\alpha\cup\gamma)\cap\mathfrak{a}_r|.$$ 
Thus, when putting the two inequalities together we have
 $c_0(P_r,D',\alpha,\gamma)<c_0(P,D,\alpha,\gamma)$ 
again, a contradiction to Assumption~\ref{ass:minimality}. This finishes the proof of the claim.
\end{proof}

Next, consider the Heegaard diagram formed by the $\{\delta,\varepsilon'\}$ and $\{\alpha,\gamma\}$ sets 
of curves.

\begin{claim}\label{cl:Vwave}
The curve $\gamma$ has two distinct subarcs that are  
$[\varepsilon'^-,\varepsilon'^+]$-edges and also has two distinct subarcs that are 
$[\delta^-,\delta^+]$-edges. 
\end{claim}

\begin{proof}[Proof of Claim~\ref{cl:Vwave}] 
As illustrated by the arrows in Figure~\ref{fig:Annulus1}, before the wave moves on 
$\varepsilon$, the orientations of $\delta$ and $\varepsilon$ are compatible within 
$\NN(\tau_D)$ (See part (2) of Proposition~\ref{pro:Paths}). The wave moves on 
$\varepsilon$, that were performed above, are disjoint from $\NN(\tau_D)$. Hence  there is an induced  
orientation on $\varepsilon'$ and the induced orientations of the arcs $\varepsilon'\cap\NN(\tau_D)$ are 
compatible with the orientation of $\delta$.
 
In the paragraph before Claim~\ref{cl:MoreThanOneArc},  we concluded that the arcs in $\varepsilon'\cap\NN(\tau_D)$ 
are all parallel to either $s_\rho$ or $s_j$ in $\NN(\tau_D)$.  By Claim~\ref{cl:MoreThanOneArc}, $\varepsilon'\cap\NN(\tau_D)$ 
contains at least two such parallel arcs. As the orientations of the arcs $\varepsilon'\cap\NN(\tau_D)$ are compatible, 
the intersection points of each arc of $\gamma\cap\NN(\tau_D)$ with $\varepsilon'$ all have the same sign. 
Thus each subarc of $\gamma\cap\NN(\tau_D)$ between two curves of $\varepsilon'\cap\NN(\tau_D)$ is an 
$[\varepsilon'^-,\varepsilon'^+]$-edge. By Claim~\ref{cl:splitting-arc}, $\gamma$ intersects each arc in 
$\varepsilon'\cap\NN(\tau_D)$ at least twice, so $\gamma$ has at least two distinct subarcs that are 
$[\varepsilon'^-,\varepsilon'^+]$-edges. Similarly, since $\varepsilon'$ does not go through all cusps, it follows 
from Claim~\ref{cl:splitting-arc} that $\gamma$ has two distinct subarcs that are   $[\delta^-,\delta^+]$-edges. 
\end{proof}

In the Heegaard diagram formed by $\{\delta,\varepsilon'\}$ and $\widehat{V}=\{\alpha,\gamma\}$, the 
existence of both $[\delta^-,\delta^+]$ and $[\varepsilon'^-,\varepsilon'^+]$-edges implies that the Whitehead
Graph $\Gamma(\delta,\varepsilon')$ must be of type (i) in Figure~\ref{graphs} with $c\ne 0$ and $d\ne 0$. 
In fact $c\ge 2$ and $d\ge 2$ by Claim~\ref{cl:Vwave}. By part (1) of Lemma~\ref{lem:+ and -}, the Heegaard 
diagram contains no wave with respect to $\{\delta,\varepsilon'\}$. By Theorems \ref{thm:HOT} and \ref{thm:NeOk}, 
it must contain a wave with respect to $\{\alpha,\gamma\}$. 

Since every essential simple closed curve is invariant by the involution $\pi$ after perhaps some isotopy, 
one may assume that  $\varepsilon'$, as well as $\delta$, $\partial_+P$, $\alpha$ and $\gamma$, 
are invariant under $\pi$. It also follows from part (1) of Proposition~\ref{pro:Paths}  
that we can fix an orientation for $\alpha$ and $\gamma$  so that the intersection points of $\alpha$ 
and $ \gamma$ with $\delta$ all have the same sign.

Now, assume that $\eta$ is a  wave with respect to $\{\gamma,\alpha\}$ and assume further that  
$|\eta\cap\partial_+P|$  is minimal among all such waves. It follows from Lemma~\ref{lem:NoAWaves} 
that $\eta\cap\partial_+P\ne\emptyset$. The curve $\partial_+P$ divides $\eta$ into a sequence of 
subarcs $\kappa_1,\dots,\kappa_s$  so that the first and last subarcs  $\kappa_1$ and $\kappa_s$ 
contain the endpoints of $\eta$.  Next, consider the  $\kappa_i, i = 1,\dots, s$, as arcs in $\widehat{\Sigma}$.  
The two arcs $\kappa_1$ and $\kappa_s$ at the ends are arcs connecting $\gamma\cup\alpha$ to 
$\partial\widehat{\Sigma}$.  All other arcs $\kappa_i$, $i=2,\dots,s-1$, are properly embedded in 
$\widehat{\Sigma}$.

For each junction $\mathcal{J}$  (see Definition~\ref{def:junctions}) in the decomposition of $\widehat{\Sigma}$, 
let $\widehat{\mathcal{J}}$ be a hexagonal neighborhood of $\mathcal{J}$ in $\widehat{\Sigma}$, see 
Figure~\ref{fig:junction-arc}(a). Any arc $\sigma$ in $\widehat{\Sigma}$ connecting a component of 
$(\gamma\cup\alpha)\cap\widehat{\Sigma}$ to $\partial\widehat{\Sigma}$, is called a 
{\it junction arc} if after isotopy on $\widehat{\Sigma}$, the arc $\sigma$ lies in a hexagonal neighborhood 
$\widehat{\mathcal{J}}$ of a junction $\mathcal{J}$ connecting a pair of opposite edges of 
$\widehat{\mathcal{J}}$, see the dashed arc in Figure~\ref{fig:junction-arc}(a).

\begin{figure}[ht]
	\begin{overpic}[width=11cm]{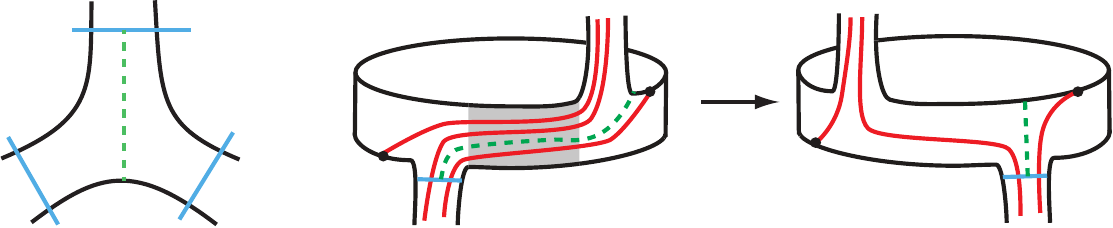}
		\put(9,-5){(a)}
		\put(45,-5){(b)}
		\put(84,-5){(c)}
		\put(15.5,11){$\partial_+P$}
		\put(62,13){{\tiny isotopy}}
		\put(89.5,8.5){{\tiny $\kappa_1$}}
	\end{overpic}
	\vspace{10pt}
	\caption{Junction arcs}
	\label{fig:junction-arc}
\end{figure}

\begin{claim}\label{cl:junction-arc}
Let  $\eta$ and $\kappa_1,\dots,\kappa_s$ be as above. Then each $\kappa_i$, $i=2,\dots, s-1$ 
is isotopic to a vertical arc of $\mathcal{A}_l$ and both $\kappa_1$ and $\kappa_s$ are junction arcs.
\end{claim} 

\begin{proof}[Proof of Claim~\ref{cl:junction-arc}]
First, consider the arcs $\kappa_i$, $i=2,\dots, s-1$, which are properly embedded in $\widehat{\Sigma}$. 
By Lemma~\ref{lem:GammaRectangles}, $\gamma$ intersects both $\mathcal{R}^u$ and
$\mathcal{R}^d$, see Figure~\ref{circular}.  As $\eta$ is a wave, $\eta\cap\delta=\emptyset$ and $\kappa_i\cap\gamma=\emptyset$ for 
$i=2,\dots, s-1$. Since 
 $\delta$ passes through $\mathcal{R}^u$ and $\mathcal{R}^d$ and since $\kappa_i\cap\gamma=\emptyset$ ($i=2,\dots, s-1$), 
$\kappa_i$  is either $\partial$-parallel in $\widehat{\Sigma}$ or an essential arc in $\mathcal{A}_l$ or 
$\mathcal{A}_r$. If $\kappa_i$ is boundary parallel, then the bigon disk between $\kappa_i$ and 
$\partial_+P$ can be eliminated by a simple  isotopy. This isotopy yields a new wave with fewer 
intersection points with $\partial_+P$, contradicting our hypothesis on $\eta$.  Thus $\kappa_i$, 
$i=2,\dots, s-1$ is isotopic to a vertical arc of $\mathcal{A}_l$ or $\mathcal{A}_r$, and since it does 
not intersect $\delta$ it must be in $\mathcal{A}_l$.

Next consider $\kappa_1$ (the proof for $\kappa_s$ is identical).  Let $\sigma$ be the component of 
$(\gamma\cup\alpha)\cap\widehat{\Sigma}$ that contains an endpoint $Z$ of $\kappa_1$  (which is also 
an endpoint of the wave $\eta$).  As illustrated in Figure~\ref{circular}, two $\gamma$-arcs in 
$\mathcal{R}^u$ and $\mathcal{R}^d$ divide $\widehat{\Sigma}$ into two topological annuli. 
Moreover, since $\Int(\kappa_1) \cap \gamma = \emptyset$,  both $\sigma$ and $\kappa_i$ lie in one 
of the two annuli. Let $\sigma'$ and $\sigma''$ be the two components of $\sigma\ssm Z$. If both 
$\kappa_1\cup\sigma'$ and $\kappa_1\cup\sigma''$ are nontrivial arcs in these two annuli,  as illustrated 
in Figure~\ref{fig:junction-arc}(b, c), perform an isotopy so that $\kappa_1$ connects a pair of opposite 
edges of a hexagonal neighborhood of a junction and hence $\kappa_1$ is a junction arc. Note that 
since $\eta$ is a wave, if $\kappa_1$ is as in Figure~\ref{fig:junction-arc}(b), then the shaded region 
in Figure~\ref{fig:junction-arc}(b) cannot contain any $\alpha$- or $\gamma$-arc. Hence the 
isotopy from Figure~\ref{fig:junction-arc}(b) to (c) is in fact the isotopy in Figure~\ref{fig:isotopy}.

\begin{figure}[ht]
	\begin{overpic}[width=8cm]{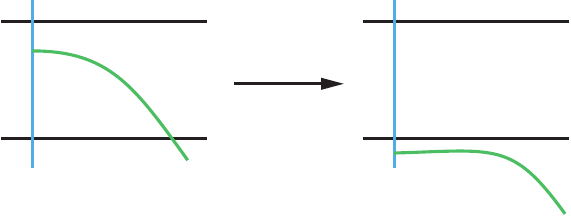}
		\put(18,-2){(a)}
		\put(75,-2){(b)}
		\put(-10,33){$\partial_+P$}
		\put(-10,12){$\partial_+P$}
		\put(12,17){$\Delta$}
		\put(1,26.5){$Z$}
		\put(18,26.5){$\kappa_1$}
		\put(32,6){$\eta$}
		\put(41,25){isotopy}
	\end{overpic}
	\vspace{8pt}
	\caption{An isotopy on $\eta$}
	\label{fig:wave-arc1}
\end{figure}

It remains to consider the case that either $\kappa_1\cup\sigma'$ or $\kappa_1\cup\sigma''$ is a 
$\partial$-parallel arc in $\widehat{\Sigma}$. Suppose $\kappa_1\cup\sigma'$ is $\partial$-parallel, 
then $\kappa_1\cup\sigma'$ and a subarc of $\partial\widehat{\Sigma}$ bound a triangle 
$\Delta\subset\widehat{\Sigma}$, see Figure~\ref{fig:wave-arc1}(a). 

If $\Delta\cap(\delta\cup\varepsilon')=\emptyset$, one can perform an isotopy pushing $\kappa_1$ 
across  $\Delta$,  as can be seen in  Figure~\ref{fig:wave-arc1}.  This yields a new wave with fewer 
intersection points with $\partial_+P$, contradicting our hypothesis on $\eta$. 

\begin{figure}[ht]
	\begin{overpic}[width=9.5cm]{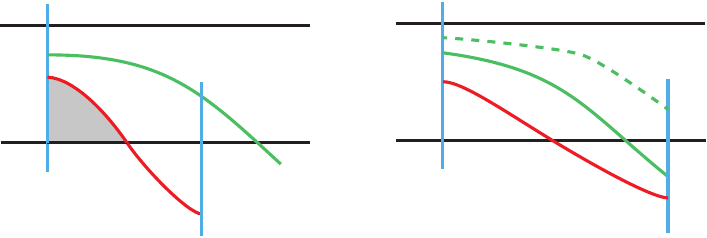}
		\put(20,-4){(a)}
		\put(78,-4){(b)}
		\put(-8,29){$\partial_+P$}
		\put(-8,12){$\partial_+P$}
		\put(8,15){{\small $\Delta'$}}
		\put(20.5,10.5){{\tiny $\pi(\Delta')$}}
		\put(29,6){${\small\pi(\sigma)}$}
		\put(19,17){$\Delta$}
		\put(65,15){$\Delta'$}
		\put(76,17){$\Delta$}
		\put(22,24){$\kappa_1$}
		\put(40,8){$\eta$}
		\put(3.5,17){$\sigma$}
		\put(60,16){$\sigma$}
		\put(95.5,8){$\pi(\sigma)$}
	\end{overpic}
\vspace{8pt}
	\caption{Possible configurations of $\kappa_1$}
	\label{fig:wave-arc2}
\end{figure}

If $\Delta\cap(\delta\cup\varepsilon')\ne\emptyset$, then $\delta$ and $\varepsilon'$ cut off a sub-triangle 
$\Delta'$ in $\Delta$, see the shaded triangle in Figure~\ref{fig:wave-arc2}(a). The triangle $\Delta'$  is rotated 
by the involution $\pi$ into a triangle  $\pi(\Delta')$ on the other side of $\partial_+P$. Since $\delta$, $\varepsilon'$  
and $\partial_+P$ are invariant under $\pi$, there are two possible configurations, corresponding to whether the vertical 
arc on the boundary of $\pi(\Delta')$ meets $\partial_+P$ in $\Delta$ or outside it. The two configurations are illustrated in 
Figures~\ref{fig:wave-arc2}(a) and (b): 

\vskip5pt

\noindent (a) Figure~\ref{fig:wave-arc2}(a) is the case where an endpoint of $\pi(\sigma)$ lies on the boundary 
of $\Delta$.  In this case the curve $\alpha$ or $\gamma$ that contains $\pi(\sigma)$ must cut into $\Delta$ 
and meet $\eta$. This contradicts the fact that $\eta$ is a wave and $\eta\cap (\alpha\cup\gamma)=\partial\eta$. 
So the situation in Figure~\ref{fig:wave-arc2}(a) cannot occur.

\vskip5pt

\noindent (b) Figure~\ref{fig:wave-arc2}(b) is the configuration where a vertex of $\Delta$ lies 
in the boundary of $\pi(\Delta')$. In this case,  the wave $\eta$ must meet $\pi(\sigma)$.  Since 
$\eta\cap(\alpha\cup\gamma)=\partial\eta$,  this means that the union of $\kappa_1$ and the 
extension of $\kappa_1$ in $\pi(\Delta')$, (see the green arc in Figure~\ref{fig:wave-arc2}(b)), 
must be the whole of $\eta$.  In this case the wave $\eta$ can be isotoped to a wave that is 
disjoint from $\partial_+P$, see the dashed arc in Figure~\ref{fig:wave-arc2}(b), contradicting 
Lemma~\ref{lem:NoAWaves}. 
\end{proof}

It follows from Claim~\ref{cl:junction-arc} that, after isotopy, any wave $\eta$ contains two junction arcs at 
its two ends. Note that in a hexagonal neighborhood of a junction, two junction arcs connecting different pairs of 
opposite edges must intersect. Hence the two junction arcs in $\eta$ belong to two different junctions.  
For each wave $\eta$, the  involution $\pi$ sends $\eta$ to a disjoint dual wave $\pi(\eta)$ (see Remark~\ref{rem:duralwave}). Hence the four junction arcs of $\eta$ and $\pi(\eta)$ must be
at different junctions.  Since   $\pi$ leaves $\mathcal{A}_l$ and $\mathcal{A}_r$ invariant, $\pi$ interchanges 
the two junctions at $\mathcal{A}_l$ and interchanges the two junctions at $\mathcal{A}_r$.  As $\eta$ and 
$\pi(\eta)$ are distinct waves, this implies that one junction arc of $\eta$ must be at a junction of 
$\mathcal{A}_r$ and the other junction arc is at a junction of $\mathcal{A}_l$.

Without loss of generality, suppose $\kappa_1$ is a junction arc in $\mathcal{A}_r$ and 
$\kappa_s$ is a junction arc in $\mathcal{A}_l$. By assumption, $\delta$ takes two short paths in 
$\mathcal{A}_r$.  Therefore, as shown in Figure~\ref{fig:wave-move}(a), since the junction arc 
$\kappa_1$ is disjoint from $\delta$, up to interchanging $\eta$ and $\pi(\eta)$, the arc $\kappa_1$ 
must connect an arc in $(\gamma\cup\alpha)\cap \mathcal{R}^u$ to $\partial\mathcal{A}_r$.

\begin{figure}[ht]
	\vspace{5pt}
	\begin{overpic}[width=11cm]{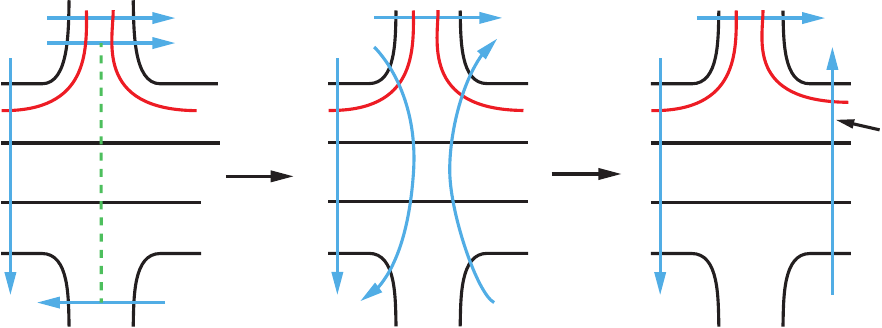}
		\put(10,36.5){$\mathcal{R}^u$}
		\put(9.5,-4){(a)}
		\put(46,-4){(b)}
		\put(83,-4){(c)}
		\put(-5,23){$\mathcal{A}_r$}
		\put(23,23){$\mathcal{A}_r$}
		\put(7.5,23.5){$\kappa_1$}
		\put(7.5,8){$\kappa_s$}
		\put(12.5,17){$\eta$}
		\put(2,16.5){$\alpha$}
		\put(38.5,16){$\alpha$}
		\put(43,18){$\alpha_0$}
		\put(52,17){$\gamma'$}
		\put(75.5,17){$\alpha$}
		\put(91,16){$\gamma'$}
		\put(100,20){\parbox{2cm}{arc next to the cusp}}
		\end{overpic}
	\vspace{8pt}
	\caption{Wave move along $\eta$ in two steps}
	\label{fig:wave-move}
\end{figure}

Next,  perform a wave move along $\eta$ and  study how the new curve, obtained by the wave move, 
intersects $\widehat{\Sigma}$, focusing on how the  curve is changed near the arc $\kappa_1$.
By Lemma~\ref{lem:FirstWave}, the wave $\eta$ must be a wave with respect to $\gamma$. 
Perform the wave move along $\eta$ in two steps as in Definition~\ref{def:wavemove}. 
The surgery-step changes $\gamma$ into two curves $\gamma'$ and $\alpha_0$ where 
$\alpha_0$ is parallel to $\alpha$, and the second step is to remove $\alpha_0$ and obtain 
a new set of meridians $\{\alpha, \gamma'\}$.

Since $\kappa_1$ is a junction arc in $\mathcal{A}_r$, as shown in Figure~\ref{fig:wave-move}(a, b), 
the surgery-step of the wave move changes the arc in $\gamma\cap\mathcal{R}^u$ that contains an 
endpoint of $\eta$ into two vertical arcs of $\mathcal{A}_r$ next to the junction, one in $\gamma'$ and 
the other in $\alpha_0$. 

By Claim~\ref{cl:junction-arc}, the intersection of both curves $\gamma'$ and $\alpha_0$ with $\widehat{\Sigma}$ 
is a collection of nontrivial arcs with respect to the decomposition of $\widehat{\Sigma}$. Furthermore, 
after slightly modifying the product structure of $\mathcal{A}_l$ and $\mathcal{A}_r$ if necessary, 
we may assume both $\gamma'$ and $\alpha_0$ meet $\widehat{\Sigma}$ in a collection of vertical 
arcs in $\mathcal{A}_l$, $\mathcal{A}_r$, $\mathcal{R}^u$, and $\mathcal{R}^d$. 

Since $\alpha_0$ is parallel to $\alpha$ and as shown in Figure~\ref{fig:wave-move}(b, c), 
after deleting $\alpha_0$ in the second step of the wave move, one of two arcs of 
$(\alpha\cup\gamma')\cap\mathcal{A}_r$ next to a junction belongs to $\gamma'$ and the 
other is in $\alpha$.  These  two arcs are depicted in Figure~\ref{fig:wave-move}(c) 
and are called  \emph{arcs in $\mathcal{A}_r$ next to the cusp of $\delta$}. 
Moreover, since the intersection points  of $\delta$ with $\alpha$ and $\gamma$ all 
have the same sign, the direction of the two arcs next to the cusp must be opposite in 
$\mathcal{A}_r$, see Figure~\ref{fig:wave-move}(c).

Let  $\mathcal{J}_r^u$ and $\mathcal{J}_r^d$ be the two junctions in $\mathcal{A}_r$
and denote the two components of $\mathcal{A}_r\ssm(\mathcal{J}_r^u\cup\mathcal{J}_r^d)$ 
by  $R_\alpha$ and $R_\gamma$. Since the involution $\pi$ interchanges the  junctions 
$\mathcal{J}_r^u$ and $\mathcal{J}_r^d$, the induced action of $\pi$  on each of 
$R_\alpha$ and $R_\gamma$ is a $180^\circ$-rotation. As a consequence of this symmetry, 
the two outermost arcs of $(\alpha\cup\gamma')\cap R_\alpha$ must belong the the same 
curve $\alpha$ or $\gamma'$.  Since the two arcs in $\mathcal{A}_r$ next to the cusp of 
$\delta$ are in different curves, see Figure~\ref{fig:wave-move}(c),  the two outermost 
arcs of $(\alpha\cup\gamma')\cap R_\alpha$ must be in the same curve, say $\alpha$, and the 
two outermost  arcs of $(\alpha\cup\gamma')\cap R_\gamma$ are in $\gamma'$. Note that it 
is possible that $(\alpha\cup\gamma')\cap R_\alpha$ or $(\alpha\cup\gamma')\cap R_\gamma$ is 
a single arc. Thus, after the first wave move along $\eta$, the curves $\gamma'$ and 
$\alpha$ have the following \emph{two properties}:

\begin{enumerate}
	\item $(\alpha\cup\gamma')\cap R_\gamma\ne\emptyset$ and each outermost arc of 
	$(\alpha\cup\gamma')\cap R_\gamma$ is in $\gamma'$.

	\item $(\alpha\cup\gamma')\cap R_\alpha\ne\emptyset$ and each outermost arc of 
	$(\alpha\cup\gamma')\cap R_\alpha$ is in $\alpha$.
\end{enumerate}

If $\gamma'\cup\alpha$ still intersects  $\mathcal{R}^u$ and $\mathcal{R}^d$, then similar to 
Claim~\ref{cl:Vwave}, two subarcs of $\gamma'\cup\alpha$ are $[\delta^+,\delta^-]$-edges 
and two subarcs of $\gamma'\cup\alpha$ are  $[\varepsilon'^+,\varepsilon'^-]$-edges,
ruling out a wave with respect to $\{\delta, \varepsilon'\}$. The same argument, as above, implies 
that there must be a wave $\eta'$ with respect to $\{\alpha,\gamma'\}$.  

 The next step is to perform a wave move along $\eta'$. Before proceeding, notice that 
$\eta'$ has similar properties as $\eta$.  A property of $\eta$ that was used in the proof of  
Claim~\ref{cl:junction-arc} is that $\eta\cap\partial_+P\ne\emptyset$, which follows from 
Lemma~\ref{lem:NoAWaves}. The next claim says that the same is true for the following wave
$\eta'$.

\begin{claim}\label{cl:NoAWaves}
Any wave with respect to $\{\alpha,\gamma'\}$ must intersect $\partial_+P$.
\end{claim}

\begin{proof}[Proof of Claim~\ref{cl:NoAWaves}]
We prove this claim using the property that, before the wave move along $\eta$, every wave with 
respect to $\{\alpha,\gamma\}$ must intersect $\partial_+P$.  

Assume in contradiction that  there is a wave $\omega$ with respect to $\{\alpha,\gamma'\}$ which 
is disjoint from $\partial_+P$. We first show that $\omega$ must be a wave with respect to $\gamma'$.  
To see this, consider a pair of pants $Q$ which is a neighborhood of $\gamma\cup\eta$ before the wave 
move along $\eta$. So we may view the three boundary curves of $Q$ as $\gamma$, $\alpha$, and $\gamma'$. 
If $\omega$ is a wave with respect to $\alpha$ in the Heegaard diagram formed by $\{\alpha,\gamma'\}$ and 
$\{\delta,\varepsilon'\}$, then $\omega\cap\gamma'=\emptyset$. If $\omega$ is also disjoint from $\gamma$, 
then $\omega$ is a wave with respect to $\{\alpha,\gamma\}$ and disjoint from $\partial_+P$, a contradiction 
to the property that there is no such wave before the wave move along $\eta$. Thus, $\omega$ must intersect 
$\gamma$ and hence $\omega$ has subarcs in $Q$ as well as outside $Q$. As $\omega\cap\gamma'=\emptyset$, this 
means that a subarc of $\omega$ must be a wave with respect to $\gamma$ in the Heegaard diagram formed by 
$\{\alpha,\gamma\}$ and $\{\delta,\varepsilon'\}$. This again contradicts the property above since 
$\omega\cap\partial_+P=\emptyset$. Therefore,  $\omega$ must be a wave with respect to $\gamma'$.

As described by the surgery-step of the wave move, we may view $\gamma'$ as the union of $\eta$ and 
a subarc $\gamma_1$ of $\gamma$, i.e.~$\gamma'=\eta\cup\gamma_1$ and $\gamma_1\subset\gamma$.
Since the original $\gamma$ has no wave disjoint from $\partial_+P$, there are two possibilities: 
\begin{enumerate}
\item $\partial\omega\subset\eta$.

\vskip5pt

\item $\omega$ has one endpoint in $\eta$ and the other endpoint in $\gamma_1$,
\end{enumerate}

\vskip5pt

\noindent (1) If $\partial\omega\subset\eta$, then the union of $\omega$ and two subarcs of $\eta$ is a wave 
with respect to $\gamma$, see the dashed arc in Figure~\ref{fig:omega}(a).  Since $\omega\cap\partial_+P=\emptyset$, 
we may view $\omega$ as an arc in $\widehat{\Sigma}$. Note that the endpoints of $\omega$ cannot be on the same 
vertical segment of $\eta\ssm\partial_+P$ as in this case $\omega$ must be trivial. Thus, as illustrated in 
Figure~\ref{fig:omega}(a), this new wave has fewer intersection points with $\partial_+P$, contradicting 
our assumption on $\eta$.

\vskip5pt

\noindent (2) If $\omega$ has one endpoint in $\eta$ and the other endpoint in $\gamma_1$, then 
the endpoint of $\omega$ in $\eta$ divides $\eta$ into two subarcs $\eta_1$ and $\eta_2$. 
Both $\omega\cup\eta_1$ and $\omega\cup\eta_2$ are waves with respect to $\gamma$, see the  
two dashed arcs in Figure~\ref{fig:omega}(b). Moreover, since $\omega\cap\partial_+P=\emptyset$, 
at least one of $\omega\cup\eta_1$ and $\omega\cup\eta_2$ has fewer intersection points with 
$\partial_+P$, contradicting our assumption on $\eta$. 

The two contradictions prove this claim.
\end{proof}

\begin{figure}[ht]
	\begin{overpic}[width=6cm]{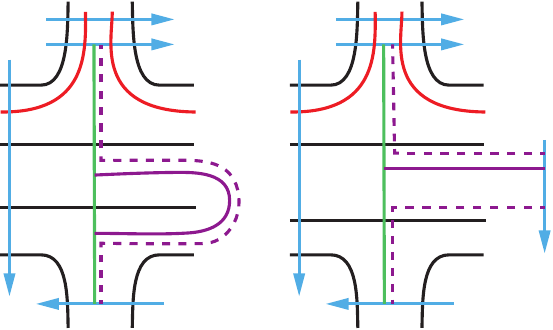}
		\put(15,-6){(a)}
		\put(67,-6){(b)}
		\put(13,24){$\eta$}
		\put(30,24){$\omega$}
		\put(64,36.5){{\footnotesize $\eta_1$}}
		\put(64,14){{\footnotesize $\eta_2$}}
		\put(85,25){$\omega$}
		\put(100,25){$\gamma_1$}
	\end{overpic}
	\vspace{8pt}
	\caption{Construct new waves using $\omega$}
	\label{fig:omega}
\end{figure}

Claim~\ref{cl:NoAWaves} implies that Claim~\ref{cl:junction-arc} also holds for $\eta'$. 
In particular, $\eta'$ is divided by $\partial_+P$ into a collection of vertical arcs in $\mathcal{A}_l$ 
and a pair of junction arcs. Moreover, one of the junction arcs is in $\mathcal{A}_r$ and the 
other is in $\mathcal{A}_l$, see the dashed arc $\eta$ in Figure~\ref{fig:wave-move}(a) for 
a picture of $\eta'$. 

Now, perform a wave move along $\eta'$. If $\partial\eta'\subset\alpha$,  call the curve 
resulting from the wave move an \emph{$\alpha$-curve}, and if $\partial\eta'\subset\gamma'$, 
 call the resulting curve a \emph{$\gamma$-curve}.  To simplify notation, we also call $\alpha$ 
and $\gamma'$ an $\alpha$- and a $\gamma$-curve respectively.

Next, apply wave moves repeatedly enabled by Claims~\ref{cl:NoAWaves} and \ref{cl:junction-arc}. 
To simplify  notation, if a wave move is along a wave with respect to an $\alpha$- or a $\gamma$-curve, 
we always name the resulting curve $\alpha'$ or $\gamma'$ respectively. 

If $(\alpha'\cup\gamma')\cap\mathcal{R}^u\ne\emptyset$ (and  $(\alpha'\cup\gamma')\cap\mathcal{R}^d\ne\emptyset$ 
by symmetry), then similar to Claim~\ref{cl:Vwave}, two subarcs of $\alpha'\cup\gamma'$ are $[\delta^+,\delta^-]$-edges 
and two subarcs are $[\varepsilon'^+,\varepsilon'^-]$-edges.  Moreover, the two properties after the first wave move that were 
emphasized before Claim~\ref{cl:NoAWaves} always hold for $\{\alpha',\gamma'\}$.  This implies that the Heegaard diagram formed by 
$\{\delta,\varepsilon'\}$ and  $\{\alpha',\gamma'\}$ is not a standard Heegaard diagram of $M$.  
Because of the $[\delta^+,\delta^-]$ and  $[\varepsilon'^+,\varepsilon'^-]$ blocking-edges, 
 there is no wave with respect to $\{\delta, \varepsilon'\}$. So, Theorems~\ref{thm:HOT} and \ref{thm:NeOk} 
 imply that there  must be a wave with respect to $\{\alpha',\gamma'\}$.  Furthermore,  Claims~\ref{cl:NoAWaves} 
 and \ref{cl:junction-arc}  also hold for the new wave with respect to $\{\alpha',\gamma'\}$ and we can continue 
 with the wave moves described above.

Since each wave move reduces  the number of intersection points, this process will eventually end. Stop 
performing the wave moves whenever $(\alpha'\cup\gamma')\cap(\mathcal{R}^d\cup\mathcal{R}^u)=\emptyset$. 
 This means that, at this stage, every arc in $(\alpha'\cup\gamma')\cap\widehat{\Sigma}$ is a vertical arc in either 
 $\mathcal{A}_l$ or $\mathcal{A}_r$.

As above,  call each component of $\alpha'\cap \widehat{\Sigma}$ an $\alpha'$-arc and each component of  
$\gamma'\cap \widehat{\Sigma}$ a $\gamma'$-arc. It follows from our construction that the two properties 
before Claim~\ref{cl:NoAWaves} also hold for $\{\alpha',\gamma'\}$. Hence for the rectangles $R_\alpha$ 
and $R_\gamma$ we have both 
\begin{enumerate} 
\item $(\alpha'\cup\gamma')\cap R_\alpha\ne\emptyset$ and each outermost arc of $(\alpha'\cup\gamma')\cap R_\alpha$ 
is in $\alpha'$ and

\item $(\alpha'\cup\gamma')\cap R_\gamma\ne\emptyset$ and each outermost arc of $(\alpha'\cup\gamma')\cap R_\gamma$ 
is in $\gamma'$. 
\end{enumerate} 

\begin{claim}\label{cl:one-arc}
If $(\alpha'\cup\gamma')\cap R_\gamma$ (respectively $(\alpha'\cup\gamma')\cap R_\alpha$) contains more 
than one arc before the wave move, then there is more than one arc in $R_\gamma$ (respectively $R_\alpha$), after the 
wave move. 	

If $R_\gamma$  contains exactly one $\gamma'$-arc (respectively, $R_\alpha$ contains exactly one $\alpha'$-arc) 
before the wave move and if the next wave move is with respect to $\gamma'$ (respectively $\alpha'$), then after this 
wave move, $R_\gamma$ (respectively $R_\alpha$) contains more than one arc. 
\end{claim}

\begin{proof}[Proof of Claim~\ref{cl:one-arc}]
Denote the wave by $\eta''$. It follows from Claims~\ref{cl:NoAWaves} and \ref{cl:junction-arc} that, similar 
to the first wave $\eta$,  the curve $\partial_+P$ divides $\eta''$ into a colleciton of vertical arcs in $\mathcal{A}_l$ 
and a pair of junction arcs. As in the argument for $\eta$, one of the junction arcs is at a junction in $\mathcal{A}_r$. 
Again,  focus on how curves are changed by this wave move along $\eta''$ near this junction.

First, assume that  $(\alpha'\cup\gamma')\cap R_\gamma$ contains more than one arc.  Then  $R_\gamma$ has 
two outermost arcs and both are $\gamma'$-arcs. As before, perform the next wave move along $\eta''$ in 
two steps: 

If the wave is with respect to $\alpha'$, then the surgery step of the wave move changes 
$\alpha'$ into two curves $\alpha_1$ and $\alpha_2$, where $\alpha_2$ is parallel to $\gamma'$.  
As illustrated in Figure~\ref{fig:Rgamma}(a), near the junction arc of $\eta''$ in $\mathcal{A}_r$, 
the new arc from the surgery step that lies in $R_\gamma$ must belong to $\alpha_2$ since $\alpha_2$ 
is parallel to $\gamma'$. So after  $\alpha_2$ is deleted in the second step of the wave move, 
the two outermost $\gamma'$-arcs  in $R_\gamma$ remain outermost arcs in $R_\gamma$. 
Hence $R_\gamma$ still has at least  two $\gamma'$-arcs as claimed. 

\begin{figure}[ht]
	\begin{overpic}[width=13cm]{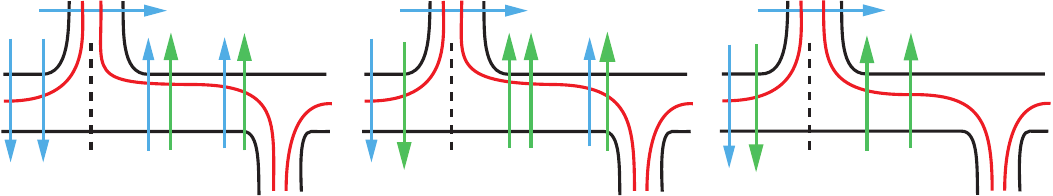}
		\put(15,-4){(a)}
		\put(50,-4){(b)}
		\put(84,-4){(c)}
		\put(7.5,2){$\eta''$}
		\put(-0.5,1){$\alpha'$}
		\put(33.5,1){$\alpha'$}
		\put(68,1){$\alpha'$}
		\put(17,7.5){$R_\gamma$}
		\put(51.5,7.5){$R_\gamma$}
		\put(3.5,1.5){$\alpha_1$}
		\put(12,2.5){$\alpha_2$}
		\put(37.5,0.5){$\gamma_2$}
		\put(47,3){$\gamma_1$}
		\put(45,8){$\kappa_\gamma$}
		\put(50,2.5){$\gamma'$}
		\put(16,2.5){$\gamma'$}
		\put(24,13){$\gamma'$}
		\put(58.5,13){$\gamma'$}
		\put(79,8){$\kappa_\gamma$}
		\put(86.8,7.3){$\mu$}
	\end{overpic}
	\vspace{8pt}
	\caption{Local pictures of a wave move at a junction in $\mathcal{A}_r$}
	\label{fig:Rgamma}
\end{figure}

If the wave $\eta''$ is with respect to $\gamma'$, then the surgery step of the wave move changes 
$\gamma'$ into two curve $\gamma_1$ and $\gamma_2$, where $\gamma_2$ is parallel to 
$\alpha'$ and will be removed in the second step of the wave move. Similar to the argument above 
and as illustrated in Figure~\ref{fig:Rgamma}(b),  two new vertical arcs are created near the 
junction arc of the wave by the surgery step, one in $R_\alpha$ and the other in 
$R_\gamma$. Since $\gamma_2$ is parallel to $\alpha'$, the new arc in $R_\alpha$ is in 
$\gamma_2$ and the new arc in $R_\gamma$ is in $\gamma_1$, see Figure~\ref{fig:Rgamma}(b). 

Let $\kappa_\gamma$ denote  this new vertical arc in $R_\gamma$  next to the junction, see 
Figure~\ref{fig:Rgamma}(b). The second step of the wave move deletes $\gamma_2$ and the 
curve $\gamma_1$ becomes the new $\gamma'$-curve. So $\kappa_\gamma$ becomes a 
new outermost $\gamma'$-arc in $R_\gamma$.  

Note that $\kappa_\gamma$ cannot be the only arc left because the deleted curve $\gamma_2$
 is parallel to $\alpha'$, which means that if an arc in $R_\gamma$ is deleted in the second step, 
there is at least one $\alpha'$-arc in $R_\gamma$. Hence there are at least two vertical arcs in 
$R_\gamma$ after the wave move on $\eta''$. In fact, by the symmetry from  $\pi$, both outermost 
arcs in $R_\gamma$ must be $\gamma'$-arcs after the wave move.

For the second part of the claim: In the argument above, suppose $(\alpha'\cup\gamma')\cap R_\gamma$ contains 
only one arc, denote it by $\mu$. So $\mu\subset\gamma'$. As the next wave move is with respect to $\gamma'$, 
after the surgery step of the wave move, $R_\gamma$ contains two arcs $\kappa_\gamma$ and $\mu$, 
see Figure~\ref{fig:Rgamma}(c).  Neither $\kappa_\gamma$ nor $\mu$ will be removed in the second 
step of the wave move because there is no $\alpha'$-arc in $R_\gamma$ and the removed curve is 
parallel to $\alpha'$.  A  similar argument applies to the $R_\alpha$ and $\alpha'$ case. This proves 
the latter part of the claim. 
\end{proof}

Recall that  a sequence of wave moves is performed on $\{\alpha',\gamma'\}$ until 
$(\alpha'\cup\gamma')\cap(\mathcal{R}^d\cup\mathcal{R}^u)=\emptyset$. 
Once that is achieved,  consider the Heegaard diagram given by $\{\delta,\varepsilon'\}$ 
and $\{\alpha',\gamma'\}$. 

As illustrated in Figure~\ref{fig:OneTwo2}, the induced orientations  (from a  fixed orientation on 
$\partial_+P$) on the two arcs in $\partial\widehat{\Sigma}\cap\partial\mathcal{A}_l$ are the same along the 
annulus $\mathcal{A}_l$. Hence the two endpoints of each arc of  $(\alpha'\cup\gamma')\cap\mathcal{A}_l$ 
represent intersection points of $\alpha'\cup\gamma'$ and $\partial_+P$ with the same sign. Similarly, the 
two endpoints of each arc of  $(\alpha'\cup\gamma')\cap\mathcal{A}_r$ represent points of 
$(\alpha'\cup\gamma')\cap\partial_+P$ with the same sign.  Since $\alpha'\cup\gamma'$ does not intersect 
$\mathcal{R}^d\cup\mathcal{R}^u$ at this stage, this implies that any two points of $\alpha'\cap\partial_+P$ 
adjacent along $\alpha'$ must have the same sign. Thus the intersection points of $\alpha'\cap\partial_+P$ 
all have the same sign. Similarly, the intersection points of $\gamma'\cap\partial_+P$ all have the same sign.

As claimed above, after the first wave move along $\eta$, $R_\alpha$ contains at least one 
$\alpha'$-arc and $R_\gamma$ contains at least one $\gamma'$-arc. Since the intersection points of 
$\delta\cap (\alpha\cup\gamma)$ all have the same sign before the wave moves, it follows from the 
wave-move operation that intersection points of $\delta\cap (\alpha'\cup\gamma')$ all have the same sign. Hence 
the $\alpha'$-arc and $\gamma'$-arc next to the cusp in $\mathcal{A}_r$ must have opposite directions with 
respect to the core curve of $\mathcal{A}_r$, see Figure~\ref{fig:OneTwo2} and Figure~\ref{fig:wave-move}(c). 
This means that the signs of $\alpha'\cap\partial_+P$ and $\gamma'\cap\partial_+P$ are opposite. 
Moreover, since the intersection points of $\delta$ with $\alpha'\cup\gamma'$ have the same sign, 
all the arcs of $R_\alpha\cap(\alpha'\cup\gamma')$ must be $\alpha'$-arcs and all the arcs in
$R_\gamma\cap(\alpha'\cup\gamma')$ must be $\gamma'$-arcs, at this stage.

\begin{figure}[ht]
	\begin{overpic}[width=8cm]{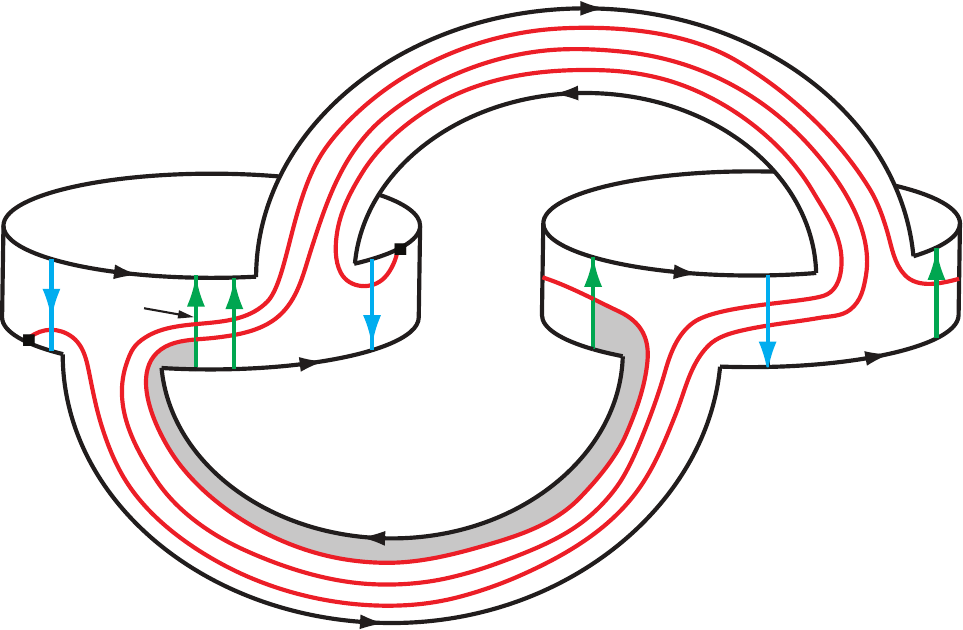}
		\put(32,33){$\delta$}
		\put(-3,24.5){$X^d$}
		\put(37.5,25.5){$\gamma'$}
		\put(78,23){$\gamma'$}
		\put(4,40){$\gamma'$}
		\put(12,32.5){$\kappa$}
		\put(23,23){$\alpha'$}
		\put(60,40){$\alpha'$}
		\put(96,26.5){$\alpha'$}
		\put(62.5,29.5){{\scriptsize $T$}}
	\end{overpic}
	\caption{Intersection of $\alpha'$ and $\gamma'$ with $\delta$ and $\widehat{\Sigma}$ after the wave moves}
	\label{fig:OneTwo2}
\end{figure}

Next, consider $\mathcal{A}_l\cap(\alpha'\cup\gamma')$. Suppose there is a component $\kappa$ 
of $\mathcal{A}_l\cap(\alpha'\cup\gamma')$ that meets $\delta$, see Figure~\ref{fig:OneTwo2}. 
The $\delta$-curve from the $\kappa$ arc to $\mathcal{A}_r$ along $\partial_+P$ determines a 
rectangle $T$,  depicted as the shaded region in Figure~\ref{fig:OneTwo2}.  Two opposite edges 
of $\partial T$ are subarcs of $\kappa$ and $(\alpha'\cup\gamma')\cap\mathcal{A}_r$. Another edge  
of $\partial T$ is a subarc of $\delta$, and the fourth edge  is a subarc of $\partial\widehat{\Sigma}$. 

Denote by $\kappa'$ the arc of $(\alpha'\cup\gamma')\cap\mathcal{A}_r$ containing the edge of 
$\partial T$ opposite to $\kappa$. Since the intersection points of $\delta$ with $\alpha'\cup\gamma'$ 
have the same sign, the $\kappa$ and $\kappa'$ must have compatible orientation.  Since the 
signs of $\alpha'\cap\partial_+P$ and $\gamma'\cap\partial_+P$ are opposite, the arcs $\kappa$ and 
$\kappa'$ are either both $\alpha'$-arcs or both $\gamma'$-arcs, see Figure~\ref{fig:OneTwo2}.  
For the same reason, any arc of $(\alpha'\cup\gamma')\cap\widehat{\Sigma}$ that intersects the 
rectangle $T$ must also belong to the same curve $\alpha'$ or $\gamma'$ as $\kappa$ and $\kappa'$.  

The discussion above shows  that the subarc of $\delta$ in $\partial T$ determines an 
$[\alpha'^-,\alpha'^+]$ blocking-edge if $\kappa$ is an $\alpha'$-arc and a 
$[\gamma'^-,\gamma'^+]$ blocking-edge if $\kappa$ is a $\gamma'$-arc.

\begin{claim}\label{cl:AlphaGamma}
Either the knot $K$ is doubly primitive and Proposition~\ref{pro:OneRTwoL} holds  or 
$(\alpha'\cup\gamma')\cap\widehat{\Sigma}$ has two $\alpha'$-arcs and two 
$\gamma'$-arcs that meet $\delta$ and the Heegaard diagram contains both
 $[\alpha'^-,\alpha'^+]$ and $[\gamma'^-,\gamma'^+]$ blocking-edges.
\end{claim}

\begin{proof}[Proof of Claim~\ref{cl:AlphaGamma}]
 Suppose  that $\delta$ contains no $[\gamma'^-,\gamma'^+]$ blocking-edge. The goal is to show that $K$ 
 must be doubly primitive. The proof for the case that $\delta$ contains no  $[\alpha'^-,\alpha'^+]$ blocking-edge is similar. 

Since there is no $[\gamma'^-,\gamma'^+]$ blocking-edge, the argument before the claim implies 
that no arc of $\mathcal{A}_l\cap(\alpha'\cup\gamma')$ that meets $\delta$ can be a 
$\gamma'$-arc.  
Furthermore, if 
$R_\gamma\cap(\alpha'\cup\gamma')$ contains more than one arc, since all the arcs 
of $R_\gamma\cap(\alpha'\cup\gamma')$ are $\gamma'$-arcs, two $\gamma'$-arcs in 
$R_\gamma$ are connected by a subarc of $\delta$ in $R_\gamma$,  which is a 
$[\gamma'^-,\gamma'^+]$-edge and this contradicts the assumption. 
Thus, $R_\gamma$ contains exactly one $\gamma'$-arc.

Since the arcs in $R_\alpha\cap(\alpha'\cup\gamma')$ consist of $\alpha'$-arcs and since no 
arc of $\mathcal{A}_l\cap(\alpha'\cup\gamma')$  that meets $\delta$ can be a $\gamma'$-arc, 
one concludes that exactly one arc in $\gamma'\cap\widehat{\Sigma}$  
intersects $\delta$.

Recall that, after the first wave move, there is at least one arc in $R_\gamma$.  
By Claim~\ref{cl:one-arc}, at any stage of the sequence of wave moves, if $R_\gamma$ contains 
more than one arc, then $R_\gamma$ will contain more than one arc regardless of the type of wave 
move. Since $R_\gamma$ contains exactly one arc in the end of the process, $R_\gamma$ must 
contain exactly one arc after each wave move in the sequence. Furthermore, if any wave move in 
the sequence,  after the first wave move, is with respect to $\gamma'$, then  by the second part of 
Claim~\ref{cl:one-arc}, $R_\gamma$ will contain more than one arc after that wave move, contradicting 
the fact that at the last stage there is a single arc. By Lemma~\ref{lem:FirstWave} the first wave 
 is with respect to $\gamma$. Thus the only possible way to have only one arc in $R_\gamma$ is
that

\begin{enumerate} [(a)] 
\item$R_\gamma$ contains only one arc after the first wave move on $\gamma$ and
\vskip5pt
\item all the subsequent wave moves are with respect to the $\alpha$-curve.
\end{enumerate}

We conclude that  the curve $\gamma'$ will not change  and stay fixed after the first wave move. 
To emphasize this fact, let $\beta$  denote the image of  $\gamma'$ after the first wave move.
Restating the conclusion above: \emph{There is exactly one arc of $\beta\cap\widehat{\Sigma}$ that 
intersects $\delta$}.

\vskip5pt

Since the curve $\beta$ is obtained by a wave move on $\gamma$, we can also obtain $\beta$ 
by band-summing $\alpha$ with $\gamma$.  Since $\gamma=\partial C$ and $\alpha$ is the outer 
boundary of the annulus $A$, $\beta$ is the outer boundary of an annulus $A'$ obtained by 
band-summing  $A$ with the disk $C$. There is exactly one arc of $\beta\cap\widehat{\Sigma}$ 
that  intersects $\delta$. It now follows from Definition \ref{def:complexity} that  $c_0(P, D, \beta) = 1$ 
and thus, by  Proposition~\ref{pro:Tao}, $K$ is doubly primitive.

The remaining case is that $\delta$ contains no $[\alpha'^-,\alpha'^+]$-edge at the end of the 
wave-move process. Then by the argument above, we have:

\begin{enumerate}
\item There is exactly one arc of $\alpha'\cap\widehat{\Sigma}$ that 
intersects $\delta$, and 

\vskip5pt

\item $R_\alpha$ contains only one arc after the first wave move and all the 
subsequent wave moves are with respect to the $\gamma'$. 
\end{enumerate}

Since the first wave move is on $\gamma$, this means that none of the wave moves in the 
sequence is on $\alpha$ and the curve $\alpha$ is never changed.  So there is exactly one 
arc of $\alpha\cap\widehat{\Sigma}$ that intersects $\delta$. As before, this means that 
$c_0(P, D, \alpha) = 1$ and it follows from  Proposition~\ref{pro:Tao} that $K$ is doubly primitive. 
\end{proof}

\begin{claim}\label{cl:DeltaEpsilon}
Suppose $K$ is not doubly primitive. Then, after the sequence of wave moves, there is no wave with respect to $\{\delta,\varepsilon'\}$.  
\end{claim}

\begin{proof}[Proof of Claim~\ref{cl:DeltaEpsilon}]
First consider the Heegaard diagram formed by curves $\{\delta, \varepsilon'\}$ and $\{\alpha',\gamma'\}$ 
before the sequence of wave moves has finished. So, at this stage, the curves $\alpha'\cup\gamma'$ 
still intersect the rectangles $\mathcal{R}^u$ and $\mathcal{R}^d$. As in the argument above and 
similar to  Claim~\ref{cl:Vwave}, this implies that $\alpha'\cup\gamma'$ contains at least two  
$[\delta^-,\delta^+]$-edges and two $[\varepsilon'^-,\varepsilon'^+]$-edges. Thus, the Whitehead graph 
$\Gamma(\{\delta,\varepsilon'\})$ is of type (i) in Figure~\ref{graphs} with $c\ge 2$, $d\ge 2$. If $a$ and $b$ 
in Figure~\ref{graphs}  (i) are both zero, then the Heegaard diagram is a standard Heegaard diagram of 
$L(c,p)\# L(d,q)$. Since $c\ge 2$ and $d\ge 2$, this contradicts our hypothesis that $M$ is either $S^3$ or 
$(S^2\times S^1)\# L(s,t)$. Thus at least one of $a$ and $b$ is nonzero and in particular, the Whitehead 
graph $\Gamma(\{\delta,\varepsilon'\})$ is connected. 

Since the intersection points of $\delta$ with $\alpha'$ and $\gamma'$ all have the same sign, 
no subarc of $\delta$ can be a wave with respect to $\{\alpha',\gamma'\}$.  It follows from 
Lemma~\ref{lem:dual} that, if no subarc of $\varepsilon'$ is a wave with respect to 
$\{\alpha',\gamma'\}$, then the Heegaard diagram is as in Figure~\ref{fig:octagon}(a), and if a subarc 
of $\varepsilon'$ is a wave with respect to $\{\alpha',\gamma'\}$ then the Heegaard diagram is as 
in Figure~\ref{fig:HeegDiagram}(c).

\begin{figure}[ht]
	\begin{overpic}[width=11cm]{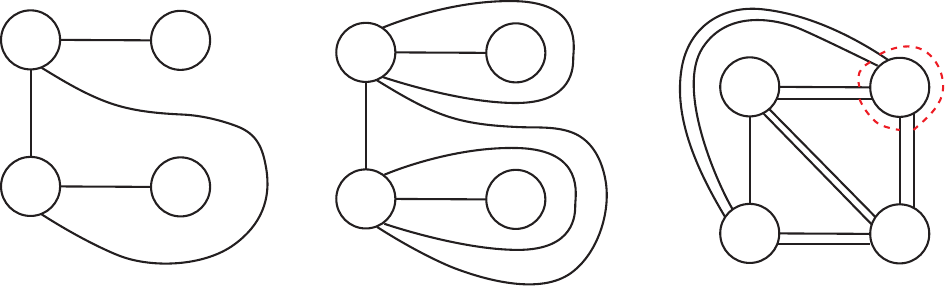}
	\put(10, -3){(a)}
\put(47, -3){(b)}
\put(85, -3){(c)}
		\put(11,26.5){$a'$}
		\put(46,25){$a'$}
		\put(11,11){$a'$}
		\put(47,9.5){$a'$}
		\put(43,13.3){$b'$}
		\put(43,28.8){$b'$}
		\put(0.5,16){$c'$}
		\put(36,15.5){$c'$}
		\put(21,18.5){$d'$}
		\put(61,16.5){$d'$}
		\put(1,25){$v_1^+$}
		\put(1,9){$v_1^-$}
		\put(36.5,23.5){$v_1^+$}
		\put(36.5,8){$v_1^-$}
		\put(77.5,4){$\delta^-$}
		\put(77.5,19.5){$\delta^+$}
		\put(93,19.5){$\varepsilon'^-$}
		\put(93,4){$\varepsilon'^+$}
	\end{overpic}
	\vspace{6pt}
	\caption{Whitehead graphs of $\Gamma(\{\delta,\varepsilon'\})$}\label{fig:HeegDiagram}
\end{figure}

Next,  consider the Heegaard diagram determined by  $\{\delta, \varepsilon'\}$ and $\{\alpha',\gamma'\}$ after
the sequence of wave moves has finished. 
Assume in contradiction that the claim is false and there is a wave with respect to $\{\delta,\varepsilon'\}$.

Since $K$ is assumed not to be doubly primitive, by Claim~\ref{cl:AlphaGamma}, the curve $\delta$ contains both $[\alpha'^-,\alpha'^+]$ and $[\gamma'^-,\gamma'^+]$ 
blocking-edges. Since $\delta$ meets both $\alpha'$ and $\gamma'$ and since the intersection points of $\delta$ with 
$\alpha'$ and $\gamma'$ all have the same sign, there are subarcs of $\delta$ connecting $\alpha'^+$ to $\gamma'^-$ 
and $\gamma'^+$ to $\alpha'^-$. This means that the Whitehead graph $\Gamma(\{\delta,\varepsilon'\})$ is still connected. 

Since there is a wave with respect to $\{\delta,\varepsilon'\}$, it follows from  Lemma~\ref{lem:dual} that 
$\Gamma(\{\delta,\varepsilon'\})$ must be as in  Figure~\ref{fig:HeegDiagram}(a) or (b), where $v_1$ 
denotes either $\delta$ or $\varepsilon'$ and all the labels $a'$, $b'$, $c'$, $d'$ are nonzero.   In both 
cases, there are two types  of arcs (i.e.~edges marked $c'$ and $d'$ in Figure~\ref{fig:HeegDiagram}(a, b)) 
connecting $v_1^+$ to $v_1^-$, where $v_1$ is either $\delta$ or $\varepsilon'$. 

This means that the last wave move (in the sequence of wave moves on $\{\alpha',\gamma'\}$) changes the configuration 
of $\Gamma(\{\delta,\varepsilon'\})$ from Figure~\ref{fig:octagon}(a) or Figure~\ref{fig:HeegDiagram}(c) 
to Figure~\ref{fig:HeegDiagram}(a) or (b). We will show next that it is impossible to make such a change 
on the configurations of $\Gamma(\{\delta,\varepsilon'\})$ via wave moves on $\{\alpha',\gamma'\}$. 
There are  two possible cases: 

\vspace{5pt}

\noindent Case (1). Before the last wave move, a subarc of $\varepsilon'$ is a wave with respect to 
$\{\alpha',\gamma'\}$, i.e.,  $\Gamma(\{\delta,\varepsilon'\})$ is as shown in Figure~\ref{fig:HeegDiagram}(c), 
see Lemma~\ref{lem:dual}.

\vspace{5pt}

In this case, by Lemma~\ref{lem:dual},  $\Gamma(\{\alpha',\gamma'\})$ is the right picture of Figure~\ref{dual2}. 
In particular, $\Sigma\ssm(\alpha'\cup\gamma'\cup\delta'\cup\varepsilon')$ consists of four hexagons and a 
collection of quadrilaterals. Since the intersection points of $\delta$ with $\alpha'$ and $\gamma'$ all have 
the same sign, no subarc of $\delta$ can be a wave.  So, the edge labelled $b'$ in Figure~\ref{dual2} 
denotes $b'$ copies of $\varepsilon'$-arcs. Moreover,  all the waves with respect to $\{\alpha',\gamma'\}$ 
are parallel to these $\varepsilon'$-arcs.

Let $\zeta_1$ and $\zeta_2$ be the two outermost arcs among these $b'$ parallel $\varepsilon'$-arcs.
This means that $\zeta_i\subset\varepsilon'$ ($i=1, 2$) and each $\zeta_i$  is an edge of a hexagonal 
region of  $\Sigma\ssm(\alpha'\cup\gamma'\cup\delta'\cup\varepsilon')$.  
Now we view $\zeta_i$ in the Whitehead graph $\Gamma(\{\delta,\varepsilon'\})$.
The three dashed red arcs in
 Figure~\ref{fig:HeegDiagram}(c) are possible pictures of waves parallel and next to $\zeta_i$.

The new meridian after a wave move along $\zeta_i$ can be obtained by a surgery using one of the 
three dashed arcs in Figure~\ref{fig:HeegDiagram}(c). The effect of this surgery on the configuration of 
$\Gamma(\{\delta,\varepsilon'\})$ is merging the two edges of a hexagonal region of  
$\Sigma\ssm(\alpha'\cup\gamma'\cup\delta'\cup\varepsilon')$ that are adjacent to $\zeta_i$ into an 
arc parallel to the edge opposite to $\zeta_i$ in the hexagon.  In particular, the  wave move along $\zeta_i$ does 
not create any new edge type in  $\Gamma(\{\delta,\varepsilon'\})$. However, Figure~\ref{fig:HeegDiagram}(a or b) 
has two types of edges (i.e.~edges marked $c'$ and $d'$ in Figure~\ref{fig:HeegDiagram}(a, b)) connecting $v_1^+$ 
to $v_1^-$ but there is only one type of $[\delta^+,\delta^-]$- or $[\varepsilon'^+,\varepsilon'^-]$-edge in 
Figure~\ref{fig:HeegDiagram}(c).  As the wave move does not create any new edge type in 
$\Gamma(\{\delta,\varepsilon'\})$,  the wave move can never change the configuration of 
$\Gamma(\{\delta,\varepsilon'\})$ from Figure~\ref{fig:HeegDiagram}(c) to Figure~\ref{fig:HeegDiagram}(a or b)

\vspace{5pt}

\noindent Case (2). No subarc of $\varepsilon'$ is a wave with 
respect to $\{\alpha',\gamma'\}$ before the last wave move.
\vspace{5pt}

By  Lemma~\ref{lem:Compatibility}, each wave in the sequence of wave moves must connect a pair of 
opposite sides of an octagaon, see the wave $\eta$ in Figure~\ref{fig:octagon}(b).  Similar to Case (1), 
the effect of this wave move on the configuration of $\Gamma(\{\delta,\varepsilon'\})$ is merging the two 
opposite edges of an octagon region of  $\Sigma\ssm(\alpha'\cup\gamma'\cup\delta'\cup\varepsilon')$ 
into another edge of this octagon. In particular, the  wave move does not create any new edge type 
in $\Gamma(\{\delta,\varepsilon'\})$. Thus, similar to Case (1), the wave move can never change the configuration of 
$\Gamma(\{\delta,\varepsilon'\})$ from Figure~\ref{fig:octagon}(a) to Figure~\ref{fig:HeegDiagram}(a or b).

Therefore, in both cases, $\Gamma(\{\delta,\varepsilon'\})$ cannot be changed to 
Figure~\ref{fig:HeegDiagram}(a or b) by the sequence of wave moves and this is a 
contradiction which proves the claim.
\end{proof}

By Lemma~\ref{lem:+ and -}, Claim~\ref{cl:AlphaGamma} and Claim~\ref{cl:DeltaEpsilon} we 
conclude  that there is no wave in the Heegaard diagram after these wave moves. Since 
the intersection pattern of the curves rules out a standard Heegaard splitting, this is a 
contradiction  to Theorems~\ref{thm:HOT} and \ref{thm:NeOk}.  The only possibility  left
by  Claim ~\ref{cl:AlphaGamma} is that $K$ is doubly primitive which finishes the proof
of the proposition.  
\end{proof}


\section{$\delta$ takes four short paths}\label{sec:AllPaths}

This section deals with Configuration 3 where $\delta$ takes two short paths in each annulus. 

\begin{proposition}\label{pro:AllPaths}
If the curve $\delta$ takes  two short paths in each of the annuli, then $K$ is doubly primitive.
\end{proposition}

\begin{proof}
The train track $\tau_D$ corresponding to this case is obtained from the train track in Proposition~\ref{pro:OneRTwoL} 
by attaching another segment. Thus,  as explained at the beginning of Chapter~\ref{cpt:ObtainingTheContradiction}, 
the complement $\Sigma\ssm\NN(\tau_D)$ is a disk with six cusps on its boundary. So one can think of 
$\Sigma\ssm\NN(\tau_D)$ as a hexagon.  

As in the  Section \ref{sec:OneArTwoAl}, the two cusps at $\partial\rho_x$ are called the \emph{$\rho$-cusps} 
and  the other four cusps called the \emph{$j$-cusps}. The $j$-cusps are located at the four junctions of 
$\mathcal{A}_l$ and $\mathcal{A}_r$ with the cusp directions pointing into the rectangles $\mathcal{R}^d$ 
and $\mathcal{R}^u$. The configuration in each annulus is similar to the picture of the right annulus in 
Figure~\ref{fig:OneTwo1}.  

The train track $\tau_D$ has a special segment $\rho_x$ containing the intersection point $X$ and the 
weight of $\delta$ at $\rho_x$ is one (as in  Section~\ref{sec:OneArTwoAl}).  As before, assume that 
$\rho_x$ and the hexagon $\Sigma\ssm\NN(\tau_D)$ are invariant under the involution $\pi$. Similar to 
Section~\ref{sec:OneArTwoAl}, the complement of $\tau_D\ssm\rho_x$ is an annulus as in shown 
Figure~\ref{fig:Annulus1} and Figure~\ref{fig:AllPathConfig}(a), and $\partial_+P$  must be a core curve of 
this annulus. Since $\tau_D$ and $\partial_+P$ are invariant under $\pi$ the intersection point 
$X=\rho_x\cap\partial_+P$ is a fixed point of $\pi$ and hence $\pi$ interchanges the two $\rho$-cusps, 
see Figure~\ref{fig:AllPathConfig}(a).  Moreover, as described in Remark~\ref{rem:Symmetries}, the 
involution $\pi$ interchanges the two junctions in $\mathcal{A}_l$ and the two junctions in $\mathcal{A}_r$.

As indicated by Figure~\ref{fig:AllPathConfig}(a), the two sides of a small neighborhood  $\NN(\rho_x)$ 
are parts of a pair of opposite edges of the hexagon $\Sigma\ssm\NN(\tau_D)$. Hence the action 
of $\pi$ on the hexagon $\Sigma\ssm\NN(\tau_D)$ interchanges these opposite edges.  
Moreover, as $\partial_+P$ is invariant under $\pi$, the arc $\partial_+P\cap (\Sigma\ssm\NN(\tau_D))$ 
connects a pair of opposite boundary edges of the hexagon $\Sigma\ssm\NN(\tau_D)$.

\begin{figure}[ht]
	\begin{overpic}[width=8cm]{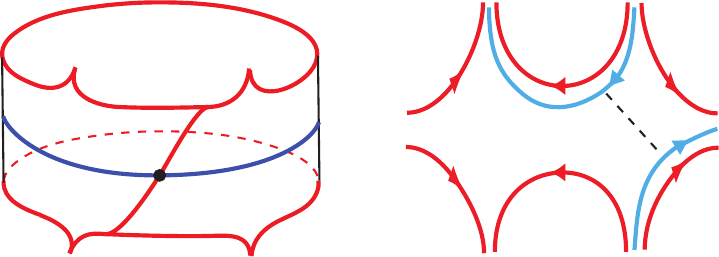}
		\put(21,-5){(a)}
		\put(78,-5){(b)}
		\put(12,22.5){$\tau_D$}
		\put(20,6){$\rho_x$}
		\put(19,12.5){{\small $X$}}
		\put(32,8){{\small $\partial_+P$}}
		\put(63,20){$\delta$}
		\put(78,17.5){$\varepsilon$}
	\end{overpic}
	\vspace{8pt}
	\caption{Configurations of $\Sigma\ssm\tau_D$ and $s$-wave for $\varepsilon$}
	\label{fig:AllPathConfig}
\end{figure}

In $\NN(\tau_D)$ there are three splitting arcs (defined in Definition~\ref{def:splitting-arc}) connecting 
the cusps of $\NN(\tau_D)$ in pairs which are disjoint from $\delta$. When $\NN(\tau_D)$ is cut open 
along these three splitting arcs, the resulting surface is a product neighborhood of $\delta$.  Denote the three 
splitting arcs by $s_\rho$, $s_l$ and $s_r$.
 
By connecting the endpoints of each splitting arc via an arc in the disk $\Sigma\ssm\NN(\tau_D)$, we obtain a 
nontrivial simple closed curve. Since every nontrivial simple closed curve is invariant under the involution $\pi$ after 
isotopy, the two cusps at the endpoints of each splitting  arc must be invariant under $\pi$. Since the involution 
interchanges the two $\rho$-cusps and interchanges the two junctions in each annulus $\mathcal{A}_l$ 
and $\mathcal{A}_r$, this implies that $s_\rho$ connects the two $\rho$-cusps, $s_l$ connects the two 
$j$-cusps at $\mathcal{A}_l$, and $s_r$ connects the two $j$-cusps at $\mathcal{A}_r$.

As the cusp directions of the $j$-cusps point into the two rectangles $\mathcal{R}^d$ and  $\mathcal{R}^u$, 
both $s_l$ and $s_r$ contain subarcs that are core arcs of $\mathcal{R}^d$ and $\mathcal{R}^u$.  Similar 
to sections~\ref{sec:OneArTwoAl}, this implies that each of $s_l$ and $s_r$ intersects $\gamma$  at least twice.

Now we consider $s_\rho$. Similar to the proof of Claim~\ref{cl:splitting-arc} in Section~\ref{sec:OneArTwoAl}, each end of $\rho_x$ is attached either to a long path or to a component of $\delta\cap(\mathcal{R}^d\cup\mathcal{R}^u)$. 
Since there is no long path in the current configuration, the two ends of $\rho_x$ are attached to two components of $\delta\cap(\mathcal{R}^d\cup\mathcal{R}^u)$ that are core arcs of $\mathcal{R}^d$ and $\mathcal{R}^u$ respectively.  As in the proof of Claim~\ref{cl:splitting-arc}, this implies that $s_\rho$ also contains subarcs that are core 
arcs of $\mathcal{R}^d$ and $\mathcal{R}^u$, and $s_\rho$ intersects $\gamma$ at least twice. 

By part (2) of Proposition~\ref{pro:Paths}, we may assume that the orientations of $\delta$ and $\varepsilon$ 
are compatible in $\NN(\tau_D)$, see the arrows of the arcs in in 
Figure~\ref{fig:AllPathConfig}(b).  The surface  $\Sigma\ssm\NN(\tau_D)$ is a disk with six cusps at its boundary. 
Similar to the $s$-wave illustrated in Figure~\ref{fig:Annulus1}, perform wave moves on $\varepsilon$ along 
all possible $s$-waves in $\Sigma\ssm\NN(\tau_D)$, (such $s$-waves are indicated by the dashed arc in 
Figure~\ref{fig:AllPathConfig}(b)). 

Notice that if there are two distinct $\varepsilon$-arcs in the disk $\Sigma\ssm\NN(\tau_D)$ that go into two 
adjacent cusps, then there is an $s$-wave connecting these two arcs, see  Figure~\ref{fig:AllPathConfig}(b). 
We perform all possible  wave moves along such $s$-waves on $\varepsilon$ and denote the resulting 
meridian of $W$ by $\varepsilon'$.

Suppose $\varepsilon'$ no longer admits any $s$-wave as in Figure~\ref{fig:AllPathConfig}(b). 
In particular, similar to the argument before Claim~\ref{cl:MoreThanOneArc}, $\varepsilon'$ does 
not pass through all cusps of $\NN(\tau_D)$.  Consider different types of $\varepsilon'$-arcs in the 
cusped disk $\Sigma\ssm\tau_D$. The induced action by $\pi$ on the disk $\Sigma\ssm\tau_D$ is a 
 $180^\circ$-rotation, so the $\varepsilon'$-arcs in the cusped disk $\Sigma\ssm\tau_D$ are symmetric 
 with respect to the rotation. Since $\varepsilon'$ admits no  $s$-wave in the cusped disk 
 $\Sigma\ssm\tau_D$, the $\varepsilon'$-arcs  in $\Sigma\ssm\tau_D$ have three possible configurations, 
 up to rotation and reflection, which are indicated in Figure~\ref{fig:AllPathCusp}: 

\begin{enumerate}
\item There are three types of $\varepsilon'$-arcs in $\Sigma\ssm\tau_D$: diagonals connecting a pair 
	of opposite cusps and arcs parallel to a pair of opposite edges next to the diagonal, see the blue 
	arcs in Figure~\ref{fig:AllPathCusp}(a), or	
		
\item All the $\varepsilon'$-arcs in $\Sigma\ssm\tau_D$ are parallel to a diagonal connecting a 
	pair of opposite cusps, see the blue arcs in Figure~\ref{fig:AllPathCusp}(b) (ignore the dashed arcs for now), or
		
\item The $\varepsilon'$-arcs in $\Sigma\ssm\tau_D$ are all parallel to a pair of opposite edges, 
see the blue arcs in Figure~\ref{fig:AllPathCusp}(c) (ignore the dashed arcs for now).
\end{enumerate}
 
 \begin{figure}[ht]
 	\begin{overpic}[width=11cm]{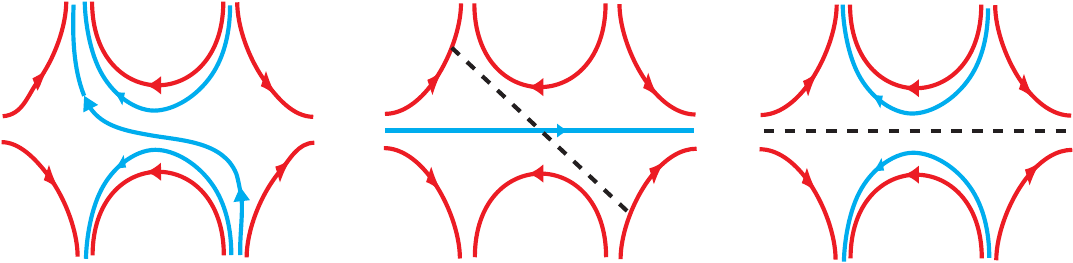}
 		\put(13,-5){(a)}
 		\put(48,-5){(b)}
 		\put(83,-5){(c)}
 		\put(27,15){$\delta$}
 		\put(76,13){$\varepsilon''$}
 		\put(18,11.5){$\varepsilon'$}
 		\put(90.5,15){$\varepsilon'$}
 		\put(91.5,5.8){$\varepsilon'$}
 		\put(41,13){$\varepsilon'$}
 	\end{overpic}
 	\vspace{10pt}
 	\caption{Possible configurations of $\varepsilon'$ in $\Sigma\ssm\NN(\tau_D)$}
 	\label{fig:AllPathCusp}
 \end{figure}

Now consider the Heegaard diagram given by $\{\delta,\varepsilon'\}$ and $\{\alpha,\gamma\}$. There are three 
cases to discuss. 

\vskip5pt

\noindent{\bf Case (a)}: Assume there is a cusp of $\NN(\tau_D)$ that contains two $\varepsilon'$-arcs. In other words, 
two components of $\varepsilon'\cap\NN(\tau_D)$ are parallel to the same splitting curve $s_\rho$, $s_l$ or $s_r$. 
\vskip5pt

Note that this case includes the configuration in Figure~\ref{fig:AllPathCusp}(a), since there are multiple curves going 
into a pair of opposite cusps and since each splitting arc $s_\rho$, $s_l$ or $s_r$ connects a pair of opposite cusps.

Similar to Claim~\ref{cl:Vwave}, 
since each $s_\rho$, $s_l$ or $s_r$ intersects $\gamma$ at least twice and since the orientations of arcs of 
$\varepsilon'\cap\NN(\tau_D)$ are compatible, $\gamma$ contains two distinct subarcs that are both 
$[\varepsilon'^-,\varepsilon'^+]$-edges.  As explained before Case (a), there is a cusp that contains no 
$\varepsilon'$-arcs.  This means that $\gamma$ also contains two distinct subarcs that are 
$[\delta^-,\delta^+]$-edges.  

Now the proof for Case (a) is identical to the proof in Section~\ref{sec:OneArTwoAl}. Recall that in Claim~\ref{cl:Vwave} 
of Section~\ref{sec:OneArTwoAl}, we have the same conclusion, which implies that there is no wave with respect to 
$\{\delta,\varepsilon'\}$ and hence there must a wave $\eta$ with respect to $\{\alpha,\gamma\}$. As in the 
argument between Claim~\ref{cl:Vwave} and Claim~\ref{cl:junction-arc}, we can isotope $\eta$ to intersect $\partial_+P$ 
efficiently.  In the  current setting, $\mathcal{A}_r$ also has two short paths, which is the same situation as in 
Section~\ref{sec:OneArTwoAl}. Thus Claims ~\ref{cl:junction-arc}--\ref{cl:DeltaEpsilon} also hold and hence 
Proposition~\ref{pro:AllPaths} holds in Case (a).

Note that the setting in this section is simpler than Section~\ref{sec:OneArTwoAl}. The annulus $\mathcal{A}_l$ 
has two short paths in this setting while $\mathcal{A}_l$ has only one short path in Section~\ref{sec:OneArTwoAl}. 
The additional short path in the current setting could be ignored and the proof can proceed as in
Section~\ref{sec:OneArTwoAl}.

Next, assume that Case (a) does not occur, i.e.,  each cusp of $\NN(\tau_D)$ has at most one 
$\varepsilon'$-arc. In particular, the first possible configuration of $\varepsilon'$-arcs in 
$\Sigma\ssm\tau_D$, i.e.~Figure~\ref{fig:AllPathCusp}(a) does not occur. Thus, it remains to 
discuss the second and third possible configurations, i.e.~Figure~\ref{fig:AllPathCusp}(b) 
and (c). 

\vskip5pt
\noindent{\bf Case (b)}: The configuration for $\varepsilon'\cap  (\Sigma\ssm\NN(\tau_D))$ is as in 
Figure~\ref{fig:AllPathCusp}(b). 
\vskip5pt

Since Case (a) does not occur, there is only a single $\varepsilon'$-arc in Figure~\ref{fig:AllPathCusp}(b) going into 
the cusp. So the intersection $\varepsilon'\cap  (\Sigma\ssm\NN(\tau_D))$ is a single arc forming a main 
diagonal of the hexagon $\Sigma\ssm\NN(\tau_D)$. 
We concluded earlier in the proof of Proposition~\ref{pro:AllPaths} that $\partial_+P\cap(\Sigma\ssm\NN(\tau_D))$ 
is an arc connecting a pair of opposite edges of $\Sigma\ssm\NN(\tau_D)$, indicated by the dashed arc in
Figure~\ref{fig:AllPathCusp}(b), also see Figure~\ref{fig:AllPathConfig}(a) for a picture of how $\partial_+P$ 
intersects $\Sigma\ssm\NN(\tau_D)$. As shown in Figure~\ref{fig:AllPathCusp}(b), this conclusion on $\partial_+P\cap(\Sigma\ssm\NN(\tau_D))$ implies that $\varepsilon'$ 
intersects $\partial_+P$ in exactly one point. Let $E'$ is the disk in $W$ bounded by $\varepsilon'$, 
then $(P, E')$ is a $(\mathcal{P},\mathcal{D})$-pair. The configuration of $\varepsilon'$ implies that 
$$|(\alpha\cup\gamma)\cap\varepsilon'|<|(\alpha\cup\gamma)\cap\delta|.$$ 
Since $\delta$ takes two paths in each annulus, every component of $(\alpha\cup\gamma)\cap\widehat{\Sigma}$ 
must intersect $\delta$ and
$$c_0(P,E',\alpha,\gamma)\le |(\alpha\cup\gamma)\cap\widehat{\Sigma}|= c_0(P,D,\alpha,\gamma).$$ 
Hence $c(P,E',\alpha,\gamma)< c(P,D,\alpha,\gamma)$,
contradicting Assumption~\ref{ass:minimality} regarding the choice of the pair $(P, D)$.

\vskip5pt
\noindent{\bf Case (c)}: The configuration for $\varepsilon'\cap  (\Sigma\ssm\NN(\tau_D))$ 
is as in  Figure~\ref{fig:AllPathCusp}(c).
\vskip5pt

In this case, we can use a main diagonal of the 
hexagon $\Sigma\ssm\NN(\tau_D)$, see the dashed 
arc in Figure~\ref{fig:AllPathCusp}(c), to connect the two ends of the splitting arc of $\NN(\tau_D)$ 
which connects (in $\NN(\tau_D)$) the two cusps that $\varepsilon'$ does not pass through.  This gives 
a nonseparating simple closed curve $\varepsilon''$ that is disjoint from both  $\delta$ and $\varepsilon'$,  
see the dashed arc in Figure~\ref{fig:AllPathCusp}(c) for a picture of $\varepsilon''$ within $\Sigma\ssm\NN(\tau_D)$. 
As $\{\delta,\varepsilon'\}$ is a complete set of meridians for $W$, $\varepsilon''$ must also bound a disk 
in $W$.  Now consider $\{\delta,\varepsilon''\}$. The configuration of $\varepsilon''$ is the same as that of 
$\varepsilon'$ in Case (b).  So we can apply the argument in Case (b) to $\varepsilon''$ and obtain the 
same contradiction to Assumption~\ref{ass:minimality}.
 
This finishes the proof of the proposition.
\end{proof}


\section{The curve $\delta$ takes a long path}\label{sec:NoLongPath}

In this section we consider Configuration 4, namely  when $\delta$ takes a long path in 
one (or both) of the annuli.  

\begin{proposition}\label{pro:LongPath}
If $\delta$ takes a long path in an annulus, then $K$ is doubly primitive
or there is an isotopy of the annulus which converts
the long path to a short path.
\end{proposition}

\begin{proof}
Without loss of generality, suppose $\delta$ takes a long path in $\mathcal{A}_r$.  

First note that in some  special configurations one can perform an isotopy on $\mathcal{A}_r$ so that $\delta$ 
only takes short paths in $\mathcal{A}_r$. For example, if $\delta$ only wraps around  $\mathcal{A}_r$ once and if 
$\alpha\cup\gamma$ does not intersect the shaded region in Figure~\ref{fig:longisotopy}, then after an 
isotopy in Figure~\ref{fig:longisotopy}, $\delta$ only takes short paths in $\mathcal{A}_r$. This isotopy 
is equivalent to choosing a different product structure for $\mathcal{A}_r$ provided that arcs of 
$(\alpha\cup\gamma)\cap\mathcal{A}_r$ remain vertical in the new product structure. So assume 
that $\delta$ takes a long path in $\mathcal{A}_r$ after all such  isotopies have been performed. 
In other words, assume that $\delta \cap \mathcal{A}_r$ (similarly $\delta \cap \mathcal{A}_l$ 
if it contains a long path) 
cannot be converted to 
short paths by such an isotopy. 

\begin{figure}[!ht]
	\begin{overpic}[width=10cm]{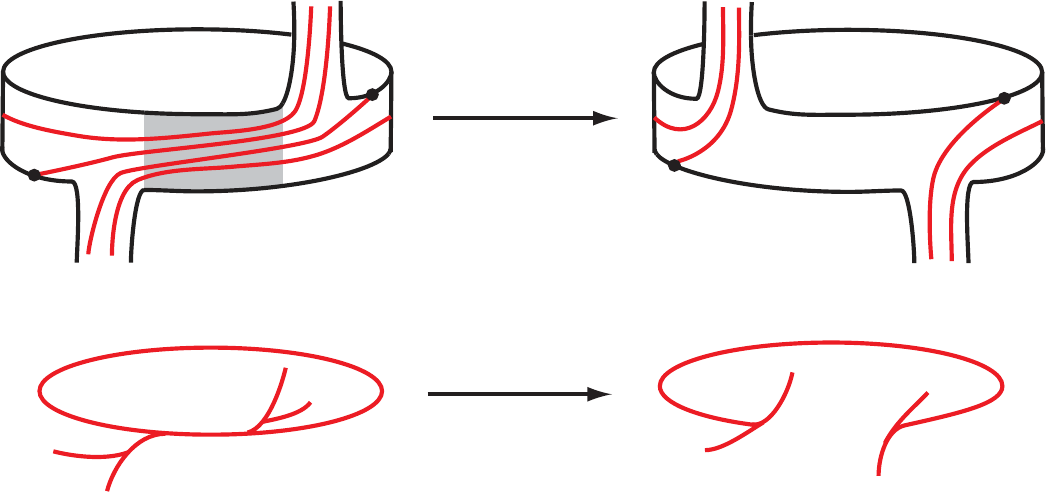}
		\put(43,38){isotopy}
		\put(41,11){splitting}
		\put(2,2){$\rho_x$}
		\put(30,9){$\rho_x$}
	\end{overpic}
	\caption{An isotopy converting a long path to a short path}
	\label{fig:longisotopy}
\end{figure}

\begin{claim}\label{cl:longpath}
Either Proposition~\ref{pro:LongPath} holds or there is a  meridian  $\varepsilon'$  of $W$, disjoint from $\delta$,
 such that there is no wave with respect to $\{\delta,\varepsilon'\}$ in the Heegaard diagram determined by 
 $\{\alpha,\gamma\}$ and  $\{\delta,\varepsilon'\}$.
\end{claim}

\begin{proof}[Proof of the Claim~\ref{cl:longpath}]
Consider the possible configurations of $\delta\cap\mathcal{A}_l$, i.e., the intersection of $\delta$ with 
the other annulus. There are two cases to consider:

\vskip5pt

\noindent \underline{Case (1)}. The curve $\delta$ either takes two short paths, or a long path in $\mathcal{A}_l$.  

Note that it is still assumed that  if $\delta$ takes a long path in $\mathcal{A}_l$, and there is no isotopy 
as in Figure~\ref{fig:longisotopy} that converts the long path into a short path. Furthermore, as explained 
at the beginning in Chapter~\ref{cpt:ObtainingTheContradiction}, the surface $\Sigma\ssm\tau_D$ in this case is a 
disk with six cusps along its boundary.  The discussion in this case is similar to  Section~\ref{sec:AllPaths}. 
For convenience here is a summary of the setup in Section~\ref{sec:AllPaths}:

There are three splitting curves in  $\NN(\tau_D)$  which are disjoint from $\delta$ and which connect the 
cusps of $\NN(\tau_D)$ in pairs.  When the surface $\NN(\tau_D)$ is cut open along these three splitting 
curves, the resulting surface is a product neighborhood of $\delta$. As in Section~\ref{sec:AllPaths},   each 
splitting arc is invariant under $\pi$. Hence the involution interchanges the two cusps at the endpoints of 
each splitting arc.  The involution interchanges the two $\rho$-cusps, and also interchanges the other two 
cusps in $\mathcal{A}_l$ and $\mathcal{A}_r$ respectively.  Thus a splitting arc, denoted by $s_\rho$, 
connects the two $\rho$-cusp.  The other two splitting arcs, denoted by $s_l$ and $s_r$, connect two 
cusps in $\mathcal{A}_l$ and two cusps in $\mathcal{A}_r$, respectively.

\vskip5pt

 The  \emph{difference} between the current setting and Section~\ref{sec:AllPaths} is that in the current setting
 the splitting arcs   $s_\rho$, $s_l$ and $s_r$ may not go through the two rectangles $\mathcal{R}^u$ and 
 $\mathcal{R}^d$. In fact, it is not immediately clear that the three splitting arcs  even intersect $\alpha\cup\gamma$.

\vskip5pt

First consider $s_\rho$. 
The proof of Claim~\ref{cl:splitting-arc} on $\rho_x$ implies that either there is a long path or $(\alpha\cup\gamma)\cap s_\rho\ne\emptyset$. 
Moreover, it follows from Definition~\ref{def:traintracks} that the only configuration where  $(\alpha\cup\gamma)\cap s_\rho=\emptyset$ is as shown in the left picture of Figure~\ref{fig:longisotopy}, where $s_\rho$ is a short arc cutting through the shaded region of Figure~\ref{fig:longisotopy} and the shaded region does not contain any $\alpha$- or $\gamma$-arc.  
Thus, if $(\alpha\cup\gamma)\cap s_\rho=\emptyset$, 
we can perform an isotopy as shown in Figure~\ref{fig:longisotopy}.  The isotopy can be viewed as a change of the product structure of the
annulus. Since $(\alpha\cup\gamma)\cap s_\rho=\emptyset$, $\alpha\cup\gamma$ does not 
intersect the shaded region in Figure~\ref{fig:longisotopy}. Hence, the isotopy does not affect 
$\alpha$ and $\gamma$. 
Note that the effect of this isotopy on the train track is a splitting along $s_\rho$, see Figure~\ref{fig:longisotopy}. 
This isotopy converts a long path to a short path, contradicting 
our assumption, before the claim, that there is no such isotopy. Thus, $(\alpha\cup\gamma)\cap s_\rho\ne\emptyset$.

Similarly, if $(\alpha\cup\gamma)\cap s_l=\emptyset$ (or $(\alpha\cup\gamma)\cap s_r=\emptyset$), then 
the configuration of  $s_l$ (or $s_r$) must be a short arc cutting through the shaded region of Figure~\ref{fig:longisotopy} and the shaded region does not contain any $\alpha$- or $\gamma$-arc. 
We can perform a 
similar isotopy (without changing $\alpha\cup\gamma$) which changes the long path into a short 
path, also contradicting the assumption that there is no such isotopy.

Therefore, $s_\rho$, $s_l$ and $s_r$ must all meet $\alpha\cup\gamma$.
Now the argument is the same as  Section~\ref{sec:AllPaths}.  Performing wave moves on $\varepsilon$ 
in the hexagon $\Sigma\ssm\tau_D$ along $s$-waves,  a meridian $\varepsilon'$ is obtained. As in 
Section~\ref{sec:AllPaths}, there are three possible configurations for the $\varepsilon'$-arcs in the 
hexagon $\Sigma\ssm\tau_D$. If the configuration is as in Case (a) in Section~\ref{sec:AllPaths}, since 
$s_\rho$, $s_l$ and $s_r$ all meet $\alpha\cup\gamma$, $\alpha\cup\gamma$ contains subarcs that are 
$[\delta^-,\delta^+]$ and $[\varepsilon'^-,\varepsilon'^+]$ blocking-edges. By Lemma~\ref{lem:+ and -}, 
there is no wave with respect to $\{\delta,\varepsilon'\}$ and the claim holds. If the configuration is 
as in Case  (b) and Case (c) in Section~\ref{sec:AllPaths}, then   Assumption~\ref{ass:minimality} is contradicted.

\vskip5pt
\noindent \underline{Case (2)}. The curve $\delta$ takes one short path in $\mathcal{A}_l$. 

In this case, as explained at the beginning of Chapter~\ref{cpt:ObtainingTheContradiction}, the surface
$\Sigma\ssm\NN(\tau_D)$  is an annulus with two cusps on each boundary curve. As in the argument 
before Claim~\ref{cl:MoreThanOneArc},  perform wave moves on $\varepsilon$ along $s$-waves 
as in Figure~\ref{fig:Annulus1} to  obtain a new meridian $\varepsilon'$ of $W$. In particular, the arcs 
in  $\varepsilon'\cap\NN(\tau_D)$, if there are any, are all parallel to one of the two splitting arcs of 
$\NN(\tau_D)$. Furthermore, Claim~\ref{cl:MoreThanOneArc} also holds for $\varepsilon'$ in the 
current setting, since $\NN(\tau_D)$ in Section~\ref{sec:OneArTwoAl} is of the same topological type. 
So the intersection $\varepsilon'\cap\NN(\tau_D)$ contains at least two arcs, all parallel to a splitting arc.

Since it is assumed that a long path cannot be converted to a short path, the argument in Case (1) 
implies that both splitting arcs of $\NN(\tau_D)$ must intersect $\alpha\cup\gamma$. Similar 
to Claim~\ref{cl:Vwave}, this implies that $\alpha\cup\gamma$ contains subarcs that are 
$[\delta^-, \delta^+]$ and $[\varepsilon'^-,\varepsilon'^+]$ blocking-edges.  By Lemma~\ref{lem:+ and -}, 
this means that there is no wave with respect to $\{\delta,\varepsilon'\}$ as claimed. 
\end{proof}

We are now in position to finish the proof of Proposition~\ref{pro:LongPath}. 
By Claim~\ref{cl:longpath}, we may assume that there is no wave with respect to $\{\delta,\varepsilon'\}$. By 
part (1) of Proposition~\ref{pro:Paths}, the intersection points of $\delta\cap\gamma$ all have the same sign. 
Next, we consider the core curve $\mathfrak{a}_r$ of $\mathcal{A}_r$.  Since $\delta$ takes a long 
path in $\mathcal{A}_r$, an arc in $\delta\cap\mathcal{A}_r$ wraps around the annulus at least once. 
So the core curve $\mathfrak{a}_r$ is carried by the train track $\tau_D$. Hence $\mathfrak{a}_r$ admits 
a direction induced from the orientation of $\delta$ and $\tau_D$. Since points of $\delta\cap\gamma$ all 
have the same sign, the intersection points of $\mathfrak{a}_r\cap\gamma$ all have the same sign. 

Recall that in our construction (see Remark~\ref{rem:core}), an arc in $\gamma\cap\mathcal{R}^u$ 
determines a $\partial$-compression disk for $P$. When  $P$  is $\partial$-compressed  along this disk,  
we obtain two planar surfaces $P_l$ and $P_r$ in the compression body $U$ with $\mathfrak{a}_l=\partial_+P_l$ 
and $\mathfrak{a}_r=\partial_+P_r$. 
Since $\gamma$ bounds a disk $C$ in $U$, if $C\cap P_r\ne\emptyset$, then the two endpoints of each arc in 
$C\cap P_r$ are points in $\gamma\cap\mathfrak{a}_r$ with  opposite signs of intersection,  a contradiction to
the conclusion that the intersection points of $\mathfrak{a}_r\cap\gamma$ all have the same sign. Thus 
$C\cap P_r=\emptyset$ and hence $\gamma\cap\mathfrak{a}_r=\emptyset$. By Lemma~\ref{lem:PCIntersection}, 
this means that $P_r$ must be a vertical annulus in the compression body $U$.   Moreover, since the surgery 
slope is an integer,  after isotopy, $P_r\cap A$ is a single vertical arc, where $A$ is the annulus with 
$\partial_+A=\alpha$. Hence $\alpha\cap\mathfrak{a}_r$ is a single point.  So we have  
$\mathcal{A}_r\cap\gamma=\emptyset$ and $\mathcal{A}_r\cap\alpha$ is a single arc.

Since there is no wave with respect to $\{\delta,\varepsilon'\}$, by Theorems~\ref{thm:HOT} and \ref{thm:NeOk}, 
there must be a wave with respect to $\{\alpha,\gamma\}$  and by Lemma~\ref{lem:FirstWave}, this wave must 
be with respect to $\gamma$.   If a component of $\delta\cap\mathcal{A}_r$ in the long path intersects the arc 
$\alpha\cap\mathcal{A}_r$ more than once, since $\mathcal{A}_r\cap\gamma=\emptyset$, a subarc of 
$\delta\cap\mathcal{A}_r$ must be a spiral around $\mathcal{A}_r$ connecting the plus-side of $\alpha$ to its 
minus-side. This means that a subarc of $\delta$ is an $[\alpha^+,\alpha^-]$ blocking-edge. However, by the proof of 
Lemma~\ref{lem:+ and -}, the existence of an $[\alpha^+,\alpha^-]$ blocking-edge means that there is no wave with 
respect to $\gamma$, a contradiction.  So each component of $\delta\cap\mathcal{A}_r$ in the long path 
intersects the arc $\alpha\cap\mathcal{A}_r$ only once.  This means that the shaded region in 
Figure~\ref{fig:longisotopy} does not contain any $\alpha$-arc. As $\mathcal{A}_r\cap\gamma=\emptyset$, 
the shaded region in Figure~\ref{fig:longisotopy} intersects neither $\alpha$ or $\gamma$.  Hence we can 
perform an isotopy as in Figure~\ref{fig:longisotopy},  changing the long path into a short path with respect 
to the new product structure. This contradicts our assumption at the beginning of the proof. Therefore, 
Proposition~\ref{pro:LongPath} holds.
\end{proof}


\section{$\delta$ takes one short path in each annulus, a special case}\label{sec:OneAndOne-S}

In this and the next section we rule out the remaining possible configuration, namely, when 
$\delta$  takes only two short paths: One in $\mathcal{A}_l$ and one in  $\mathcal{A}_r$. 
As $\delta$ takes \emph{only two short paths} it imposes the least restrictions on the possible 
configurations and hence this is the most complicated case. We therefore begin by considering a special 
case in this section and will deal with the general case in the following section. 
The configurations and proof for the special case are similar to those of Sections~\ref{sec:OneArTwoAl} 
and \ref{sec:NoLongPath}. 

As before, we assume $\delta$ and $\partial_+P$ are invariant under  $\pi$ 
and $X=\delta\cap\partial_+P$ is a fixed point of $\pi$.  By symmetry (2) in Remark~\ref{rem:Symmetries}, 
we may assume $X^d\cup X^u$ lies in either $\partial A_l$ or $\partial\mathcal{R}^d\cup\partial\mathcal{R}^u$.
Moreover, since symmetry  (3) in Remark~\ref{rem:Symmetries} interchanges the two short paths in 
$\mathcal{A}_l$, we may assume $\delta$ takes a fixed short path in $\mathcal{A}_l$.

Consider the train track $\tau_D$.  As explained at the beginning of this chapter, the shape of $\tau_D$ 
must be as shown in Figure~\ref{fig:traintrack}. So both $\NN(\tau_D)$ and $\Sigma\ssm\NN(\tau_D)$ 
are once-punctured tori. The train track $\tau_D$ has a special segment $\rho_x$ which contains the 
point $X=\delta\cap\partial_+P$ and the weight of $\delta$ at $\rho_x$ is one. 

\vskip5pt

\begin{figure}[!ht]
	\begin{overpic}[width=7cm]{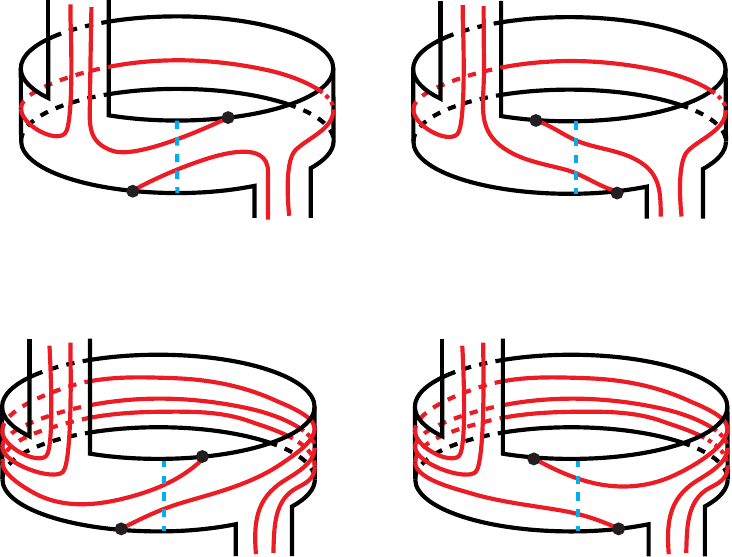}
		\put(19,38){(a)}
		\put(77,38){(b)}
		\put(19,-5){(c)}
		\put(77,-5){(d)}
	\end{overpic}
	\vspace{10pt}
	\caption{Possible configurations of $\rho_x$ and $\kappa$}
	\label{fig:double}
\end{figure}

The special case  that we discuss in this section is: 

\vskip5pt
\noindent\underline{Special Case}: 
There is a component $\kappa$ of $(\alpha\cup\gamma)\ssm\partial_+P$ such that 
$\kappa\cap \tau_D=\kappa\cap\rho_x$ contains more than one point, see Figure~\ref{fig:double} 
where the dashed arc is $\kappa$.

Consider the arc $\kappa$. By the construction of $\tau_D$,  the intersection $\kappa\cap\rho_x$ consists 
of exactly two points,  see Figure~\ref{fig:double}.   There are four possible configurations for 
$\delta\cap\mathcal{A}_l$ as shown in Figure~\ref{fig:double}.  

\begin{figure}[!ht]
	\begin{overpic}[width=7cm]{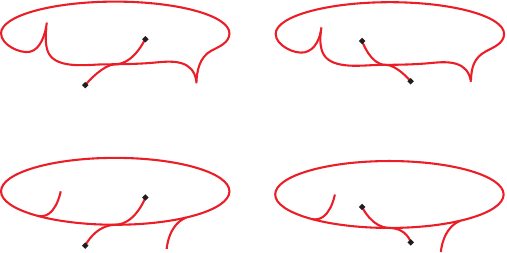}
		\put(20,25){(a)}
		\put(75,25){(b)}
		\put(20,-5){(c)}
		\put(75,-5){(d)}
		\put(29,41){{\footnotesize $X^u$}}
		\put(9,28.5){{\footnotesize $X^d$}}
	\end{overpic}
	\vspace{5pt}
	\caption{Possible configurations of $\tau_D'$ at $\mathcal{A}_l$}
	\label{fig:doubletraintrack}
\end{figure}

Let $\kappa'$ be the subarc of $\kappa$ between the two points of $\kappa\cap\rho_x$.  By collapsing the arc 
$\kappa'$ to a point, we can pinch $\tau_D$ to a new train track $\tau_D'$ that also fully carries $\delta$.  
The train tracks $\tau_D'\cap\mathcal{A}_l$ corresponding to the four possible configurations  
in Figure~\ref{fig:double} are shown in Figure~\ref{fig:doubletraintrack} respectively. Notice that, if  $\mathcal{A}_l$ 
and $\mathcal{A}_r$ are switched, the train tracks $\tau_D'$ in Figure~\ref{fig:double}(a, b) are basically 
the same as  the train track $\tau_D$ for two short paths in Section~\ref{sec:OneArTwoAl}. Also, the train tracks 
$\tau_D'$ for Figure~\ref{fig:double}(c, d) are basically the same as the train track $\tau_D$  for a long path in
 Section~\ref{sec:NoLongPath}. The proofs for the two types of configurations will also be the same as the arguments 
 in Sections~\ref{sec:OneArTwoAl} and \ref{sec:NoLongPath} respectively.

Recall the definition of a splitting arc in Definition~\ref{def:splitting-arc}. The original train track $\tau_D$ 
has two cusps which survive the pinching of $\kappa'$. The splitting arc in $\NN(\tau_D)$ that connects 
the two cusps and which is disjoint from $\delta$ in $\NN(\tau_D)$ was denoted by $s_j$.  Note that 
Lemma~\ref{lem:TwoDeltaArcs} implies that $s_j$ passes through both $\mathcal{R}^d$ and $\mathcal{R}^u$.
Denote the image of $s_j$ in $\NN(\tau_D')$ by  $s_j'$. Clearly, $s_j'$ passes through both rectangles 
$\mathcal{R}^u$ and $\mathcal{R}^d$. For the train track $\tau_D'$, we have an additional splitting arc, 
which we denote by $s_\rho'$, corresponding to the pinching along $\kappa$, i.e.~cutting $\NN(\tau_D')$ 
along $s_\rho'$ will undo the pinching operation along $\kappa$. Since $\kappa$  is an arc in 
$(\alpha\cup\gamma)\cap\widehat{\Sigma}$, this new splitting arc $s_\rho'$ intersects 
$\alpha\cup\gamma$ at least once. Furthermore as $\varepsilon\cap\delta=\emptyset$, 
the intersection $\varepsilon\cap\NN(\tau_D')$ consists of arcs parallel to $s_j'$ and $s_\rho'$.
In all four possible configurations in the Figure~\ref{fig:doubletraintrack}, the surface 
$\Sigma\ssm\NN(\tau_D')$ is an annulus with two cusps at each boundary curve, see 
Figure~\ref{fig:Annulus1}  in  Section~\ref{sec:OneArTwoAl}, for a picture.  Now follow the 
argument in Section~\ref{sec:OneArTwoAl}.  The first step is to perform a maximal number of 
wave moves on $\varepsilon$ along $s$-waves in the annulus $\Sigma\ssm\NN(\tau_D')$, see 
Figure~\ref{fig:Annulus1} for a picture of such $s$-waves. Denote the resulting meridian by $\varepsilon'$.  
As in Section ~\ref{sec:OneArTwoAl}, $\varepsilon'$ does not pass through all the cusps, in other words,
 $\varepsilon'\cap \NN(\tau_D')$ consists of arcs all of which are parallel to either $s_j'$ or $s_\rho'$. 

First note that Claim~\ref{cl:MoreThanOneArc} is also valid for $\NN(\tau_D')$ and $\varepsilon'$, that is 
$\varepsilon'\cap \NN(\tau_D')$ contains at least two arcs. The slight difference  is that, under the setting of
 Sections~\ref{sec:OneArTwoAl}, the proof of Claim~\ref{cl:MoreThanOneArc} needs the fact that the core 
curve $\mathfrak{a}_r$ of the annulus $\mathcal{A}_r$ can be isotoped into $\NN(\tau_D)$, while in the 
current setting, the core curve $\mathfrak{a}_l$ of $\mathcal{A}_l$ can be isotoped into $\NN(\tau_D')$ 
in all four configurations of Figure~\ref{fig:doubletraintrack}. 

Similar to Claim~\ref{cl:Vwave}, since both $s_j'$ and $s_\rho'$ intersect $\alpha\cup\gamma$, this means 
that $\alpha\cup\gamma$ contains both  $[\varepsilon'^-,\varepsilon'^+]$- and $[\delta^-,\delta^+]$-edges. 
So, by Lemma~\ref{lem:+ and -}, there is no wave with respect to $\{\delta,\varepsilon'\}$. Moreover, since $s_j'$ 
passes through both $\mathcal{R}^d$ and $\mathcal{R}^u$, $s_j'$ meets $\gamma$ at least twice and there are 
two $[\varepsilon'^-,\varepsilon'^+]$- or $[\delta^-,\delta^+]$-edges. This means that the Heegaard diagram 
formed by $\{\alpha,\gamma\}$ and $\{\delta,\varepsilon'\}$ is not a standard Heegaard diagram of $S^3$ or 
$(S^2\times S^1)\# L(p,q)$. Hence there must be  a wave with respect to $\{\alpha,\gamma\}$. 
By Lemma~\ref{lem:FirstWave}, this wave must be with respect to $\gamma$.

\begin{figure}[!ht]
	\begin{overpic}[width=10cm]{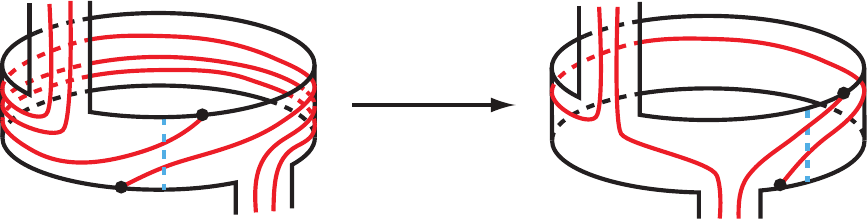}
		\put(43,9){isotopy}
	\end{overpic}
	\caption{An isotopy converting Figure~\ref{fig:double}(c) to Figure~\ref{fig:double}(b) up to symmetry}
	\label{fig:longisotopy2}
\end{figure}

If the configuration of $\rho_x$ is as shown in Figures~\ref{fig:double}(c or d), then the corresponding 
configuration of the train track $\tau_D'$ is Figure~\ref{fig:doubletraintrack}(c or d) respectively. In both 
configurations, the core curve $\frak{a}_l$ of $\mathcal{A}_l$ is carried by the train track $\tau_D'$ and 
hence has an orientation compatible with the orientation of $\delta$.  The argument for these two 
configurations is the same as in Section~\ref{sec:NoLongPath}:  Since there is no wave 
with respect to $\{\delta,\varepsilon'\}$, by part (1) of Proposition~\ref{pro:Paths}, the intersection 
points of $\delta\cap\gamma$ all have the same sign. As in Section~\ref{sec:NoLongPath}, this implies 
that $\frak{a}_l$ and hence the planar surface $P_l$ is disjoint from $\gamma$ and the disk $C$.  
Similar to Section~\ref{sec:NoLongPath}, this means that $P_l$ is an annulus and $\frak{a}_l\cap\alpha$ 
is a single point. Thus  the arc $\kappa$ is the only arc in $(\alpha\cup\gamma)\cap\mathcal{A}_l$.  
Similar to Section~\ref{sec:NoLongPath} and Figure~\ref{fig:longisotopy}, we can perform an isotopy in 
$\mathcal{A}_l$ which fixes $\kappa$ and twists $\mathcal{A}_l$ along the short path, such that the cusp 
directions point into  the rectangles $\mathcal{R}^d$ and $\mathcal{R}^u$ after the isotopy, see 
Figure~\ref{fig:longisotopy2}.  So this isotopy basically converts the configurations of  Figure~\ref{fig:double}(c, d)
to the configurations of Figure~\ref{fig:double}(b, a), up to symmetry. Moreover the isotopy does not affect 
$\alpha\cup\gamma$.   

So it remains to consider the configurations in Figure~\ref{fig:double}(a, b).  
The next step in the proof is to convert the configurations in Figure~\ref{fig:double}(a, b) to 
the setup in Section~\ref{sec:OneArTwoAl} and apply the arguments used in that section. 

\begin{remark}\label{rem: TooShrt}\hfill
 
 \noindent(1) Recall that, for the train track in Section~\ref{sec:OneArTwoAl}, both splitting arcs  $s_j$ 
 and $s_\rho$  pass through both rectangles $\mathcal{R}^d$ and $\mathcal{R}^u$. In the present case 
 $\NN(\tau_D')$ also has  two splitting arcs $s_j'$ and $s_\rho'$ where $s_\rho'$ is the short 
 splitting arc corresponding to the pinching operation. In the current setting for $\NN(\tau_D')$, 
 the splitting arc  $s_j'$ indeed passes through both rectangles $\mathcal{R}^d$ and $\mathcal{R}^u$ but 
 the splitting arc $s_\rho'$  is too short and does \emph{not pass} through the rectangles.
\vskip5pt
 \noindent(2) The reason we need the property that both splitting arcs pass through both rectangles 
 is the following: In the argument in Section~\ref{sec:OneArTwoAl} (after Claim~\ref{cl:MoreThanOneArc}), 
 we perform a sequence of wave moves on $\{\alpha,\gamma\}$ and obtain a new set of meridians 
 $\{\alpha',\gamma'\}$. Each $\alpha'$- or $\gamma'$-arc in $\mathcal{R}^d$ and $\mathcal{R}^u$ intersects the 
 splitting arcs $s_j$ and $s_\rho$ and hence contain $[\varepsilon'^-,\varepsilon'^+]$- and $[\delta^-,\delta^+]$-edges.
Thus we can continue the sequence of wave moves \emph{until $\alpha'$ and $\gamma'$ do not intersect 
$\mathcal{R}^d$ and $\mathcal{R}^u$ at all.} This implies that the intersection points of $\alpha'$ and 
$\gamma'$ with $\partial_+P$ have consistent signs after these wave moves. So, the property that both 
splitting arcs pass through the rectangles plays a crucial role in the final part of the argument in 
Section~\ref{sec:OneArTwoAl}. 
\end{remark}

\begin{figure}[!ht]
	\begin{overpic}[width=10cm]{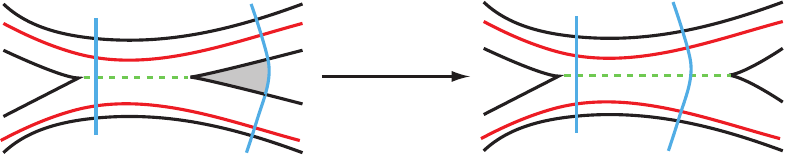}
		\put(43,11){enlarge}
		\put(21,7.3){{\scriptsize$s_\rho'$}}
		\put(89,7){{\scriptsize$s_\rho'$}}
	\end{overpic}
	\caption{Enlarge $\NN(\tau_D')$}
	\label{fig:enlarge}
\end{figure}

To overcome this obstacle we next  ``enlarge'' $\NN(\tau_D')$ and $s_\rho'$ so that the new $s_\rho'$ has 
the desired property. That is, each $\alpha$- or $\gamma$-arc in  $\mathcal{R}^d$ or $\mathcal{R}^u$ gives 
rise to an intersection arc with the new $s_\rho'$. The ``enlarging operation'' is  as follows: 

If there is an $\alpha$- or a $\gamma$-arc that forms a triangle with a cusp of $\NN(\tau_D')$ 
(see the shaded triangle in Figure~\ref{fig:enlarge}), then we  enlarge $\NN(\tau_D')$ to contain 
this arc. The operation is illustrated in Figure~\ref{fig:enlarge}, where the dashed arcs denote 
$s_\rho'$. This operation on $\NN(\tau_D')$ is similar in spirit to the pinching operation on a train track.
 
In order to do that, we will consider the gluing map $\varphi\colon\partial_+P^u\to\partial_+P^d$ 
which glues the two  components of $\partial\widehat{\Sigma}$ together to form $\Sigma$. The 
only thing we know for sure about $\varphi$ is that  $\varphi(X^u)=X^d$. Otherwise it can stretch 
and contract various segments  when gluing the arc $\partial P_+^d \ssm X^d$ to the arc 
$\partial P_+^d \ssm X^d$.  We will consider the various possibilities of the gluing by $\varphi$ 
by studying the configuration of $\delta$ near the two sides of $\partial_+P$.  

\begin{figure}[!ht]
	\begin{overpic}[width=8cm]{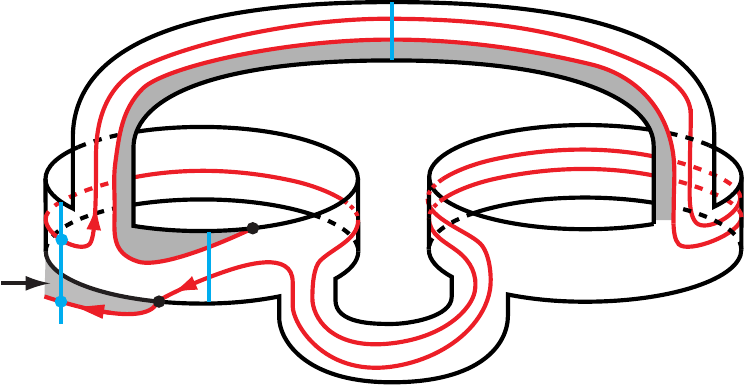}
		\put(17.5,17){{\small $\Delta$}}
		\put(-5,12.5){$\Delta'$}
		\put(33,17.3){{\scriptsize $X^u$}}
		\put(20,7.3){{\scriptsize $X^d$}}
		\put(6,5){$\gamma$}
		\put(27,7.5){$\kappa$}
		\put(51,40){$\gamma$}
	\end{overpic}
	\vspace{5pt}
	\caption{A global picture for Figure~\ref{fig:double}(a)}
	\label{fig:case1b1}
\end{figure}

If the configuration of $\rho_x$ is as in Figure~\ref{fig:double}(a), then consider the shaded triangle region 
$\Delta$ in  Figure~\ref{fig:case1b1} which contains $X^u$ as a vertex.  The gluing map $\varphi$ sends 
a boundary edge of $\Delta$ to a subarc of $\partial_+P^d$.  We can view how things are identified by 
$\varphi$ in the picture of a cyclic cover of $\Sigma$ dual to $\partial_+P$, see Figure~\ref{fig:case1b1}, 
where the shaded region $\Delta'$ is a translation of $\Delta$.  The location of $\Delta'$ depends on 
the gluing map $\varphi$. Since $\gamma$ intersects $\mathcal{R}^u$, there must be a $\gamma$-arc 
meeting the  interior of $\Delta$, see the top blue arc in Figure~\ref{fig:case1b1}.  So there is a $\gamma$-arc 
meeting the  interior of $\Delta'$ as illustrated in Figure~\ref{fig:case1b1}.

Let $R$ be the rectangle in $\mathcal{A}_l\ssm(\mathcal{J}_l^u\cup\mathcal{J}_l^d)$ between the two 
junctions  which  contains the points $X^u$, $X^d$ and the arc $\kappa$. Note that $R$ is invariant 
under the involution $\pi$. Moreover, $\pi|_{R}$ is a $180^\circ$-rotation around the center of $R$, 
and the vertical arc of $R$ that contains its center is invariant under $\pi$.  Let $r^u$ be the boundary 
edge of $R$ in $\partial_+P^u$ (i.e.~the top edge of $R$). $\pi$ rotates $r^u$ to the bottom edge 
of $R$. 

Consider the subarc of $r^u$ that lies in the boundary of $\Delta$, i.e.~the arc $r^u\cap\partial\Delta$, 
and the two $\delta$-arcs in $R$ attached to $X^u$ and $X^d$. If the length of the arc $r^u\cap\partial\Delta$ 
is less than half of the length of $r^u$, then the symmetry from $\pi$ implies that the vertical arc of $R$ which 
contains its center must separate the two $\delta$-arcs attached to $X^u$ and $X^d$. In this case no vertical 
arc of $R$ intersects both $\delta$-arcs, which contradicts the hypothesis that $\kappa$ intersects both 
$\delta$-arcs, see Figure~\ref{fig:case1b1}. Therefore, the length of the arc $r^u\cap\partial\Delta$  
must be at least half of the length of $r^u$ and this implies that part of the arc $\varphi(r^u\cap\partial\Delta)$ 
must be outside $R$. In other words, $\varphi(r^u\cap\partial\Delta)$ must reach the left junction in 
Figure~\ref{fig:case1b1}. By the configuration of $\Delta$ and $\Delta'$, this means that the 
$\gamma$-arc that intersects $\Delta'$ (the left blue arc in  Figure~\ref{fig:case1b1}) must intersect 
the $\delta$-curve in $\widehat{\Sigma}$ further to the left of the left junction of $\mathcal{A}_l$. 

As illustrated by the two blue dots in the $\gamma$-arc that meets $\Delta'$ in Figure~\ref{fig:case1b1},  
two intersection points of $\delta$ with this $\gamma$-arc have opposite signs, contradicting part (1) of
 Proposition~\ref{pro:Paths}.  This means  that the configuration in Figure~\ref{fig:double}(a) cannot happen.

\begin{figure}[!ht]
	\begin{overpic}[width=8cm]{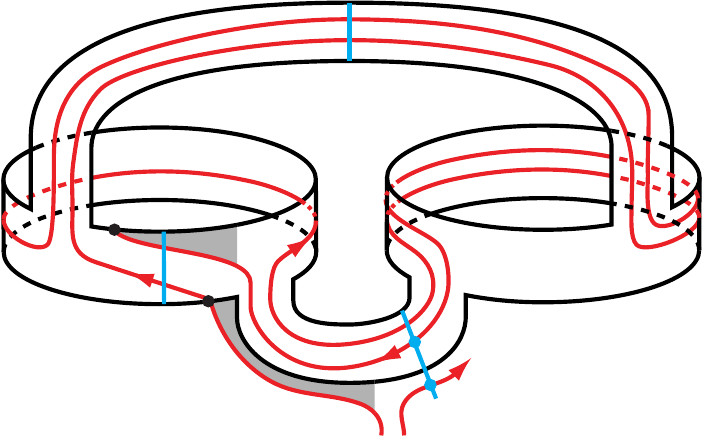}
		\put(30,25.7){{\small $\Delta$}}
		\put(15,30.4){{\tiny $X^u$}}
		\put(28.3,20.5){{\tiny $X^d$}}
		\put(61,2.5){$\alpha$- or $\gamma$-arc}
		\put(22,15){$\kappa$}
		\put(48,50){$\gamma$}
	\end{overpic}
	\vspace{5pt}
	\caption{The first  possible global picture for Figure~\ref{fig:double}(b)}
	\label{fig:case1b2}
\end{figure}

It remains to consider the configuration of $\rho_x$ as in Figure~\ref{fig:double}(b).   
Let $\Delta$ be the triangular region in Figure~\ref{fig:case1b2} and we study the gluing 
map $\varphi$ from  the translation of $\Delta$ in a cyclic cover of $\Sigma$ dual to $\partial_+P$. 
Again there are two possibilities:

\vskip5pt

 The first possibility is that the shaded triangular region marked $\Delta$ in Figure~\ref{fig:case1b2} 
 plus the junction next to $\Delta$ are glued along $\partial_+P^d$ as illustrated in Figure~\ref{fig:case1b2}.
 Figure~\ref{fig:case1b2} is the case where the junction next to $\Delta$ is glued to $\partial\mathcal{A}_l$
 or  $\partial\mathcal{R}^d$ and there is an $\alpha$- or a $\gamma$-arc that intersects $\mathcal{R}^d$ 
 and  meets a $\delta$-arc on the other side of the junction, see the bottom blue arc in Figure~\ref{fig:case1b2}.  
 As illustrated by the two blue dots on this arc in Figure~\ref{fig:case1b2}, this means that two intersection 
 points of $\delta$ with this $\alpha$- or $\gamma$-arc have opposite signs, contradicting part (1) of 
 Proposition~\ref{pro:Paths}.  

\begin{figure}[!ht]
	\begin{overpic}[width=8cm]{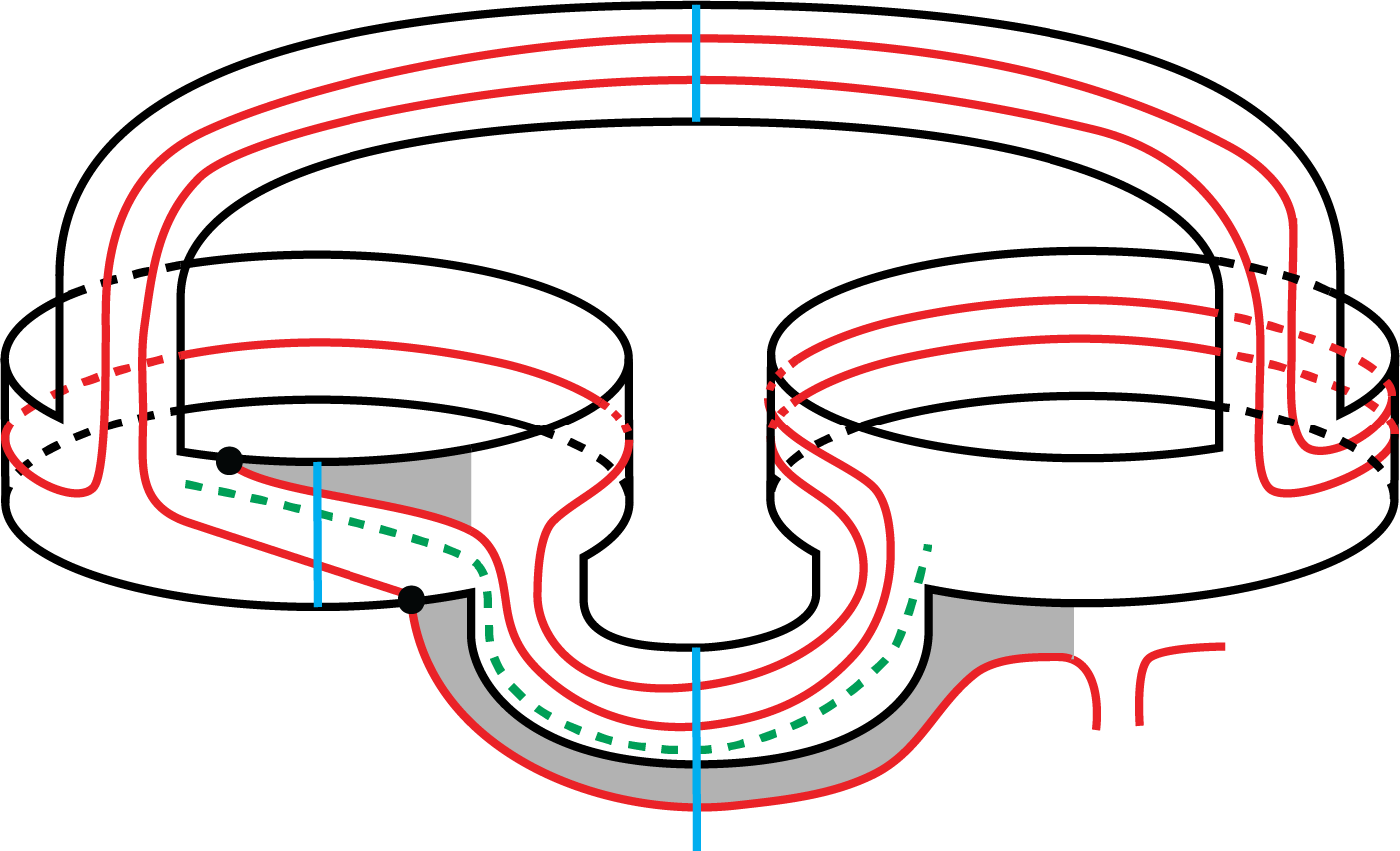}
		\put(30,24){$\Delta$}
		\put(65,21.5){$s_\rho '$}
		\put(21,13.5){$\kappa$}
		\put(48,-3){$\gamma$}
		\put(48,48){$\gamma$}
	\end{overpic}
	\vspace{5pt}
	\caption{The second possible global picture for Figure~\ref{fig:double}(b)}
	\label{fig:case1b3}
\end{figure}

The second possibility is that there is no such $\alpha$- or $\gamma$-arc.  Then after an isotopy 
(fixing $\alpha\cup\gamma$), we may assume that the triangular region $\Delta$ is 
glued to $\partial_+P^d$ and covers an entire boundary edge of $\mathcal{R}^d$, as illustrated 
in Figure~\ref{fig:case1b3}.  

Now consider $\NN(\tau_D')$ and the splitting arc $s_\rho'$ described earlier.  
As illustrated in Figure~\ref{fig:enlarge}, if an $\alpha$- or a $\gamma$-arc  in 
$\Sigma\ssm\NN(\tau_D')$ is a $\partial$-parallel arc that cuts off a neighborhood of a cusp, then 
we can enlarge $\NN(\tau_D')$ to include this neighborhood of the cusp and extend the corresponding 
splitting arc $s_j'$ or $s_\rho'$.  As the $\delta$-edge of $\partial\Delta$ is parallel to $\partial_+P$ 
and since $\Delta$ is glued to cover a whole edge of $\partial\mathcal{R}^d$, we can enlarge 
$\NN(\tau_D')$ to contain $\Delta$, and in particular, contain this boundary edge of $\partial\mathcal{R}^d$. 
Let $\NN(\tau_D')^+$  be the surface obtained by enlarging $\NN(\tau_D')$ as above.   As illustrated in 
Figure~\ref{fig:case1b3}, the splitting arc $s_\rho'$ can be extended to pass through the rectangle 
$\mathcal{R}^d$.  This picture is symmetric under $\pi$, so we can enlarge $\NN(\tau_D')$ and extend 
the splitting arc $s_\rho'$ in the other direction to pass through $\mathcal{R}^u$ as well.   As in 
Section~\ref{sec:OneArTwoAl}, $\Sigma\ssm\NN(\tau_D')^+$ is still an annulus with two cusps at 
each boundary component.

Now the configuration  becomes the same as Section~\ref{sec:OneArTwoAl}.  In particular, both splitting
arcs $s_j'$ and $s_\rho'$ pass through the two rectangles $\mathcal{R}^d$ and $\mathcal{R}^u$.  
Thus we have all the ingredients needed  for the argument in Section~\ref{sec:OneArTwoAl}  which 
can be applied to conclude that  $K$ is doubly primitive. This finishes the proof for the special case.


\section{$\delta$ takes  one short path in each annulus, the general case}\label{sec:OneAndOne}

In this section we deal with the remaining general case of  Configuration 1.
Namely: The curve  $\delta$  takes only two short paths, one in $\mathcal{A}_l$ and one in $\mathcal{A}_r$, 
and the configurations in section~\ref{sec:OneAndOne-S} and Figure~\ref{fig:double} do not occur.

As before, we assume $\delta$ and $\partial_+P$ are invariant under the involution $\pi$ 
and $X=\delta\cap\partial_+P$ is a fixed point of $\pi$.  By symmetry (2) in Remark~\ref{rem:Symmetries}, 
we may assume $X^d\cup X^u$ lies in either $\partial A_l$ or $\partial\mathcal{R}^d\cup\partial\mathcal{R}^u$.

\begin{proposition}\label{pro:OneAndOne} 
Suppose that the curve $\delta$ takes one short path in $\mathcal{A}_l$ and one short path in 
$\mathcal{A}_r$. Then $K$ is doubly primitive.
\end{proposition}

\begin{proof}
Suppose the configuration is not the one in Section~\ref{sec:OneAndOne-S}.
Hence for each component $\kappa$ of $(\alpha\cup\gamma)\cap\widehat{\Sigma}$ so that 
$\kappa\cap \tau_D=\kappa\cap\rho_x$, $\kappa\cap\rho_x$  contains at most one point.

Consider the train track $\tau_D$.  As explained at the beginning of this chapter, the shape of $\tau_D$ 
must be as in Figure~\ref{fig:traintrack}.  We may assume $\tau_D$ is also invariant under $\pi$. 
The train track $\tau_D$ consists of three segments one of which has cusp directions at both endpoints  
pointing into the segment. Call this segment the $x$-arc of $\tau_D$. The  other two segments 
have cusps directions pointing away from the segments and we call them the $y$-arcs, see 
Figure~\ref{fig:traintrack}.  Note that the special segment $\rho_x$ is a $y$-arc and, as before, the 
weight of $\delta$ at $\rho_x$ is one. 

Let $\NN(\tau_D)$  be a small neighborhood of $\tau_D$. It is a once-punctured torus with two cusps 
on its boundary corresponding to the two cusps of $\tau_D$.   We may assume $\varepsilon\cap\NN(\tau_D)$ 
consists of essential arcs in $\NN(\tau_D)$.  As $\varepsilon\cap\delta=\emptyset$ and $\delta$ is a simple 
closed curve in $\NN(\tau_D)$, $\varepsilon\cap\NN(\tau_D)$ is a collection of parallel arcs going into one 
cusp of $\NN(\tau_D)$ and coming out of the other cusp.  It follows from part (2) of Proposition~\ref{pro:Paths} 
that the orientations of these parallel arcs $\varepsilon\cap\NN(\tau_D)$ (induced from an orientation of 
$\varepsilon$) are all the same.  

Let $T$  denote the complement $\Sigma\ssm\NN(\tau_D)$.  So $T$ is also a once-punctured torus 
(with two cusps on the boundary), see Figure~\ref{fig:PuncturedTorus} for a picture of $\varepsilon$ near 
$\partial T$.
 
\begin{figure}[!ht]
\begin{overpic}[width=4cm]{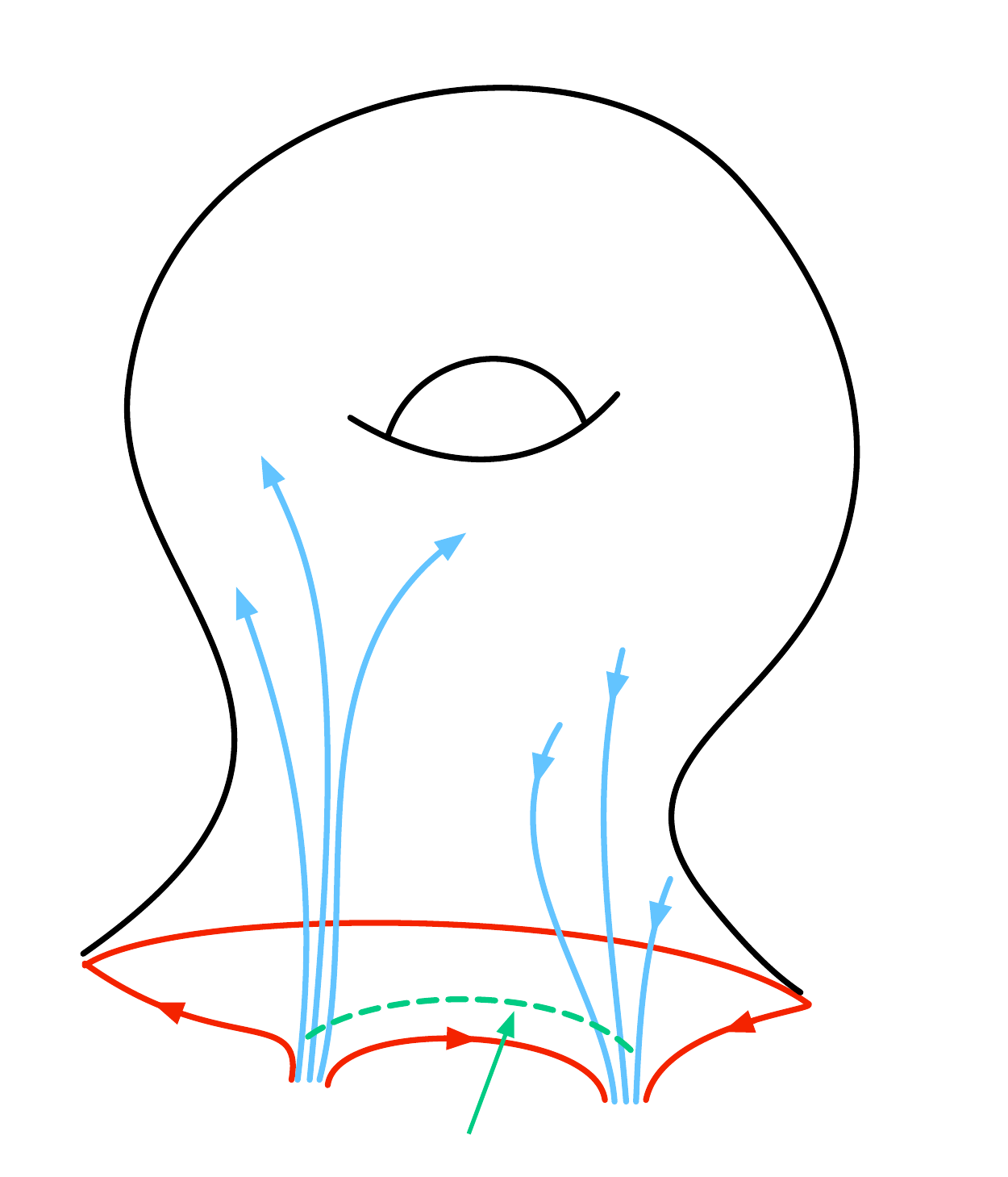}
\put(40,50){$\varepsilon$}
\put(37,-1){$\rho$}
\put(59,4.5){$\delta$}
\end{overpic}
\caption{The once punctured torus $T$ and the $s$-wave $\rho$ with respect to $\varepsilon$.}
\label{fig:PuncturedTorus}
\end{figure}

If $\varepsilon$ contains arcs in the cusps of $\NN(\tau_D)$, then as shown in Figure~\ref{fig:PuncturedTorus}, 
there is an $s$-wave $\rho$ with respect to $\varepsilon$, parallel to one of the $\delta$-arcs connecting the 
two cusps on $\partial T$. If there are  $k_e$ components of $\varepsilon\cap\NN(\tau_D)$, then we 
can perform $k_e$ consecutive wave moves on $\varepsilon$ along such $s$-waves $\rho$. The
result of the wave moves is a meridional curve $\varepsilon'$ of $W$ which is completely contained in $T$. 
 Let $E'$ be the disk in $W$ bounded by $\varepsilon'$.

\vskip5pt

The proof will now proceed by the following sequence of claims:

\vskip7pt

\begin{claim}\label{cl:PD} There is a properly embedded planar surface $P_{\varepsilon'}$ in $U$ 
such that $(P_{\varepsilon'}, E')$ forms a $(\mathcal{P},\mathcal{D})$-pair, and $\partial_-P_{\varepsilon'}$ 
has the same slope as $\partial_-P$ in the boundary torus $\partial_{-}U$.
\end{claim}

\begin{proof}[Proof of Claim~\ref{cl:PD}]
As explained in Remark~\ref{rem:core}, the result of a $\partial$-compression on $P$ is a pair of 
planar surface $P_l$ and $P_r$ with $\partial_+ P_l$ and $\partial_+ P_r$ isotopic to the core 
curves $\frak{a}_l$ and $\frak{a}_r$ of $\mathcal{A}_l$ and $\mathcal{A}_r$ respectively.  

Since $\delta$ takes one short path in both $\mathcal{A}_l$ and $\mathcal{A}_r$,  
after isotopy, each of $\frak{a}_l\cap\NN(\tau_D)$ and $\frak{a}_r\cap\NN(\tau_D)$ is a single nonseparating 
arc in $\NN(\tau_D)$. Note that the two arcs $\frak{a}_l\cap\NN(\tau_D)$ and $\frak{a}_r\cap\NN(\tau_D)$ 
cannot be parallel in $\NN(\tau_D)$ because this would imply that a band sum of $\frak{a}_l=\partial_+ P_l$ 
and $\frak{a}_r=\partial_+ P_r$ produces a curve disjoint from $\NN(\tau_D)$ but parallel to $\partial_+P$. 
This is a contradiction to the fact that $\delta\cap\partial_+P$ is a single point.  So  $\frak{a}_l\cap\NN(\tau_D)$ 
and $\frak{a}_r\cap\NN(\tau_D)$ are non-isotopic arcs in the once-punctured torus $\NN(\tau_D)$.  Hence 
the endpoints of  $\frak{a}_l\cap\NN(\tau_D)$ and $\frak{a}_r\cap\NN(\tau_D)$ alternate along the boundary 
of the once-punctured torus (one can also see this by checking the possible configurations geometrically).  
Now, consider $ T = \Sigma \ssm \NN(\tau_D)$. The endpoints of $\frak{a}_l\cap T$ and  $\frak{a}_r\cap T$ 
alternate along $\partial T$, which implies that  $\frak{a}_l\cap T$ and  $\frak{a}_r\cap T$ 
are also non-isotopic essential arcs in the once-punctured torus $T$.

Use the arcs $\frak{a}_l\cap T$ and  $\frak{a}_r\cap T$ as representatives of a basis 
for  $H_1(T,\partial T)\cong\mathbb{Z}\oplus\mathbb{Z}$, and suppose they represent elements 
with  slope $1/0$ and $0/1$ respectively.  Since $\varepsilon'\subset T$, there is a properly
embedded arc in $T$ which intersects $\varepsilon'$ in a single point.  Suppose this arc has slope $p/q$,
where $p$ and $q$ are coprime. Take $p$ parallel copies of  $P_l$ and $q$ parallel copies of $P_r$, 
and perform a sequence of band  sums of these planar surfaces along  $\partial T$ to obtain a planar
surface  $P_{\varepsilon'}$ so that the intersection of $T$ with the resulting curve $\partial_+P_{\varepsilon'}$
is an arc of slope $p/q$. Hence the planar surface $P_{\varepsilon'}$ is such that  
$\varepsilon'\cap\partial_+P_{\varepsilon'}$ is a single point. Hence $(P_{\varepsilon'}, E')$ 
forms a $(\mathcal{P},\mathcal{D})$-pair.
\end{proof}

Isotope $\varepsilon'$ to intersect $\alpha$ and $\gamma$ minimally. 
The following is an immediate corollary of Claim~\ref{cl:PD}.

\begin{corollary}\label{cor:noepsilonwave}
No subarc of $\varepsilon'$ is a wave with respect to $\{\gamma,\alpha\}$.
\end{corollary}
\begin{proof}[Proof of Corollary~\ref{cor:noepsilonwave}]
The corollary is basically the same as Corollary~\ref{cor:NoWaves}. 

By Claim~\ref{cl:PD}, $(P_{\varepsilon'}, E')$ forms a $(\mathcal{P},\mathcal{D})$-pair. 
If  $P_{\varepsilon'}$ is an annulus, then by Proposition~\ref{pro:DPEquivalentToAD}, $K$ 
is doubly primitive and we are done. So we may assume that $P_{\varepsilon'}$ is not an annulus. 
After replacing the planar surface $P$ in Lemma~\ref{lem:NoAWaves} with 
$P_{\varepsilon'}$, one can conclude that any $s$-wave with respect to $\{\gamma,\alpha\}$ must intersect 
$\partial_+ P_{\varepsilon'}$: As in the proof of Corollary~\ref{cor:NoWaves}, if a subarc $\eta$ of $\varepsilon'$ is a 
wave with respect to $\{\gamma,\alpha\}$, then $\eta\cap\partial_+ P_{\varepsilon'}\ne\emptyset$. By the symmetry 
induced by  $\pi$, the curve $\varepsilon'$ has another subarc $\eta'=\pi(\eta)$ which is also a wave with respect to
$\{\gamma,\alpha\}$ and $\eta'\cap\partial_+ P_{\varepsilon'}\ne\emptyset$. This contradicts the fact that 
$\varepsilon'\cap \partial_+ P_{\varepsilon'}$ is a single point as $(P_{\varepsilon'}, E')$ is a  
$(\mathcal{P},\mathcal{D})$-pair.
\end{proof}

As before, consider the Heegaard diagram formed by $\widehat{W}=\{\delta,\varepsilon'\}$ and 
$\widehat{V} = \{\gamma,\alpha\}$.  Let $s_j$ be the splitting arc in $\NN(\tau_D)$ which connects the two 
cusps,  so that $s_j\cap\delta=\emptyset$ and also $\NN(\tau_D)\ssm s_j$ is a product neighborhood of $\delta$.   

Lemma~\ref{lem:TwoDeltaArcs} 
implies that $s_j$ passes through both $\mathcal{R}^d$ and $\mathcal{R}^u$. Since $\gamma$ 
intersects both $\mathcal{R}^d$ and $\mathcal{R}^u$ and since  $\varepsilon'\subset\Sigma\ssm\NN(\tau_D)$, 
our assumption on the  orientation of $\delta$ implies that $\gamma$ has two distinct subarcs in 
$\mathcal{R}^d$ and $\mathcal{R}^u$ that are $[\delta^-,\delta^+]$ blocking-edges.

In a once-punctured torus if two essential embedded curves or arcs intersect minimally then the 
algebraic intersection number  is equal to the geometric intersection number. Hence the 
intersection points of $\varepsilon'$ with each arc of $(\alpha\cup\gamma) \cap T$ all have the 
same sign. 

For each arc of $(\alpha\cup\gamma) \cap T$,  by collapsing all its intersection 
points with $\varepsilon'$ into one point along this arc, we can construct a  train track fully 
carrying $\varepsilon'$.  In fact, we can extend such collapsing/pinching  as much as possible 
 to construct a train track $\tau_E$ which is either a train track as in  Figure~\ref{fig:traintrack} 
 or a circle (in which case each arc of $(\alpha\cup\gamma) \cap T$ intersects $\varepsilon'$ in
 at most one point).  In either case, each arc of $(\alpha\cup\gamma) \cap T$ intersects 
 the train track $\tau_E$ in at most one point. 

\begin{claim}\label{cl:tauE}
 $|(\alpha\cup\gamma)\cap\tau_E|\le c_0(P,D,\alpha,\gamma)$.
\end{claim}

\begin{proof}[Proof of Claim~\ref{cl:tauE}]
Let $\zeta$ be a component of $(\alpha\cup\gamma)\ssm\partial_+P$ and view $\zeta$ as an 
arc in $(\alpha\cup\gamma)\cap\widehat{\Sigma}$. We first consider the possibility that 
$\zeta\cap\tau_D$ contains more then one point. In the construction of $\tau_D$, all the parallel
$\delta$-arcs, except for the special arc $\rho_x$, were pinched into a single segment. 
Thus, $\zeta$ intersects $\tau_D\ssm\rho_x$ in at most one point.  As $\zeta\cap\tau_D$ 
contains more then one point, $\zeta$ must intersect $\rho_x$. 

It is assumed at the beginning of Proposition~\ref{pro:OneAndOne} that the configuration
 is not in the special case in Section~\ref{sec:OneAndOne-S}.  
That is, for each component $\zeta$ of $(\alpha\cup\gamma)\ssm\partial_+P$ so that 
$\zeta\cap \tau_D=\zeta\cap\rho_x$ then $\zeta\cap\rho_x$  contains at most one point.   
Thus, if $\zeta\cap(\tau_D\ssm\rho_x)=\emptyset$, then 
$\zeta\cap\tau_D=\zeta\cap\rho_x$ contains at most one point. Since $\zeta\cap\tau_D$ 
contains more then one point, $\zeta\cap(\tau_D\ssm\rho_x)\ne\emptyset$. Figure~\ref{fig:arcs} 
is a local picture when this happens. Moreover, a subarc of $\zeta$ between two points of 
$\zeta\cap\tau_D$ must connect $\rho_x$ to $\tau_D\ssm\rho_x$.  

\begin{figure}[!ht]
	\vspace{5pt}
	\begin{overpic}[width=3.3cm]{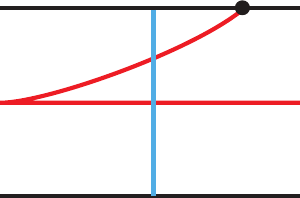}
		\put(65,47){$\rho_x$}
		\put(20,22){$\tau_D$}
		\put(43,14){$\zeta$}
		\put(-23,-2){$\partial_+P$}
		\put(-23,60){$\partial_+P$}
		\put(77.5,66.5){$X$}
	\end{overpic}
	\caption{Intersection of $\zeta$ with $\tau_D$}
	\label{fig:arcs}
\end{figure}

Recall that in our case the train track $\tau_D$ is as in Figure~\ref{fig:traintrack}. Hence, such a 
subarc of $\zeta$ between two points of $\zeta\cap\tau_D$ is an arc in $\Sigma\ssm\tau_D$ 
that ``cuts off" a cusp, as illustrated in Figure~\ref{fig:arcs}. In particular, such an arc corresponds 
to a $\partial$-parallel arc in $\Sigma\ssm\NN(\tau_D)$ which can be assumed to be disjoint from  
$\varepsilon'$ after isotopy. 

Thus, for each component $\zeta$ of $(\alpha\cup\gamma)\ssm\partial_+P$ that meets $\delta$, 
we can find a subarc $\zeta'  \subset \zeta$ so that $\zeta'\cap\delta = \zeta\cap\delta$, i.e. $\zeta'$
contains all the points of $\zeta\cap\delta$, and so that $\zeta'$ is disjoint from $\varepsilon'$. 
Let $\mathcal{C}'$ be the collection of all such subarcs $\zeta'$.  Hence, by definition
 $$|\mathcal{C}'|=c_0(P,D,\alpha,\gamma).$$ 
Since each component of $(\alpha\cup\gamma)\ssm\NN(\tau_D)$ intersects $\tau_E$ in at 
most one point, we have $$|(\alpha\cup\gamma)\cap\tau_E|\le |(\alpha\cup\gamma)\ssm\mathcal{C}'|.$$
As  $|(\alpha\cup\gamma)\ssm\mathcal{C}'|=|\mathcal{C}'|=c_0(P,D,\alpha,\gamma)$, 
we have 
$|(\alpha\cup\gamma)\cap\tau_E|\le c_0(P,D,\alpha,\gamma)$ 
as stated. 
\end{proof}

\begin{claim}\label{cl:NoCircle}
The train track $\tau_E$ cannot be a circle.
\end{claim}

\begin{proof}[Proof of Claim \ref{cl:NoCircle}] 
Suppose to the contrary that $\tau_E$ is a circle. So $\tau_E=\varepsilon'$ and each arc of 
$(\alpha\cup\gamma) \cap T$ intersects $\varepsilon'$ in at most one point.  By Claim~\ref{cl:tauE}, 
we have  $$|\,(\alpha\cup\gamma) \cap\varepsilon'\,|\le c_0(P,D,\alpha,\gamma). $$  

Let $P_{\varepsilon'}$ be as in Claim~\ref{cl:PD} and let $E'$ be the disk in $W$ bounded by $\varepsilon'$.  
So $(P_{\varepsilon'},E')$ is a $(\mathcal{P},\mathcal{D})$-pair. By the definition of complexity, 
$$c_0(P_{\varepsilon'}, E',\alpha,\gamma)\le |\,(\alpha\cup\gamma)\cap \varepsilon'\,|.$$ 
Hence $$c_0(P_{\varepsilon'}, E',\alpha,\gamma)\le  c_0(P,D,\alpha,\gamma).$$ 
Moreover, Lemma~\ref{lem:TwoDeltaArcs} implies that 
$$|\,(\alpha\cup\gamma)\cap \tau_D\,| < |\,(\alpha\cup\gamma)\cap\delta\,|$$ 
By the construction of $\tau_E$, we have  
$$|\,(\alpha\cup\gamma)\cap\tau_E\,|\le |\,(\alpha\cup\gamma)\cap\tau_D\,|.$$ 
As $\tau_E=\varepsilon'$, we have  
$$|\,(\alpha\cup\gamma)\cap\varepsilon'\,|=|\,(\alpha\cup\gamma)\cap\tau_E\,|\le 
|\,(\alpha\cup\gamma)\cap\tau_D\,| < |\,(\alpha\cup\gamma)\cap\delta\,|.$$ 
Thus 
$ c(P_{\varepsilon'}, E',\alpha,\gamma) < c(P,D,\alpha,\gamma)$ 
 and this contradicts Assumption~\ref{ass:minimality}.
\end{proof}

Claim~\ref{cl:NoCircle} means that there is an arc of $(\alpha\cup\gamma) \cap T$ that intersects $\varepsilon'$ 
in more than one point in $T$. Hence a subarc of this arc is an $[\varepsilon'^-,\varepsilon'^+]$ blocking-edge.  
We have concluded before Claim~\ref{cl:tauE} that   $\gamma$ contains two subarcs in $\mathcal{R}^d$ and 
$\mathcal{R}^u$ that are $[\delta^-,\delta^+]$ blocking-edges. So by part (1) of Lemma~\ref{lem:+ and -}, the 
Heegaard diagram has no wave with respect to $\{\delta,\varepsilon'\}$.  
Moreover, since $\gamma$ contains two distinct  $[\delta^-,\delta^+]$-edges, the Heegaard diagram is 
not a standard Heegaard diagram of $S^3$ or $(S^2\times S^1)\# L(r,s)$. Thus, by Theorems~\ref{thm:HOT} 
and \ref{thm:NeOk}, there must be a wave with respect to $\{\alpha,\gamma\}$.  Since by Lemma~\ref{lem:FirstWave}, 
there is no wave with respect to $\alpha$, there must be a wave with respect to $\gamma$.  

By part (1) of Proposition~\ref{pro:Paths}, we may assign orientations to $\alpha$ and $\gamma$ 
so that the intersection points of $\delta$ with $\alpha\cup\gamma$ all have the same sign. 
By Corollary~\ref{cor:noepsilonwave}, the curves $\alpha$, $\gamma$, $\delta$, and $\varepsilon'$ 
satisfy the conditions in Lemma~\ref{lem:Compatibility}. Thus, by Lemma~\ref{lem:Compatibility}, the 
Heegaard diagram must be as depicted in Figure~\ref{fig:octagon}(a) and the $\gamma$-wave $\eta$ must 
connect a $[\delta^-,\delta^+]$-edge to an $[\varepsilon'^-,\varepsilon'^+]$-edge as shown in 
Figure~\ref{fig:octagon}(b). 

As illustrated in Figure~\ref{fig:octagon}(a) all the $[\delta^-,\delta^+]$-edges are parallel. Hence we 
may assume that all the $[\delta^-,\delta^+]$-edges are in the cusp of $\NN(\tau_D)$. Similarly, we may 
assume that all the $[\varepsilon'^-,\varepsilon'^+]$-edges are in the cusp of $\NN(\tau_E)$.
Figure~\ref{fig:2cusp}(a) is a picture of the octagon in Figure~\ref{fig:octagon}(b) that contains  
 the wave $\eta$ together with a picture of the cusps of $\NN(\tau_D)$ and 
$\NN(\tau_E)$. Since the intersection points of $\delta$ with $\alpha\cup\gamma$ all have the same sign, 
the two edges of the octagon in Figure~\ref{fig:2cusp}(a) next to the cusp must have opposite orientation, 
see the orientation of the two arcs marked $\alpha$ and $\gamma$ in Figure~\ref{fig:2cusp}(a).
This means that these two arcs must belong to different curves of $\{\alpha,\gamma\}$, that is, 
one $\alpha$-arc and one $\gamma$-arc as labelled in Figure~\ref{fig:2cusp}(a), since, if they belong 
to the same curve, then an arc in the octagon connecting these two arcs is a wave and this 
wave is as shown in Figure~\ref{fig:octagon}(c). This  contradicts part (2) of Lemma~\ref{lem:Compatibility}
which claims that the wave should be as in  Figure~\ref{fig:octagon}(b) and the two are not compatible.

\begin{figure}[!ht]
	\begin{overpic}[width=5.5cm]{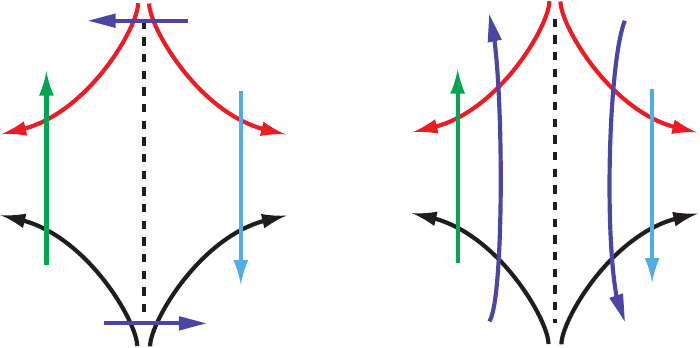}
		\put(16,-8){(a)}
		\put(75,-8){(b)}
		\put(27,38.5){$\delta$}
		\put(27,7){$\varepsilon'$}
		\put(1,23){$\alpha$}
		\put(35,23){$\gamma$}
		\put(21,23){$\eta$}
		\put(60,23){$\alpha$}
		\put(93.5,23){$\gamma$}
		\put(82,45){$\delta$}
		\put(99.5,15){$\varepsilon'$}
	\end{overpic}
	\vspace{3pt}
	\caption{Wave move along $\eta$}
	\label{fig:2cusp}
\end{figure}

\begin{remark}\label{rem:twosides}
Let $\eta$ be the wave as in Figure~\ref{fig:2cusp}(a) (also see Figure~\ref{fig:octagon}(b). By 
Lemma~\ref{lem:FirstWave}, the wave $\eta$ must be with respect to $\gamma$.  As in 
Definition ~\ref{def:wavemove} and similar to the argument in Section~\ref{sec:OneArTwoAl}, 
a wave move along $\eta$ can be done in two steps. The first step is to perform a surgery: connect 
the endpoints of  $\gamma\ssm\NN(\partial\eta)$ using two parallel copies of $\eta$.  It can be 
seen from Figure~\ref{fig:2cusp}(b) that the resulting two 
curves do not form any bigon with $\delta$ and $\varepsilon'$.  One of the resulting curves, which we denote 
by $\alpha_0$, is parallel to $\alpha$, and the other resulting curve, which we denote by $\gamma'$, is 
the final curve in the wave move. The configuration of Figure~\ref{fig:2cusp}(b) implies that the intersection 
of $\alpha_0$ (resp. $\gamma'$) with this octagon is an arc parallel to the edge marked $\alpha$ 
(resp. $\gamma$) in Figure~\ref{fig:2cusp}(b).  The curve $\alpha_0$ is removed in the second step 
of the wave move.  Thus after the wave move, the two subarcs of $\alpha\cup\gamma'$ next to a 
cusp of $\tau_D$ belong to different curves in $\{\alpha,\gamma'\}$, similar to the original 
configuration of $\{\alpha,\gamma\}$ before this wave move.  Thus, as long as $\{\alpha,\gamma'\}$ 
contains both $[\delta^-,\delta^+]$ and $[\varepsilon'^-,\varepsilon'^+]$ edges, the component of 
$\Sigma\ssm(\delta\cup\varepsilon'\cup\alpha\cup\gamma')$ that contains the previous wave remains an octagon, 
see Figure~\ref{fig:octagon}(b) and Figure~\ref{fig:2cusp}(a). As in Lemma~\ref{lem:Compatibility}, the next wave 
with respect to $\{\alpha,\gamma'\}$ must also connect a $[\delta^-,\delta^+]$-edge to an 
$[\varepsilon'^-,\varepsilon'^+]$-edge as depicted in Figure~\ref{fig:2cusp}(a) as well.
\end{remark}

Remark~\ref{rem:twosides} implies that there is a sequence of wave moves $\eta_1,\dots,\eta_k$ along 
waves as shown in Figure~\ref{fig:2cusp}(a) such that after the last wave move $\eta_k$ either 
\vskip3pt
\noindent (a) $\Gamma(\delta,\varepsilon')$ contains no more $[\varepsilon'^-,\varepsilon'^+]$ edges, or 
\vskip3pt
\noindent (b) $\Gamma(\delta,\varepsilon')$ contains no more $[\delta^-,\delta^+]$ edges.
\vskip5pt

 Note that before the last wave move $\eta_k$, the Whitehead graph $\Gamma(\delta,\varepsilon')$ 
contains \emph{both}  $[\delta^-,\delta^+]$  and  $[\varepsilon'^-,\varepsilon'^+]$ edges. 
Hence,  we are left with the following three cases regarding $\tau_E$ and the Whitehead graph 
$\Gamma(\delta,\varepsilon')$:

\vskip7pt

\noindent \underline{\bf Case\,1}. The weight of $\varepsilon'$ at each segment of the train track $\tau_E$ is at least two, and 
$\Gamma(\delta,\varepsilon')$ contains no more $[\varepsilon'^-,\varepsilon'^+]$-edges after the last wave move $\eta_k$.
\vskip5pt
\noindent \underline{\bf Case\,2}.    $\Gamma(\delta,\varepsilon')$ contains no more $[\delta^-,\delta^+]$-edges after the 
last wave move $\eta_k$.
\vskip5pt
\noindent \underline{\bf Case\,3}.   The weight of $\varepsilon'$ at some segment of $\tau_E$ is one, and 
$\Gamma(\delta,\varepsilon')$ contains no more $[\varepsilon'^-,\varepsilon'^+]$-edges after the last wave move $\eta_k$.
\vskip5pt

To simplify notation, we use $\{\alpha',\gamma'\}$ to denote the set of meridional curves of $V$ obtained by the
sequence of wave moves, where $\alpha'$ and $\gamma'$ 
are derived from the original meridians $\alpha$ and $\gamma$ respectively.   
\vskip5pt

\noindent\underline{{\bf Case 1}}:  The weight of $\varepsilon'$ at each segment of the train track $\tau_E$ is at least 
two, and $\Gamma(\delta,\varepsilon')$ contains no more $[\varepsilon'^-,\varepsilon'^+]$-edges after the last wave 
move $\eta_k$.

The goal in this case is to show that either the knot $K$ is a \emph{Berge-Gabai knot} (i.e.~$K$ lies 
in a Heegaard solid torus and the solid torus remains a solid torus after a nontrivial  Dehn surgery on $K$) 
and hence it is doubly primitive by \cite{Be2} and \cite{Gabai}, or  that before the last wave move along $\eta_k$, 
subarcs of $\varepsilon'$ are $[\alpha'^+,\alpha^-]$- edges and $[\gamma'^+,\gamma'^-]$-edges. Since 
$\Gamma(\delta,\varepsilon')$ contains both  $[\delta^-,\delta^+]$  and $[\varepsilon'^-,\varepsilon'^+]$ edges 
before the last wave move $\eta_k$, the latter possibility is a contradiction to Theorems~\ref{thm:HOT} 
and \ref{thm:NeOk}.  So in Case 1 we focus on how $\varepsilon'$ intersects $\alpha'\cup\gamma'$ 
during the sequence of wave moves. This is done by studying the intersection with the train track $\tau_E$.

\vskip5pt

No subarc of $\varepsilon'$ is a wave by Corollary~\ref{cor:noepsilonwave}. 
As in the proof of Lemma~\ref{lem:Compatibility} and illustrated in Figure~\ref{fig:octagon}, we may assign  
orientations to $\varepsilon'$, $\alpha$
and $\gamma$ so that all the intersection points have the same sign. Thus the orientation of $\alpha$ 
and $\gamma$ are compatible along the train track $\tau_E$.  By Claim~\ref{cl:NoCircle}, $\tau_E$ 
is a train track with two cusps as shown in Figure~\ref{fig:traintrack}. Recall that the train track $\tau_E$ is  
constructed so that  each arc of $(\alpha\cup\gamma)\cap T$ intersects $\tau_E$ in at most one point.  

Consider the curves $\alpha'$ and $\gamma'$ during the sequence of wave moves $\eta_1,\dots,\eta_k$ as 
described above. The curves $\alpha'$ and $\gamma'$ have an induced orientation from $\alpha$ and $\gamma$. 
We will inductively require that, at any stage of the sequence of wave moves, $\alpha'$ and $\gamma'$ have the 
following property: 
\vskip5pt
Each arc of $(\alpha'\cup\gamma')\cap T$ intersects the train track $\tau_E$ in at most one point and the 
orientations of $\alpha'$ and $\gamma'$ are compatible along $\tau_E$. 
\vskip5pt
Since the involution $\pi$ leaves $\varepsilon'$, $\alpha'$ and $\gamma'$ invariant up to isotopy, we 
may assume that $\pi$ leaves the train track $\tau_E$ invariant. Hence $\pi$ interchanges the two cusps of $\tau_E$. 
The two cusp points of the train track $\tau_E$ divide $\tau_E$ into three segments: one $x$-arc and two $y$-arcs, 
as in Figure~\ref{fig:traintrack}. 

For any segment $\kappa$ of the train track $\tau_E$, the assumption on $\pi$ implies that $\pi$ leaves $\kappa$ 
invariant and interchanges the two endpoints of $\kappa$.  Since $\pi$ also leaves $\alpha'$ and 
$\gamma'$ invariant up to isotopy, this implies that the intersection points of $(\alpha'\cup\gamma')\cap\kappa$ 
that are outermost in $\kappa$ are interchanged by $\pi$ and hence must belong to the same curve 
$\alpha'$ or $\gamma'$. 

We first show that $\alpha'\cup\gamma'$ must meet at least two segments of the train track $\tau_E$. Assume 
to the contrary that  $\alpha'\cup\gamma'$ intersects only one segment $\kappa$ of $\tau_E$. Now 
consider the region of $\Sigma\ssm(\NN(\tau_D)\cup\NN(\tau_E))$ corresponding to the octagon 
of Figure~\ref{fig:2cusp}(a). We see that the two subarcs of $\alpha'\cup\gamma'$  next to 
a cusp of $\tau_E$ (see the arcs labeled $\alpha$ and $\gamma$ in Figure~\ref{fig:2cusp}(a)) must 
correspond to arcs that intersect $\kappa$ with outermost intersection points along $\kappa$ 
and thus belong to the same curve in $\{\alpha',\gamma'\}$, this implies that 
the two vertical edges of the octagon in Figure~\ref{fig:2cusp}(a) closest to the  cusps of $\tau_D$ and 
$\tau_E$ must belong to the same curve of $\{\alpha',\gamma'\}$. This contradicts the conclusion in 
Remark~\ref{rem:twosides}. Thus we conclude that $\alpha'\cup\gamma'$ must intersect at least 
two segments of train track $\tau_E$. 

\vskip5pt
Let $p:\NN(\tau_E)\to\tau_E$ be the map collapsing $\NN(\tau_E)$ to the train track $\tau_E$. 
Since the cusp direction of $\tau_E$ points into the $x$-arc, there are two parallel and adjacent 
$\varepsilon'$-arcs $t_1$ and $t_2$ in $\NN(\tau_E)$ that collapse onto the $x$-arc by the map $p$.   
If the $x$-arc intersects $\alpha'\cup\gamma'$, then each intersection point corresponds 
to an arc of $\alpha'\cup\gamma'$ between the two $\varepsilon'$-arcs $t_1$ and $t_2$. As the orientation 
is compatible, this arc is an $[\varepsilon'^-,\varepsilon'^+]$-edge.

\vskip5pt

Perform the sequence of wave moves  $\eta_1,\dots,\eta_k$, with respect to $\{\alpha',\gamma'\}$, 
as described above in Remark \ref{rem:twosides}.   There are two situations depending on whether or 
not $\alpha'$ and $\gamma'$ meet the $x$-arc of the train track.

If the $x$-arc of the train track $\tau_E$ intersects $\alpha'\cup\gamma'$, then each intersection point correspond to an 
$[\varepsilon'^-,\varepsilon'^+]$-edge.  As indicated in Figure~\ref{fig:2cusp}(a, b) and Remark~\ref{rem:twosides}, 
a wave move basically ``pushes'' an arc out of a cusp of $\tau_E$. Hence after each wave move, the number 
of intersection points of the $x$-arc with $\alpha'\cup\gamma'$ decreases by at least one.  Moreover, since 
originally each arc of $(\alpha\cup\gamma)\cap T$ intersects $\tau_E$ in at most one point, the configuration of 
Figure~\ref{fig:2cusp}(a) implies that after each wave move, each arc of $(\alpha'\cup\gamma')\cap T$ 
intersects the traintrack $\tau_E$ in at most one point as well.  So we can perform  wave moves until the  $x$-arc 
no longer intersects $\alpha'\cup\gamma'$. Note that  Remark~\ref{rem:twosides} says that the 
two arcs next to a cusp of $\tau_E$ are still one $\alpha'$-arc and one $\gamma'$-arc.  

\begin{figure}[!ht]
	\begin{overpic}[width=7.5cm]{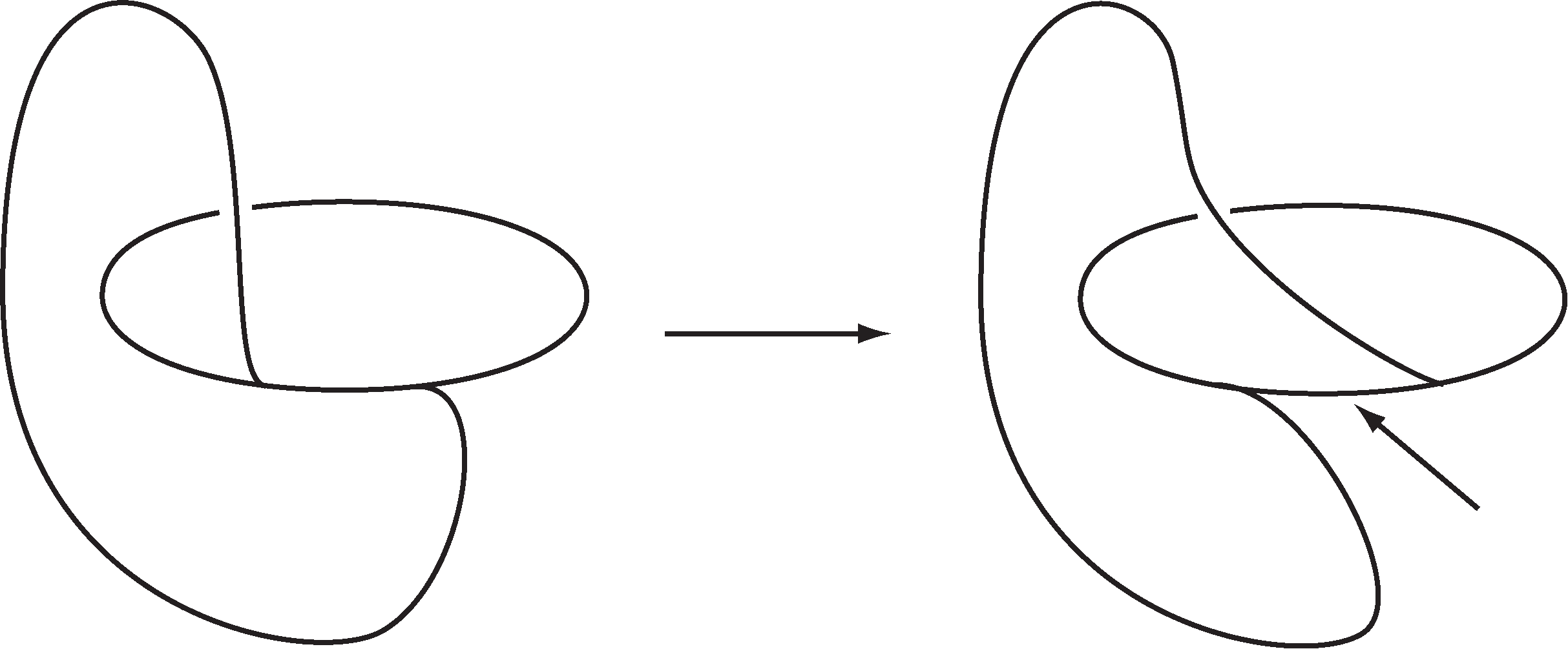}
		\put(41,22){splitting}
		\put(89,6){new $y$-arc}
	\end{overpic} 
	\caption{Splitting the train track $\tau_E$}
	\label{fig:splitting}
\end{figure}

If the $x$-arc of the train track $\tau_E$ does not intersect $\{\alpha'\cup\gamma'\}$ and we have not exhausted  
the sequence of wave moves, then split the train track $\tau_E$ as shown in Figure~\ref{fig:splitting}. 
Note that the splitting does not affect $\alpha'$ and $\gamma'$ as they do not intersect the $x$-arc. 
The splitting changes a $y$-arc of $\tau_E$ into the $x$-arc of the resulting train track, and  creates 
a new $y$-arc that does not intersect $\alpha'\cup\gamma'$.   

Call a $y$-arc of the train track a 
{\it pure} $y$-arc if it either does not intersect $\alpha'\cup\gamma'$ or intersects only one of the 
curves in $\alpha'$ or $\gamma'$. So the new, after this splitting, $y$-arc is a pure $y$-arc.   
We showed earlier that $\alpha'\cup\gamma'$ must intersect at least two segments of train track. 
Since the newly created $y$-arc does not intersect $\alpha'\cup\gamma'$, the new $x$-arc of 
the train track must intersect $\alpha'\cup\gamma'$. 

To simplify notation, we still use $\tau_E$ to denote the train track after the splitting.
Now we continue with the sequence of wave moves.

\begin{claim*}
Each pure $y$-arc of the train track $\tau_E$ remains a pure $y$-arc after the next wave move.
\end{claim*}

\begin{proof}[Proof of the Claim] There are two cases: 
The first is when  $\kappa\subset \tau_E$ is a pure $y$-arc which meets one curve in $\{\alpha',\gamma'\}$, 
say $\alpha'$, before a wave move.  Then as in the argument in Remark~\ref{rem:twosides}, the surgery-step 
of the wave move ``pushes'' an arc out of the cusp creating a pair of arcs parallel to the $\alpha'$- and 
$\gamma'$-edges inside the octagon, see Figure~\ref{fig:2cusp}(b).  

If the wave move is with respect to $\alpha'$, then $\gamma'$ is unchanged by the wave move. 
Since $\kappa$ is a pure $y$-arc and does not intersect $\gamma'$ before the wave move, it 
does not intersect $\gamma'$ after the wave move either. Hence $\kappa$ remains a pure 
$y$-arc after the wave move.

If the wave move is with respect to $\gamma'$, then the surgery-step of the wave move ``splits'' 
$\gamma'$ into two curves $\gamma_1$ and $\gamma_2$. One of $\gamma_1$ and $\gamma_2$ 
is parallel to $\alpha'$.  When viewed in the octagon, the surgery-step creates a pair of arcs parallel 
and next to the $\alpha'$- and $\gamma'$-edges of the octagon  as seen in Figure~\ref{fig:2cusp}(b). 
 Without loss of generality, suppose the arc parallel to the $\gamma'$-edge (resp.~$\alpha'$-edge) in the octagon 
 belongs to $\gamma_1$ (resp.~$\gamma_2$). Since the pure $y$-arc $\kappa$ only meets 
 $\alpha'$, this means that $\gamma_2$ intersects $\kappa$ while $\gamma_1$ does not. 
Recall that in the construction of a wave move,  the set of curves $\gamma'\cup\gamma_1\cup\gamma_2$ 
bounds a pair of pants and one of $\gamma_1$ or $\gamma_2$  is parallel and next to $\alpha'$. Since the 
arc parallel to the $\gamma'$-edge in the octagon belongs to $\gamma_1$, $\gamma_1$ cannot 
be the curve parallel to $\alpha'$.  Hence $\gamma_2$ is parallel to $\alpha'$ and is deleted in the second 
step of the wave move and $\gamma_1$ is the new $\gamma'$-curve resulting from the wave move.  
Since $\gamma_1$ does not intersect $\kappa$ and $\gamma_2$ is removed, this means that after 
this wave move the pure $y$-arc $\kappa$ only meets $\alpha'$ (after the wave move) and hence 
remains a pure $y$-arc. 

The second case is when a $y$-arc $\kappa$ of the train track $\tau_E$ does not intersect $\alpha'\cup\gamma'$ before 
 a wave move with respect to $\gamma'$. Then similar to the argument above, the surgery-step of the wave 
 move ``pushes'' an arc out of the cusp and ``splits'' $\gamma'$ into two curves $\gamma_1$ and $\gamma_2$.  
 One of the resulting curves, say $\gamma_1$, meets this $y$-arc $\kappa$ in exactly one point.  Note that 
 $\gamma_1$ cannot be the curve parallel to $\alpha'$ because this $y$-arc $\kappa$ does not meet $\alpha'$ 
 and two parallel curves must have the same sequence of intersection points with $\delta$ and $\varepsilon'$ 
 (since there is no bigon intersection by Remark~\ref{rem:twosides} and since each arc of 
 $(\alpha'\cup\gamma')\cap T$ intersects $\tau_E$ in at most one point). So $\gamma_2$ is the curve 
 parallel to $\alpha'$ and it  is removed in the second step of the wave move.  Thus, if a $y$-arc $\kappa$ 
 of the train track $\tau_E$ does not intersect $\alpha'\cup\gamma'$ before a wave move, then $\kappa$ 
 intersects the new pair of curves $\alpha'\cup\gamma'$ after the wave move in exactly one point and 
 hence is a pure $y$-arc.  

Therefore, if a $y$-arc $\kappa$ is a pure $y$-arc before a wave move, it remains a pure 
$y$-arc after the wave move.  
\end{proof}

Continue with the sequence of wave moves and as before, if the $x$-arc of the train track, after some 
wave moves, no longer intersects $\alpha'\cup\gamma'$, we split the train track as in Figure~\ref{fig:splitting} 
and get a new $x$-arc and a new pure $y$-arc.  Since the splitting in Figure~\ref{fig:splitting} is along an 
$x$-arc which does not intersect $\alpha'\cup\gamma'$, when performing the operations of wave moves 
and splittings, the following two properties are preserved:

\begin{enumerate}
	\item The curves $\alpha'\cup\gamma'$ always intersect at least two segments of the train 
	track $\tau_E$.
	
	\item Each component of  $(\alpha'\cup\gamma')\cap T$ intersects the traintrack $\tau_E$ in 
	at most one point.
\end{enumerate}
 
Since the weight of $\varepsilon'$ at each segment of the original train track $\tau_E$ is at least two, 
 the splitting arc of $\NN(\tau_E)$, see Definition~\ref{def:splitting-arc}, must ``pass through" each 
 segment. This means that, if we proceed with all the wave moves and splittings, every segment of the 
 train track  will eventually become an $x$-arc.  Since the splitting in Figure~\ref{fig:splitting}
  turns the $x$-arc into a pure $y$-arc, after a number of such operations, we reach a configuration where
both $y$-arcs are pure $y$-arcs and $\varepsilon'$ has weight one at each $y$-arc. This implies that the
 weight of $\varepsilon'$ at the $x$-arc is 2.  

Let $y_0$ and $y_1$ be the two $y$-arcs, such that $y_0$ is the pure $y$-arc right after the splitting in 
Figure~\ref{fig:splitting} and $y_1$ is the other $y$-arc which is a pure $y$-arc  \emph{from a previous splitting}.  
So $y_0$ does not intersect $\alpha'\cup\gamma'$.  Since $\alpha'\cup\gamma'$ always intersects at least two 
segments of the train track, both $y_1$ and the $x$-arc intersect $\alpha'\cup\gamma'$.  As $y_1$ is a pure 
$y$-arc, all the points of $y_1\cap(\alpha'\cup\gamma')$ belong to the same curve. Without loss of generality, 
suppose $y_1\cap(\alpha'\cup\gamma')\subset\alpha'$.  

As discussed earlier, since $\alpha'\cup\gamma'$ intersects the $x$-arc at this stage, there is an 
$[\varepsilon'^-,\varepsilon'^+]$ edge. By the hypothesis of this case, this means that 
we have not exhausted  the sequence of wave moves $\eta_1,\dots,\eta_k$.  

Suppose the next wave move is $\eta_m$ ($m\le k$).  Since $y_0$ does not intersect 
$\alpha'\cup\gamma'$ and by the argument above, after the wave move $\eta_m$, the $y$-arc 
$y_0$ intersects $\alpha'\cup\gamma'$ in exactly one point.  As we have assumed 
$y_1\cap(\alpha'\cup\gamma')\subset\alpha'$,  after the wave move $\eta_m$, 
then $y_0\cap(\alpha'\cup\gamma')$ is a point in $\gamma'$.

If the $x$-arc no longer intersects $\alpha'\cup\gamma'$ after the wave move $\eta_m$, then since 
$y_1$ only intersects $\alpha'$ and since the weight of $\varepsilon'$ at $y_0$ is one, 
$\gamma'\cap\varepsilon'=\gamma'\cap y_0$ is a single point after the wave move.  As $\alpha'$ and 
$\gamma'$ are boundary curves of disks in the handlebody $V$ and $\varepsilon'$ is the boundary of a disk in $W$, 
this means that  $\varepsilon'$ and $\gamma'$ are boundary curves of a destabilizing pair of disks for 
$\Sigma$ in $M$ (when viewing $\Sigma$ as a Heegaard surface of $M$).  By Claim~\ref{cl:PD}, 
when viewing $\Sigma$ as a Heegaard surface of the lens space by the Dehn surgery on $K$,  the curves 
$\varepsilon'$ and $\partial_+P_{\varepsilon'}$ are boundary curves of a destabilizing pair of disks for $\Sigma$ 
in the lens space. So after compressing $\Sigma$ along the disk $E'$ in the handlebody $W$ bounded by 
$\varepsilon'$, the resulting torus is a Heegaard torus for both $M$ and the lens space obtained by 
the Dehn surgery.  This means that $K$ is contained in a  Heegaard solid torus in $M$ which remains a solid torus 
after a nontrivial Dehn surgery. Knots  in solid tori  with this property are called Berge-Gabi knots.
They were classified by Berge \cite{Be2} and Gabai \cite{Gabai} who also showed that they are all doubly primitive.  

Thus, we may assume that, after the $\eta_m$ wave move, the $x$-arc still intersects $\alpha'\cup\gamma'$. 
At this point the situation is that: both $y_0$ and $y_1$ are pure $y$-arcs and all three segments of the train 
track $\tau_E$ intersect $\alpha'\cup\gamma'$.  

If the $x$-arc 
intersects $\alpha'\cup\gamma'$ in more than one point, then as the surgery-step of the wave move ``pushes" 
only one arc ``out of'' the $x$-arc,  after the surgery-step, the $x$-arc still intersects the resulting curves. 
Recall that each component of $(\alpha'\cup\gamma')\cap T$ intersects the train track $\tau_E$ in at most one point.   
This implies that the second step of the wave move (which deletes a parallel copy of $\alpha'$ or $\gamma'$) cannot 
remove all the remaining intersection points of this $x$-arc with $\alpha'\cup\gamma'$ since the two parallel 
curves have the same sequence of intersection points with the train track $\tau_E$.  Thus, if the $x$-arc intersects 
$\alpha'\cup\gamma'$ in more than one point, then this $x$-arc still intersects $\alpha'\cup\gamma'$ after one 
additional wave move.

Therefore, we can continue the sequence of wave moves until the $x$-arc intersects $\alpha'\cup\gamma'$ 
in a single point. As the $x$-arc still intersects $\alpha'\cup\gamma'$,  $\alpha'\cup\gamma'$ still contains 
an $[\varepsilon'^+,\varepsilon'^-]$-edge. The hypothesis of Case (1) says that, after we have exhausted  the 
sequence of wave moves $\eta_1,\dots,\eta_k$, there is no more $[\varepsilon'^+,\varepsilon'^-]$-edge. 
Since we still have an $[\varepsilon'^+,\varepsilon'^-]$-edge at this stage, this means that we have not exhausted 
the sequence of wave moves. Each wave in the sequence of wave moves connects a $[\delta^+,\delta^-]$-edge 
to an $[\varepsilon'^+,\varepsilon'^-]$-edge, as explained before Remark~\ref{rem:twosides}. Thus the curves 
$\alpha'\cup\gamma'$ must still contain a $[\delta^+,\delta^-]$-edge.

As in the earlier discussion, we have  
$\emptyset\ne y_1\cap(\alpha'\cup\gamma')\subset\alpha'$ and 
$\emptyset\ne y_0\cap(\alpha'\cup\gamma')\subset\gamma'$. 
Without loss of generality, suppose the  intersection point of the $x$-arc with $\alpha'\cup\gamma'$ 
belongs to $\gamma'$. So there is a subarc of $\varepsilon'$ connecting the $\gamma'$-arc that 
intersects the $x$-arc to a $\gamma'$-arc that intersect $y_0$.  As orientation of $\alpha'\cup\gamma'$ 
is compatible along the train track, this subarc of $\varepsilon'$ is a $[\gamma'^+,\gamma'^-]$ edge.  

We now have two situations to consider. If $|y_1\cap(\alpha'\cup\gamma')|\ge 2$, then a subarc of $y_1$ 
between two points of $y_1\cap(\alpha'\cup\gamma')$ corresponds to a subarc of $\varepsilon'$ which is an 
$[\alpha'^+,\alpha'^-]$ edge. By part (1) of Proposition~\ref{lem:+ and -}, there is no wave with 
respect to $\{\alpha',\gamma'\}$.  Moreover, we have concluded above that $\alpha'\cup\gamma'$ 
contains both  $[\delta^+,\delta^-]$ and $[\varepsilon'^+,\varepsilon'^-]$ blocking 
edges at this stage. So by part (1) of Proposition~\ref{lem:+ and -}, there is no wave with respect to 
$\{\delta,\varepsilon'\}$ either. This contradicts Theorems~\ref{thm:HOT} and \ref{thm:NeOk}.

If $|y_1\cap(\alpha'\cup\gamma')|= 1$, then  since the weight of $\varepsilon'$ at $y_1$ is one, 
$\alpha'$ intersects $\varepsilon'$ in just one point. Hence $\varepsilon'$ and $\alpha'$ are 
boundary curves of disks which are a destabilizing pair for $\Sigma$ in $M$ (when viewing $\Sigma$ 
as a Heegaard surface of $M$).  As before,  $\varepsilon'$ and $\partial_+P_{\varepsilon'}$ are 
boundary curves of a destabilizing pair of disks for $\Sigma$ in the lens space that was obtained by the 
Dehn surgery. As in the argument above, this implies that $K$ is a Berge-Gabai knot in a Heegaard 
solid torus and hence is  doubly primitive, by \cite{Be2} and \cite{Gabai}. This finishes the proof of  Case (1).

\vskip 7pt

\noindent\underline{\bf Case 2}. $\Gamma(\delta,\varepsilon')$ contains no more $[\delta^-,\delta^+]$-edges 
after the last wave move $\eta_k$.
\vskip 7pt

The difference between Case 2 and Case 1 is due to the difference in the configuration of $\tau_D$.  
Recall that the train track $\tau_D$ consists of a circle in $\widehat{\Sigma}$ and an arc $\rho_x$ 
that contains the point $X=\delta\cap\partial_+P$.  By construction, the weight of $\delta$ at the 
segment $\rho_x$ is one, and this is the main difference between $\tau_D$ and $\tau_E$.  
If we simply apply the argument in Case 1 to $\tau_D$, then we may not be able to conclude 
that both $y$-arcs of the train track will become pure $y$-arcs as in Case 1. The way to deal 
with Case 2 is to use the rectangle-annulus decomposition of $\widehat{\Sigma}$ 
and the proof below is similar to the proof in Section~\ref{sec:OneArTwoAl}.

By Remark~\ref{rem:twosides}, the subarcs of $\alpha'\cup\gamma'$ next to a cusp of $\tau_D$ 
are an $\alpha'$-arc and a $\gamma'$-arc. As before, $\tau_D$ is invariant under $\pi$. Similar to the 
discussion on $\tau_E$ in Case 1, for any segment of $\tau_D$, the outermost (on the segment) intersection 
points  of $\alpha'\cup\gamma'$ with the segment of $\tau_D$ belong to the same curve $\alpha'$ or $\gamma'$.  
As in Case (1), this implies that $\alpha'\cup\gamma'$ intersects at least two segments of the train track $\tau_D$.  

Similar to Case 1, continue with the wave moves and, at each stage, denote the resulting curves by 
$\alpha'$ and $\gamma'$. We first discuss how the curves $\alpha'$ and $\gamma'$ after the wave
move intersect $\partial_+P$ and $\widehat{\Sigma}$. The goal is to show that, after some isotopy, 
the intersection of $\alpha'$ and $\gamma'$ with $\widehat{\Sigma}$ consists of vertical arcs in 
$\mathcal{A}_l$, $\mathcal{A}_r$, $\mathcal{R}^u$ and $\mathcal{R}^d$. This will show in particular 
that the claims in Section~\ref{sec:OneArTwoAl} also hold in this setting. 

Since every non-separating simple closed curve in $\Sigma$ is invariant under $\pi$, after isotopy, 
we may assume that  $\varepsilon'$, $\delta$ and $\partial_+ P$ are all invariant. Consider the wave 
$\eta$. As illustrated in Figure~\ref{fig:2cusp}(a), $\eta$ connects an arc in the cusp of $\NN(\tau_D)$ 
to an arc in the cusp of $\NN(\tau_E)$. We may suppose $\eta$ intersects $\partial_+P$ minimally. 
Note that the proof of Claim~\ref{cl:junction-arc} only uses the minmality of $\eta\cap\partial_+P$ 
and the symmetry of $\delta$ and $\varepsilon'$ under the involution, which means that 
Claim~\ref{cl:junction-arc} also holds for $\eta$. That is, $\partial_+P$ divides $\eta$ into a pair 
of junction arcs at the ends and possibly some vertical arcs in $\mathcal{A}_l$ and $\mathcal{A}_r$. 
Thus, after isotopy, the new curve resulting from the wave move along $\eta$ intersects $\widehat{\Sigma}$ 
in a collection of vertical arcs in the rectangles $\mathcal{R}^u$, $\mathcal{R}^d$ and annuli 
$\mathcal{A}_l$, $\mathcal{A}_r$. Furthermore, Claim~\ref{cl:NoAWaves} also holds, that 
is, any wave with respect to the new set of meridians $\{\alpha',\gamma'\}$ must intersect $\partial_+P$. 
Therefore, by repeatedly applying Claims~\ref{cl:junction-arc} and \ref{cl:NoAWaves}, we may 
assume $\{\alpha',\gamma'\}$ always intersect $\widehat{\Sigma}$ in vertical arcs after isotopy. 

Before proceeding, we study how these wave moves and isotopies affect the $x$- and $y$-arcs of 
$\tau_D$.  Without loss of generality, assume that the cusps of $\tau_D$, which are symmetric 
under $\pi$, lie in $\mathcal{A}_l$. 

There are four possible configurations of $\delta$ in $\mathcal{A}_l$,  as shown in Figure~\ref{fig:Al}. 
In Figure~\ref{fig:Al}(a, b), the cusp direction of the two cusps of  $\tau_D$ are horizontal in $\mathcal{A}_l$ 
and the shaded regions correspond to the $x$-arc of $\tau_D$. In this case, after a number of wave moves, 
the shaded regions in Figure~\ref{fig:Al}(a, b) no longer contain any $\alpha'$- or $\gamma'$-arcs. Thus, 
the next wave arc must pass through the shaded region  in Figure~\ref{fig:Al}(a, b). 

As illustrated in Figure~\ref{fig:junction-arc}(b, c),  there is an isotopy so that, after the isotopy, the wave 
becomes a union of junction arcs and vertical arcs.  Notice that the effect of the isotopy from 
Figure~\ref{fig:junction-arc}(b) to (c) changes the configuration of $\delta\cap\mathcal{A}_l$ from 
that of Figure~\ref{fig:Al}(a, b) to that of Figure~\ref{fig:Al}(c, d). Therefore, it remains to consider the 
configurations in Figure~\ref{fig:Al}(c, d) for $\delta\cap\mathcal{A}_l$. In this case the two $y$-arcs of 
$\tau_D$ are $\rho_x$ and the short path in $\mathcal{A}_l$. Note that in the configurations of 
Figure~\ref{fig:Al}(c, d),  the cusp directions of the two cusp of $\tau_D$ point into $\mathcal{R}^u$ 
and $\mathcal{R}^d$. 

\begin{figure}[!ht]
	\begin{overpic}[width=7cm]{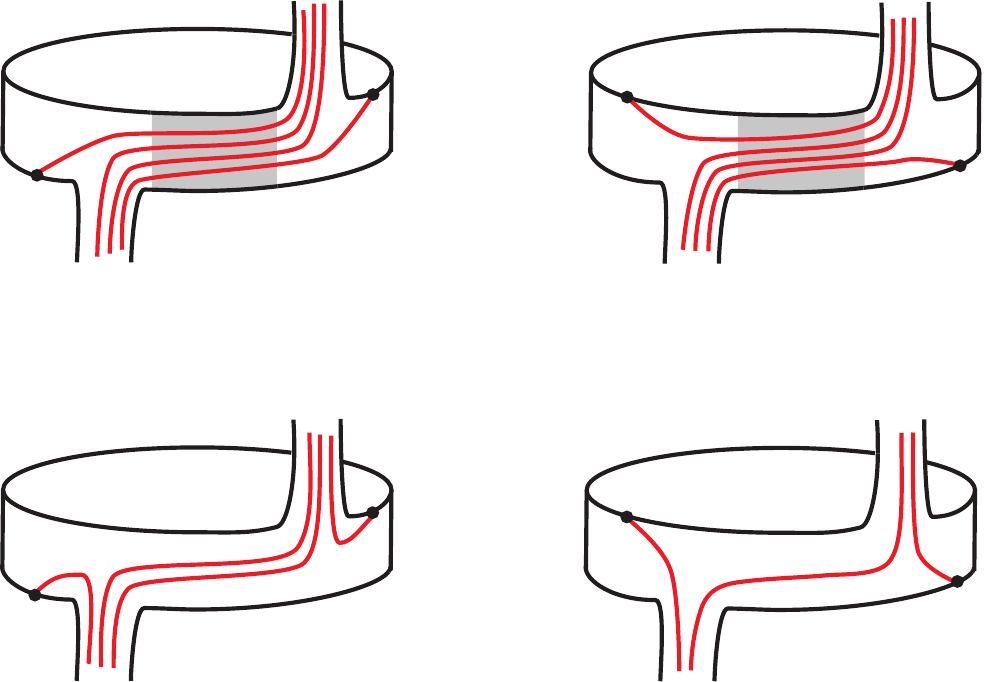}
			\put(17,38){(a)}
			\put(77,38){(b)}
			\put(17,-5){(c)}
			\put(77,-5){(d)}
		\end{overpic}
	\vspace{5pt}
	\caption{Possible configurations of $\delta$ in $\mathcal{A}_l$}
	\label{fig:Al}
\end{figure}

Similar to the discussion in Section~\ref{sec:OneArTwoAl}, it follows from Lemma~\ref{lem:TwoDeltaArcs}, 
that as long as $(\alpha'\cup\gamma')\cap (\mathcal{R}^u\cup\mathcal{R}^d)\ne\emptyset$, the curves 
$\alpha'\cup\gamma'$ contain $[\delta^+,\delta^-]$-edges inside $\mathcal{R}^u$ and $\mathcal{R}^d$. 
Thus, by the hypothesis of Case 2, after a number of wave moves, we reach a situation where 
$(\alpha'\cup\gamma')\cap (\mathcal{R}^u\cup\mathcal{R}^d)=\emptyset$. 

The continuation of the argument is similar in spirit to the argument in the proof of 
Claim~\ref{cl:AlphaGamma}. 

Fix an orientation for $\partial_+P$ and as illustrated in Figure~\ref{fig:orientation}, the two endpoints of 
each vertical arc of  $(\alpha'\cup\gamma')\cap\mathcal{A}_l$ and $(\alpha'\cup\gamma')\cap\mathcal{A}_r$ 
are a pair of intersection points of $(\alpha'\cup\gamma')\cap\partial_+P$ with the same sign.  Since 
$(\alpha'\cup\gamma')\cap (\mathcal{R}^u\cup\mathcal{R}^d)=\emptyset$  at this stage, the intersection 
$(\alpha'\cup\gamma')\cap\widehat{\Sigma}$ consists of vertical arcs in $\mathcal{A}_l$ and $\mathcal{A}_r$. 
This implies that the all the intersection points of $\alpha'\cap\partial_+P$ have the same sign and the intersection 
points of $\gamma'\cap\partial_+P$ all have the same sign.

\begin{figure}[!ht]
\begin{overpic}[width=7cm]{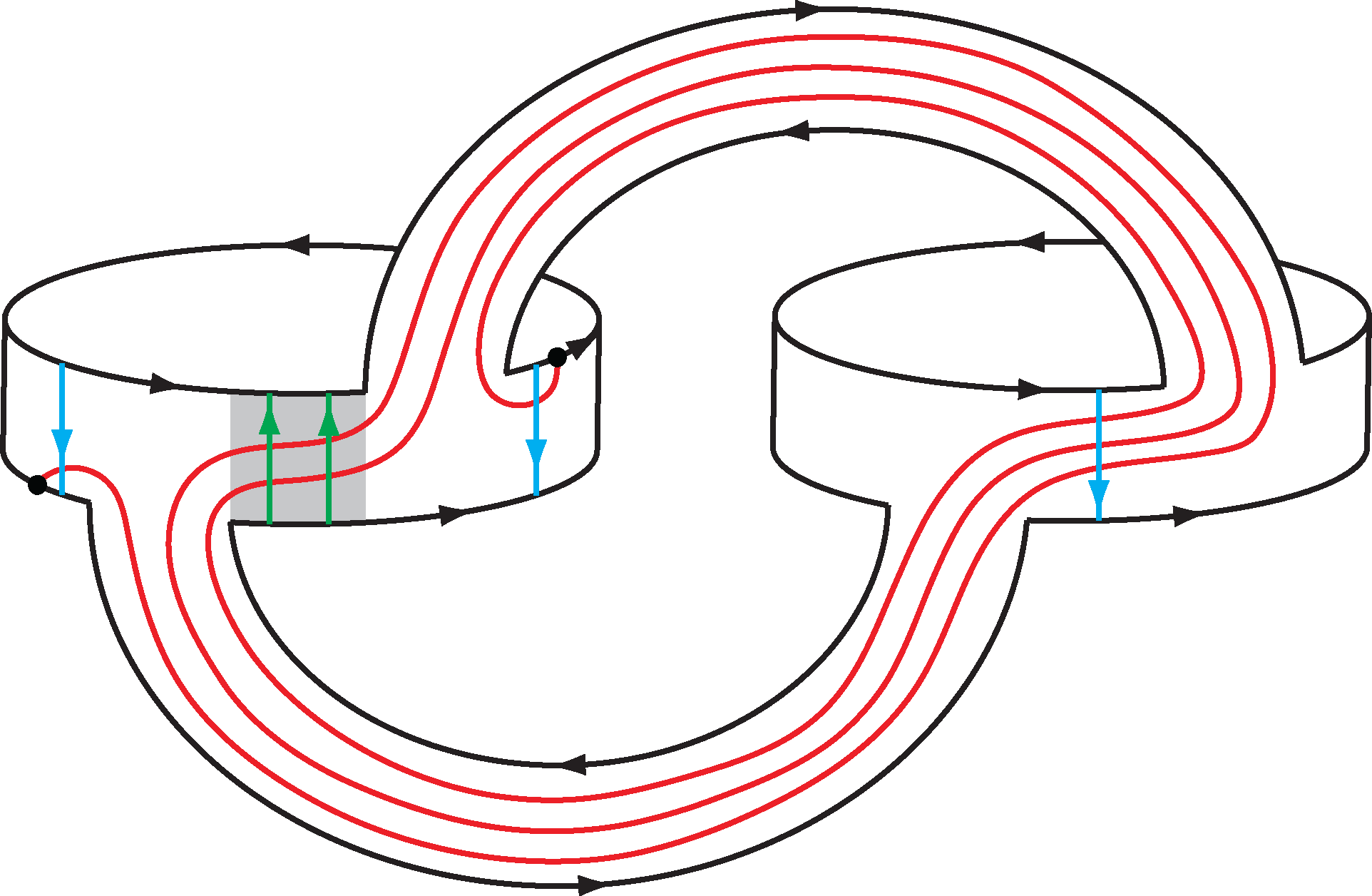}
\put(66,30){$\delta$}
\put(21,22){$\alpha'$}
\put(38,25){$\gamma'$}
\put(101,34){$\mathcal{A}_r$}
\end{overpic}
\caption{Orientations of $\partial_+P$, $\alpha'$ and $\gamma'$}
\label{fig:orientation}
\end{figure}

Let $R_x$ and $R_y$ be the two rectangles of  $\mathcal{A}_l\ssm(\mathcal{J}_l^u\cup\mathcal{J}_l^d)$ 
and suppose $X\subset\partial R_x$.  See the shaded region in Figure~\ref{fig:orientation} for a picture of 
$R_y$.  Note that $\rho_x$ is a $y$-arc of  the train track $\tau_D$ and the arcs of $\delta\cap R_y$ project 
onto the other $y$-arc of $\tau_D$ when we pinch $\delta$ onto $\tau_D$.  Moreover, since both
$\gamma\cap\mathcal{R}^d\ne\emptyset$ and  $\gamma\cap\mathcal{R}^u\ne\emptyset$ before any wave 
move, at least one wave move must be performed to reach this stage.  Similar to the proof of Case 1, 
this implies that $\alpha'\cup\gamma'$ must intersect both $y$-arcs of the train track. In other words, 
$(\alpha'\cup\gamma')\cap R_y\ne\emptyset$ and $(\alpha'\cup\gamma')\cap\rho_x\ne\emptyset$.  

Since the intersection points of $\delta$ with $\alpha'$ and $\gamma'$ all have the same sign, 
the arcs of $(\alpha'\cup\gamma')\cap R_y$ and the arcs of $(\alpha'\cup\gamma')\cap R_x$ which meet 
$\rho_x$ must have opposite directions with respect to the core curve $\mathfrak{a}_l$ of $\mathcal{A}_l$, 
see Figure~\ref{fig:orientation} for a picture.  
However,  the argument above showed that the intersection points of $\alpha'\cap\partial_+P$ all have 
the same sign and the intersection points of $\gamma'\cap\partial_+P$ all have the same sign. 
These conclusions imply that arcs in $(\alpha'\cup\gamma')\cap R_y$ must belong to only one of the curves 
$\alpha'$ or $\gamma'$. Without loss of generality, suppose that all arcs in $(\alpha'\cup\gamma')\cap R_y$  
belong to $\alpha'$. Then the arcs of $\alpha'\cup\gamma'$ that meet $\rho_x$ all belong to $\gamma'$.
  
The conclusion above says that one $y$-arc of $\tau_D$ intersects only $\alpha'$ and the other 
$y$-arc of $\tau_D$  intersects only $\gamma'$.  Therefore, during the sequence of splittings, once 
$(\alpha'\cup\gamma')\cap (\mathcal{R}^u\cup\mathcal{R}^d)=\emptyset$, then both $y$-arcs of 
$\tau_D$ are pure $y$-arcs and the situation becomes similar to Case 1.

Now we can apply the argument for $\tau_E$ and $\varepsilon'$ in Case 1 to $\tau_D$ and $\delta$ 
(i.e., we switch the roles of $\tau_D$ and $\tau_E$). As indicated in Figure~\ref{fig:orientation}, the 
union of the $\delta$-curves in $\mathcal{R}^u$, $\mathcal{R}^d$ and the short path in  $\mathcal{A}_r$ 
are pinched into the $x$-arc of $\tau_D$. By the argument in Case 1, after some wave moves, we reach 
a situation that the curves $\alpha'\cup\gamma'$ intersect the $x$-arc in a single point. Since the 
involution $\pi$ interchanges $\mathcal{R}^u$ and $\mathcal{R}^d$, the symmetry from $\pi$ implies 
that, if $\alpha'\cup\gamma'$ intersects $\mathcal{R}^u$, it must also intersect $\mathcal{R}^d$. 
As $\alpha'\cup\gamma'$ intersects the $x$-arc in a single point, at this stage, we have that 
$(\alpha'\cup\gamma')\cap (\mathcal{R}^u\cup\mathcal{R}^d)=\emptyset$ and $\alpha'\cup\gamma'$ 
intersects the short path in $\mathcal{A}_r$ exactly once, see the blue arc in $\mathcal{A}_r$ in 
Figure~\ref{fig:orientation}. Now the situation becomes exactly the same as  Case 1: 
by performing 
a sequence of wave moves with respect to $\{\alpha',\gamma'\}$ and possibly splitting the train track 
$\tau_D$, we can eventually either conclude that $K$ is a Berge-Gabai knot or obtain a contradiction 
to Theorems~\ref{thm:HOT} and \ref{thm:NeOk}.

\vskip5pt

\noindent \underline{\bf Case 3}. The weight of $\varepsilon'$ at some segment of $\tau_E$ is one, 
and $\Gamma(\delta,\varepsilon')$ contains no more $[\varepsilon'^-,\varepsilon'^+]$-edges after the last 
wave move $\eta_k$. 

The idea of the proof is to convert Case 3 to the setup of Case 2.
By the construction of $\tau_E$, the hypothesis of Case 3 implies that there is a properly embedded arc 
$\tau$ in the once-punctured torus $T=\Sigma\ssm\NN(\tau_D)$ such that $\tau\cap\tau_E=\tau\cap\varepsilon'$ 
is a single point. One can band sum parallel copies of the planar surfaces $P_l$ and $P_r$, as in Claim~\ref{cl:PD}, 
to obtain a planar surface $P_{\varepsilon'}$ such that $\partial_+P_{\varepsilon'}\cap T=\tau$.  
So, $(P_{\varepsilon'},E')$ is a $(\mathcal{P},\mathcal{D})$-pair, where $E'$ is the disk in $W$ 
bounded by $\varepsilon'$.  In particular, $\partial_+P_{\varepsilon'}\cap\tau_E=\partial_+P_{\varepsilon'}\cap\varepsilon'$
 is a single point since $\tau\cap\tau_E=\tau\cap\varepsilon'$.

By Claim~\ref{cl:tauE}, we may assume that
$$|(\alpha\cup\gamma)\cap\tau_E|\le c_0(P,D,\alpha,\gamma).$$ 
Since $\partial_+P_{\varepsilon'}$ intersects both $\varepsilon'$ and $\tau_E$ in a single point, it follows from 
the definition of the complexity that 
$$c_0(P_{\varepsilon'},E',\alpha,\gamma)\le |(\alpha\cup\gamma)\cap\tau_E|.$$
Thus $$c_0(P_{\varepsilon'},E',\alpha,\gamma)\le c_0(P,D,\alpha,\gamma).$$  
Since $c_0(P,D,\alpha,\gamma)$ is minimal among all $(\mathcal{P},\mathcal{D})$-pairs and all such 
curves $\alpha$, the equality holds. In particular, 
$$c_0(P_{\varepsilon'},E',\alpha,\gamma)= |(\alpha\cup\gamma)\cap\tau_E|.$$  
We may assume that $P_{\varepsilon'}$ is minimal in the sense that $|\partial_+P_{\varepsilon'}\cap(\alpha\cup\gamma)|$ 
is the smallest among all such planar surfaces $P_{\varepsilon'}$ with $\partial_+P_{\varepsilon'}$ intersecting 
both $\tau_E$ and $\varepsilon'$ in a single point.

Let $\widehat{\Sigma}_\varepsilon'$ be the closure (under the path metric) of $\Sigma\ssm\partial_+P_{\varepsilon'}$. 
Similar to the decomposition of $\widehat{\Sigma}$ into rectangles and annuli, we have a similar decomposition of 
$\widehat{\Sigma}_\varepsilon'$ into a pair of annuli joined by a pair of rectangles.  Moreover, each component of 
$(\alpha\cup\gamma)\cap\widehat{\Sigma}_\varepsilon'$ is a cocore arc of an annulus or a rectangle in this 
decomposition.  Since $c_0(P_{\varepsilon'},E',\alpha,\gamma)= |(\alpha\cup\gamma)\cap\tau_E|$, each component 
of $(\alpha\cup\gamma)\cap\widehat{\Sigma}_\varepsilon'$ intersects $\tau_E$ in at most one point.  Similar to the 
discussion on $\widehat{\Sigma}$ in Section~\ref{sec:NoCircularArcs}, a cocore arc of a rectangle in the 
decomposition is a boundary arc of a $\partial$-compressing disk for $P_{\varepsilon'}$, and if one performs a 
$\partial$-compression on $P_{\varepsilon'}$, one obtains two planar surfaces whose $\partial_+$-boundaries 
are the core curves of the two annuli in the decomposition.  This, plus the assumption that 
$|\partial_+P_{\varepsilon'}\cap(\alpha\cup\gamma)|$ 
is the smallest among all such planar surfaces, implies that Lemma~\ref{lem:TwoDeltaArcs} 
is also true for $\varepsilon'$ and this decomposition of $\widehat{\Sigma}_\varepsilon$. In other words, 
$\varepsilon'$ passes through each rectangle in the decomposition of $\widehat{\Sigma}_\varepsilon'$ 
at least twice. In particular, this implies that the curves $\alpha\cup\gamma$ must contain an 
$[\varepsilon'^+,\varepsilon'^-]$-edge as long as they intersect the two rectangles in the rectangle-annulus 
decomposition of $\widehat{\Sigma}_\varepsilon'$.

Since $\tau_E$ only has two cusps and since each component of 
$(\alpha\cup\gamma)\cap\widehat{\Sigma}_\varepsilon'$ intersects $\tau_E$ in at most one point, 
the structure of $\tau_E$ in $\widehat{\Sigma}_\varepsilon'$ is the same as the structure of $\tau_D$ 
in $\widehat{\Sigma}$. In particular, $\varepsilon'$ must take one short path in each annulus 
in this decomposition. Since $\varepsilon'$ passes through each rectangle in the decomposition 
of $\widehat{\Sigma}_\varepsilon'$ at least twice (similar to Lemma~\ref{lem:TwoDeltaArcs}), 
we have all the ingredients and the argument in Case 2 also works for $\varepsilon'$, 
$\partial_+P_{\varepsilon'}$ and $\widehat{\Sigma}_\varepsilon'$ in this setup.  Thus  
Case 3 follows from the argument in Case 2 after switching the roles of $\delta$, $\tau_D$ and 
$\varepsilon'$, $\tau_E$ respectively and using the rectangle-annulus structure of $\widehat{\Sigma}_\varepsilon'$ 
instead of $\widehat{\Sigma}$. So, after a sequence of wave moves with respect to $\{\alpha',\gamma'\}$,  
we can eventually either conclude that $K$ is doubly primitive or obtain a contradiction to Theorems~\ref{thm:HOT} 
and \ref{thm:NeOk}. 
\end{proof}

\section{The Proof}\label{sec:TheProof} 

\vskip10pt

In this section we prove the main  Theorem ~\ref{thm:MainTheorem}.

\begin{proof}[Proof of Theorem ~\ref{thm:MainTheorem}]  We begin with the meridional systems 
$\widehat{W}=\{\delta,\varepsilon\}$ and $\widehat{V}=\{\alpha,\gamma\}$, as chosen in Section
 \ref{sec:MeridionalDisksandWhiteheadGraphs} for the Heegaard diagram for $M$.  As described 
 at the beginning of Chapter~\ref{cpt:ObtainingTheContradiction}, we consider all possible 
configurations of $\delta$ depending on the paths that $\delta$ takes in the two annuli. 

Now Propositions~\ref{pro:OneRTwoL}, \ref{pro:AllPaths}, \ref{pro:LongPath}, and \ref{pro:OneAndOne} deal with 
all possible configurations of $\delta$. They show that either $K$ is doubly primitive or that the 
corresponding Heegaard diagram is not induced by a Heegaard splitting of $M$. This proves  
Theorem  ~\ref{thm:MainTheorem}.
\end{proof}


\end {document}